\documentclass[12pt,a4paper,final]{article}

\usepackage[margin=24.5mm]{geometry}
\usepackage[utf8]{inputenc}  
\usepackage{xcolor,scrtime,graphicx}
\usepackage{amsmath,amsfonts,amssymb}
\usepackage{lmodern} 

\numberwithin{figure}{section}
\usepackage{todonotes} 

\usepackage{hyperref}
\usepackage[notref,notcite]{showkeys}

\usepackage{tikz}
\usepackage{pgfplots}

\usepackage{textgreek,bm}
\usepackage{bef_alex}
\def\FG{\bm}
\def\mafo{\mathrm}

\numberwithin{equation}{section}

\newcommand{\HH}{h}
\newcommand{\HHhh}{H}

\newcommand{\ul}{\underline}

\newcommand{\sign}{\mathop{\mafo{sign}}}
\newcommand{\ppp}{\kappa}
\newcommand{\inveps}{\tfrac{\ds1}{\ds \eps}}
\newcommand{\mfD}{\mathfrak D}
\newcommand{\mfJ}{\mathfrak J}
\newcommand{\mfp}{\mathfrak p}

\newcommand{\weak}{\rightharpoonup}
\newcommand{\weakstrong}{\xrightarrow[\!\!{}^{\rmw{\ti}\rms}\!\!]{}}
\newcommand{\Gweak}{\overset{\Gamma}\weak}
\newcommand{\Gweakstrong}{\overset{\Gamma}{\weakstrong}}
\newcommand{\Eto}{\overset{\rmE}\to}
\newcommand{\Gto}{\overset{\Gamma}\to}
\newcommand{\Mto}{\overset{\rmM\rmo}\longrightarrow}
\newcommand{\EDPto}{\overset{\mafo{EDP}}{\longrightarrow}}
\newcommand{\relEDPto}{\xrightarrow{\mafo{relEDP}}}

\newcommand{\eff}{\mathrm{eff}} 
\newcommand{\prim}{\mathrm{prim}} 
\newcommand{\dual}{\mathrm{dual}} 
\newcommand{\wigg}{\mathrm{wig}}

\begin{document}
 
\title{ A gradient system with a wiggly energy \\
and relaxed EDP-convergence\thanks{Research partially supported by DFG
  via SFB\,1114 (project C05)}}

\author{Patrick Dondl\thanks{Universit\"at Freiburg}, Thomas
  Frenzel\thanks{WIAS Berlin}, Alexander Mielke\thanks{WIAS Berlin and
  Humboldt-Universit\"at zu Berlin}}

\date{8.\ Januar 2018} 

\maketitle
%

\section{Introduction}
\label{s:Intro}

 This paper is devoted to the general question of convergence of a
family of gradient systems $(\bfQ,\calE_\eps,\calR_\eps)$  towards an
effective gradient system $(\bfQ,\calE_0,\calR_\eff)$ when the small
parameter $\eps\to 0$. Here $\bfQ$ is the state space (e.g.\ a convex
subset of a Banach space), $\calE_\eps:[0,T]\ti \bfQ\to \R$ are the
possibly time-dependent energy functionals, and $\calR_\eps$ are the
dissipation potentials such that the gradient-flow equation reads
\[
0 = \rmD_{\dot q}\calR_\eps(q_\eps,\dot q_\eps) + \rmD_q
\calE_\eps(t,q_\eps).
\]  
The objective is to show that limits $q_0$ of solutions $q_\eps$ are
solutions of the limiting gradient system $(\bfQ,\calE_0,\calR_\eff)$,
where typically $\calE_0$ is the \textGamma-limit of the energies
$\calE_\eps$, but in some interesting cases the effective dissipation
potential $\calR_\eff$ in the limiting equation
\begin{equation}
  \label{eq:LimEqn}
  0 = \rmD_{\dot q}\calR_\eff(q_0,\dot q_0) + \rmD_q
\calE_0(t,q_0)
\end{equation} 
differs from the \textGamma-limit $\calR_0$ of the dissipation
potentials $\calR_\eps$. However, we are not so much interested in the
effective equation, but in the limiting gradient structure
$(\bfQ,\calE_0,\calR_\eff)$ that
contains additional information to the limiting equation
\eqref{eq:LimEqn}. Indeed, in \eqref{eq:GS.exa} we give four different
gradient structures for the simple ODE $\dot q =1-q$. 

A general study of \textGamma-convergence for gradient systems was
initiated in \cite{SanSer04GCGF}, which lead to a rich body of
research, see \cite{Stef08BEPD, Serf11GCGF, Brai13LMVE, Visi13VFSS,
  Miel16EGCG} and the references therein.  Several convergence notions
are covered by the general name \emph{evolutionary
  \textGamma-convergence}, which emphasizes that evolutionary problems
are treated by variational methods involving \textGamma-convergence
for the associated functionals. In this work, we want to generalize
the notion of evolutionary \textGamma-convergence in the sense of the
\emph{energy-dissipation principle} (in short EDP-convergence) 
introduced in \cite{LMPR17MOGG}, which is the first notion that
provides a method to calculate the effective dissipation potential
$\calR_\eff$ in a unique way. 
 
Our new notion of \emph{relaxed EDP-convergence} for gradient systems
is explained by studying in detail the following wiggly-energy
model
\begin{equation}
\label{eq:evol.main}
\nu\dot{u}=-\rmD\calE_\eps(t,u), \quad u(0)=u^0\in \R, 
\end{equation}
with the energy 
\[
\calE_\eps(t,u) = \Phi(u)-\ell(t)u+\eps \ppp(u,\inveps u),
\]
where $\ppp(u,\cdot)$ is a 1-periodic function,  and the
dissipation potential is simply $\calR(\dot u)=\frac\nu2 \dot u{}^2$. This model was
introduced in \cite{Jame96HPT, AbChJa96KMWE} as
a very simple model for explaining slip-stick motions in martensitic
phase transformations by starting from a linear viscosity law as in
\eqref{eq:evol.main}. See also \cite{Meno02GSWE, Sull09AGDR} for
vector-valued versions (i.e.\ $u(t)\in \R^n$) of such gradient systems. 
Earlier models for explaining dry friction
go back to Prandtl \cite{Pran28GKTF} and Tomlinson \cite{Toml29MTF},  
see also \cite{PopGra12PTMH} for historical remarks.
The general feature of such models is that a viscous evolution law
in a temporally constant, but spatially rapidly varying energetic
environment may lead to stick-slip motion, where the limit evolution
cannot be a described by the homogenized energy alone. In particular,
we find that the effective dissipation potential $\calR_\eff$ is
much bigger than $\calR_0=\calR$, where the difference depends on the
wiggly part $\ppp$ of the the energy landscape.   

Further applications of such models occur in the evolution of phase
boundaries in a heterogeneous environment is modeled in
\cite{Bhattacharya_99}, based on \cite{AbKn_88}, or in the evolution
of dislocations in a slip plane with heterogeneities like forest
dislocations~\cite{Garroni:2005ve, Garroni:2006tn, Monneau:2010te,
  Dondl:2017ut} (when neglecting lattice friction).  Applications to
crawling are studied in \cite{GidDes17GDFB}, and an extension to creep
is given in \cite{SKTO09BDSC}.

A different approach to modeling phase transforming materials by
considering connected bistable springs also leads to a complex energy
landscape and an evolution in effective wiggly potential
\cite{PugTru02MTP,PugTru05TRIP}.  A rigorous derivation of
rate-independent one-dimensional pseudo-elasticity is given in
\cite{MieTru12DVEC}. The latter papers as well as
\cite{PugTru02RIHB,Miel12ERID} are especially devoted to the
mathematical justification of the rate-independent case, where
$\nu_\eps \to 0$ as $\eps \to 0$, such that the limit dynamics doesn't
have any internal time-scale any more. 

Here we revisit the general class of scalar wiggly-energy models in the form
\begin{equation}
  \label{eq:I.GenModel}
  \pl_{\dot u}\calR(u,\dot u) = - \rmD_u\calE_\eps(t,u), \quad
  u(0)=u^0\in \R, 
\end{equation}
where $\calR:\R^2 \to {[0,\infty[}$ is a fixed dissipation potential,
i.e.\ $\calR(u,0)=0$ and $\calR(u,\cdot)$ is convex, while the energy
$\calE_\eps$ is as above. Thus, \eqref{eq:I.GenModel} is the flow
induced by the gradient system $(\R,\calE_\eps,\calR)$. Under suitable
assumptions it is well known from the above works (see e.g.\
\cite{AbChJa96KMWE, Meno02GSWE, PugTru02RIHB, Sull09AGDR}) that the
solutions $u_\eps$ of \eqref{eq:I.GenModel} converge for $\eps \to 0$
to limits $u_0$ that are solutions of the limiting gradient system
$(\R,\calE_0,\calR_\eff )$. We emphasize that $\calE_\eps$ converges
uniformly to the limit energy $\calE_0:(t,u)\mapsto \Phi(u)-\ell(t)u$,
however, the restoring forces $\rmD\calE_\eps$ do not converge because
of the wiggly part involving the non-decaying, oscillatory term
$\pl_y\ppp(u,\inveps u)$, where $y$ is used as a placeholder for the
second argument $\inveps u \in \bbS^1:= \R/_{\Z}$ of $\ppp$. The
major task is then to find the effective dissipation potential
$\calR_\eff $, which, as we will see, is larger than $\calR$ and
depends on $\pl_y\ppp$.

The purpose of this work is to show how the gradient structure of the
underlying problem can be exploited in a natural way using the method
for evolutionary \textGamma-convergence for gradient systems. 
Thus, we (i) obtain the effective dissipation potential
$\calR_\eff  $ (and as a by-product the limit evolution) by
purely energetic principles, (ii) identify a new mechanical function
$(\dot u, \xi) \mapsto \calM (u,\dot{u},\xi)$, which we call
\emph{contact potential}, that encodes the effective
dissipation law, but which is not a dual pairing in the form
$\calR_\eff (u,\dot u){+}\calR^*_\eff (u,\xi)$, and finally (iii) 
discuss the convexity properties of $\calM (u,\cdot,\cdot)$ in the
sense of bipotentials, see \cite{BuDeVa08ECBG,BuDeVa08?NMCM}.

To be more specific, we use the formulation of gradient flows
via the following energy-dissipation principle,
which originates in the work of De Giorgi \cite{DeMaTo80PEMS} and
states that \eqref{eq:I.GenModel} is equivalent to the energy
dissipation balance (EDB) stated below. The
EDB asks simply that the final
energy plus the dissipated energy equals the initial energy plus the
work of the external forces, where the dissipated energy has to be
expressed in a particular way in terms of $\calR$ and its
Legendre-Fenchel dual
$\calR^*$,  namely 
\begin{equation}
\label{eq:I.EDB}
  \calE_\eps(T,u(T)) + \mfD_\eps(u) =
  \calE_\eps(0,u(0))
 + \int_0^T\partial_t\calE_\eps(t,u(t))\dd t, 
\end{equation}
where the dissipation functional  $\mfD_\eps$ is given by 
\begin{equation}
  \label{eq:mfD.R.E}
  \mfD_\eps(u)=\int_0^T \!\! \Big(\calR(u(t),\dot{u}(t)) + 
\calR^{*}\big(u(t),-\rmD\calE_\eps(t,u(t))\big)\Big) \dd t.
\end{equation}
Several notions of \emph{evolutionary \textGamma-convergence} rely on
passing to the limit $\eps \to 0$ in \eqref{eq:I.EDB} (cf.\
\cite{Miel16EGCG}) and identifying
the limits of the four terms accordingly, see Section \ref{s:Model}. 

In our case the convergence of $u_\eps(t)\to u(t)$ immediately
implies, for all $t\in [0,T]$, the convergence
$\calE_\eps(t,u_\eps(t)) \to \calE_0(t,u(t))$ as well as $\pl_t
\calE_\eps(t, u_\eps(t)) \to \pl_t \calE_0(t,u(t))$.  Thus, it remains
to understand the limit of $\calD_\eps(u_\eps)$, and the notion of
EDP-convergence asks for the identification of the \textGamma-limit of
$\mfD_\eps$ on a suitable subset of functions $u \in \rmW^{1,p}(0,T)$
with $p\in {]1,\infty[}$. Our main technical results are in Section
\ref{se:Homog} and imply the desired statement
\[
\mfD_\eps \Gweak \mfD_0 \quad \text{with } 
\mfD_0(u)=\int_0^T \calM (u,\dot u, {-}\rmD\calE_0(t,u)) \dd t,
\]
The novelty of the notion of EDP-convergence is that we study
$\calD_\eps$ not only along the exact solutions $u_\eps$ of
\eqref{eq:I.GenModel} (or equivalently \eqref{eq:I.EDB}), but rather along
general functions. This reflects the fact that a given evolution equation
$\dot u =F(t,u)$ may have different gradient structures, and this
difference is only seen by looking at fluctuations around the
deterministic solutions, cf.\ \cite{PeReVa14LDSH, MiPeRe14RGFL, 
  LMPR17MOGG}. These fluctuations explore $\mfD_\eps$ also away from
the exact solutions of the gradient flow. 

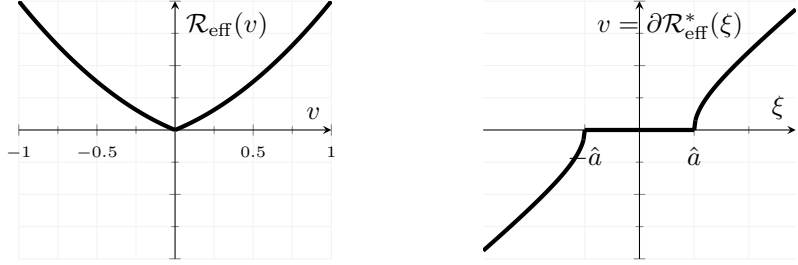
\begin{figure}
\centering
\begin{tikzpicture}
  \begin{axis}%
  [     
      grid=both,
      minor tick num=1,
      grid style={line width=.1pt, draw=gray!10},
      xlabel={\footnotesize $v$},
      axis x line=middle,
      xmin=-1,
      xmax=1,
      tick label style={font=\tiny},
      ylabel={{\footnotesize $\calR_\eff(v)$}},
      axis y line=middle,
      ymin=-1,
      ymax=1,
      yticklabels={ , , },
      no markers,
      samples=200,
      domain=-1:1,
      scale=0.6
    ]
    \addplot[ultra thick, domain=-1 : 0.0] (x,{ -0.5*(1-x)*x});
    \addplot[ultra thick, domain=0.0 : 1] (x, { 0.5*(1+x)*x });
  \end{axis}
\end{tikzpicture}
\qquad\qquad
\begin{tikzpicture}
  \begin{axis}%
  [     
      grid=both,
      minor tick num=1,
      grid style={line width=.1pt, draw=gray!10},
      xlabel={\footnotesize$\xi$},
      axis x line=middle,
      xmin=-1,
      xmax=1,
      xtick = {-0.35,0.35},
      xticklabels={$-\hat a$ ,$\hat a$ },
      tick label style={font=\footnotesize},
ylabel={{\footnotesize \hspace*{-2em}$v=\partial\calR_\eff^{*}(\xi)$}},
      axis y line=middle,
      ymin=-1,
      ymax=1,
      yticklabels={ , , },
      no markers,
      samples=200,
      domain=-1:1,
      scale=0.6
    ]
    \addplot[ultra thick,domain=-0.35:0.35] (x,{0});
    \draw[ultra thick] (axis cs: 0.35, 0) -- (axis cs: 0.358, 0.0753);
    \draw[ultra thick] (axis cs: -0.35, 0) -- (axis cs: -0.358, -0.0753);
    \addplot[ultra thick,domain=-1:-0.358] (x,{-(x^2-0.35^2)^0.5});
    \addplot[ultra thick,domain=0.358:1] (x,{(x^2-0.35^2)^0.5});
  \end{axis}
\end{tikzpicture}
\caption{The dissipation potential $\calR_\eff$ and the kinetic
  relation $v=\pl\calR^{*}_\eff (\xi)$ for the quadratic case, see
  \eqref{eq:I.Case}.} 
\label{fig:KinRel}
\end{figure}

Theorem \ref{th:GCvg.Jwigg} provides the explicit form of the effective
contact potential $\calM $, viz.\ 
\begin{equation}
  \label{eq:I.rmM0}
  \calM (u,v,\xi) := \inf\Bigset{ \int_0^1
    \!\!\Big(\calR(u,|v|\dot{z}(s))+ \calR^*\big(u, 
  \xi - \pl_y \ppp(u,z(s))\big)\!\Big) \rmd s}{ z\in\rmW_v^p(0,1)},
\end{equation}
where $\rmW^{1,p}_v :=
\set{z\in\rmW^{1,p}(0,1)}{z(1)-z(0)=\sign(v)}$.  The proof is 
a generalization of the homogenization results in
\cite{Brai02GCB} for functionals of the form $u \mapsto \int_0^T
f(t,u,\inveps u) \dd t$: 


In Section \ref{se:M0.prop} we discuss the basic properties of
$\calM $, which allows us to recover the limiting evolution and to
identify the effective dissipation potential $\calR_\eff $. In
fact, we show 
\begin{align*}
\text{(i) \ } &\calM (u,v,\xi)\geq \xi v,\\
\text{(ii) } & \calM (u,v,\xi)=\xi v \ \Longleftrightarrow \xi \in
\pl_v \calR_\eff (u,v) 
\end{align*}
for a unique effective dissipation potential $\calR_\eff$. Thus, all
ingredients of relaxed EDP-convergence (cf.\ Definition
\ref{def:relEDP}) are established.   
The main observation here is that the contact sets 
\[
\mathsf C_{\calM}(u):= \bigset{(v,\xi)\in \R^2}{
  \calM (u,v,\xi)=\xi v }
\]
can be identified directly giving a general formula for
$\calR_\eff  $ in terms of a harmonic mean of $y \mapsto \pl_\xi
\calR^*(u,\xi{-}\pl_y\ppp(u,y))$, see Lemma \ref{le:calM0}. 
Of course, we recover the classical result of
\cite{Jame96HPT,AbChJa96KMWE} for the case $\calR(u,v)=\frac{1}{2\mu}
v^2$ and $\ppp(u,y) = \hat{a}\sin(2\pi y)/(2\pi)$, namely 
\begin{equation}
  \label{eq:I.Case}
\calR_\eff (v)=\int_0^{|v|}\big(\hat a^2{+}\frac{\hat v^2}{\mu^2}\big)^{1/2}\rmd
\hat v \ \Longleftrightarrow \ 
\pl\calR^{*}_\eff (\xi) = \mu\,\sign(\xi) 
  \big(\max\{  \xi^2{-}\hat{a}^2, 0\}\big)^{1/2}. 
\end{equation}
See also Figure \ref{fig:KinRel} for $\calR_\eff$ and the kinetic relation
$v=\pl\calR^*_\eff(\xi)$. We note that for a non-degenerate wiggly
potential this leads to a 
motion of the interface that is large compared to the excess driving
force $\xi- \hat a$ near the depinning transition. This is in agreement with
experiments, where it is seen that a phase boundary propagates nearly
freely when subjected to a driving force above the critical value
\cite{EsCl_93, AbKn_97}.

 Hence, $\mathsf C_\calM (u) $ is the graph of a
subdifferential of $\calR_\eff(u,\cdot)$ which determines $\calR_\eff$
uniquely, which in the sense of \cite{Visi13VFSS} can be understood as
$\calM (u,\cdot,\cdot)$ representing the monotone operator $v\mapsto
\pl_v \calR_\eff(u, \cdot)$. However, there the function
$\calM (u,\cdot,\cdot)$ is assumed to be convex, which is not the
case in our model.

Of course, $\calM $ contains more information than $\calR_\eff$, 
and it is worth to study $\calM $ as such, as we expect it to be
relevant as rate function for suitable large deviation limits in the
sense of \cite{BonPel16QRIL}.  In Section \ref{su:Bipot} we discuss the question
whether $\calM $ is a bipotential in the sense of
\cite{BuDeVa08ECBG,BuDeVa08?NMCM}, which means that 
\begin{subequations}
\label{eq:I.Bipot}
\begin{align}
\label{eq:I.Bipot.A}
\text{(i)\ }& \calM (v,\cdot, \xi) \text{ and } \calM (u,v,\cdot)
\text{ are convex}, \\
\text{(ii)}& v \in \pl_\xi\calM (u,v,\xi) \ \Longleftrightarrow \ 
\calM (u,v,\xi)=\xi v 
\ \Longleftrightarrow \ \xi \in \pl_v\calM (u,v,\xi) . 
\label{eq:I.Bipot.B}
\end{align}
\end{subequations}
While $\calM (u,\cdot,\xi)$ is always convex, our Example
\ref{ex:M.noncvx}  shows that in general $\calM (u,v,\cdot)$ is
non-convex. For the special $p$-homogeneous case
$\calR(u,v)=r(u)|v|^p$ we are able to show that $\calM $ is indeed a
bipotential, see Theorem \ref{th:M.cvx.xi}.  

In Section \ref{se:Discuss} we discuss the results and
highlight specific properties of this limit procedure and compare
it with recent results in \cite{Visi13VFSS,Visi15?SSFE,Visi17?EGCW}
concerning related evolutionary \textGamma-convergence results based
on an extended version of the Brezis-Ekeland-Nayroles principle, see
Section \ref{su:EGCWeakT}.  We explicitly show that
$\calM(u,v,\xi)\neq \calR_\eff(u,v){+}\calR_\eff^*(u,\xi)$, which
implies that there is no EDP-convergence in the sense of
\cite{LMPR17MOGG}.

Moreover, for converging
solutions $u_\eps(t)\to u_0(t)$ of \eqref{eq:I.EDB} we easily obtain
$\mfD_\eps (u_\eps) \to \mfD_0(u_0)$, i.e.\ solutions are recovery
sequences for the dissipation functional. However, if we separate the
dissipation into its primal and its dual part, the corresponding convergences  
\begin{align*}
&\mfD^\prim_\eps(u_\eps):=\int_0^T \calR(u_\eps,\dot u_\eps)\dd t
&&\to\  \mfD^\prim_\eff(u_0):=\int_0^T
\calR_\eff(u_0,\dot u_0)\dd t \quad \text{and} \\
&\mfD^\dual_\eps(u_\eps):=\int_0^T
\calR^*(u_\eps,{-}\rmD\calE_\eps(t,u_\eps))\dd t\hspace{-1.2em} &&\to \ 
\mfD^\dual_\eff(u_0):=\int_0^T \calR^*_\eff(u,{-}\rmD\calE_0(t, u) )\dd t
\end{align*}
do not hold. Indeed, for quadratic $\calR:v\mapsto \frac\nu2v^2$ we
always have 
\[
\mfD^\prim_\eps(u_\eps)= \mfD^\dual_\eps(u_\eps)=\frac12
\mfD_\eps(u_\eps) \ \to \ \frac12 \mfD_0(u_0),
\]
but $\calR_\eff$ is such that $ \mfD^\prim_\eff(u_0)\gneqq
\mfD^\dual_\eff(u_0)$ if $\dot u_0\not\equiv 0$. This shows that the
classical approach of \cite{SanSer04GCGF} is not applicable because of
an exchange of dissipation between the dual part $\mfD^\dual$ and the
primal part $\mfD^\prim$ in the limit $\eps \to 0$. This is again
reflected in the fact that $\calR_\eff$ is larger that $\calR$ and
depends on $\pl_y \ppp$.

\section{Evolutionary  \textGamma-convergence and main results}
\label{s:Model}

\subsection{The energy-dissipation principle for gradient system}
\label{su:EDP.GS}

To explain the general structure between our special model of
\eqref{eq:I.GenModel} we use general ordinary differential equation (ODE)
$\dot q = \bfF(t,q)\in \R^n$ and general gradient systems (GS)
$(\bfQ,\bfE,\bfR)$, where $\bfQ=\R^n$ is the state space,
$\bfE:[0,T]\ti\bfQ\to \R$ is a sufficiently smooth, time-dependent
energy functional, and $\bfR:\bfQ\ti \bfQ\to {[0,\infty[}$ is a
sufficiently smooth dissipation potential. By $\bfR^*$ we denote the
(Legendre-Fenchel) dual dissipation potential defined via 
$\bfR^*(q,\bfxi)=\sup\set{\langle \bfxi,\bfv\rangle-\bfR(q,\bfv)}{\bfv \in
  \bfQ}$.  
  
We say that the ODE $\dot q = \bfF(t,q)$ \emph{has a gradient
  structure} or is a \emph{gradient flow} if there exists a GS
$(\bfQ,\bfE,\bfR)$ such that 
$\bfF(t,q)=\pl_\bfxi \bfR^*(q,{-}\rmD_q \bfE(t,q))$. In that case, we
also say that the ODE is a \emph{generated by the GS}
$(\bfQ,\bfE,\bfR)$. We emphasize that one ODE can have several
distinct gradient structures, e.g.\ $\dot q = 1- q\in \R$ is generated
by the gradient systems $({[0,\infty[},\bfE_j,\bfR_j)$ for $j=1,\ldots,4$ with  
with 
\begin{align}
\label{eq:GS.exa}
&\bfE_1(q)=\bfE_2(q)=\frac12(1{-}q)^2, \quad \bfR^*_1(\xi)=\frac12 \xi^2,\quad 
\bfR^*_2(q,\xi)=\frac{\frac12\xi^2+\frac14\xi^4}{1+(1{-}q)^2},\\ 
\nonumber
&\bfE_3(q)=\bfE_4(q)= q \log q - q +1, \
\bfR^*_3(q,\xi)=\frac{q{-}1}{2\log q}\xi^2,
\quad \bfR_4^*(q,v)=
2\sqrt q \,\big(\!\cosh(\tfrac12 \xi)-1\big).
\end{align}
We also refer to \cite{PeReVa14LDSH, MiPeRe14RGFL} for discussion of
different gradient structures for the heat equation or for
finite-state Markov processes. Thus, we emphasize that the gradient
structure of a given ODE has additional physical information, e.g.\
about the microscopic origin of the ODE, see \cite{LMPR17MOGG}. 
This is seen in the above case, since we may choose different energies
$\bfE_j$ and even for one chosen $\bfE_j$ we may choose different
dissipation functionals $\bfR_k$.

We recall that the evolution law associated with a gradient system can
be written in two equivalent ways, namely 
\begin{equation}
  \label{eq:ODE.GS}
  0 \in \pl_{\dot q}\calR(q,\dot q) + \rmD_q\calE(t,q) \quad
\Longleftrightarrow \quad \dot q \in \pl_{\xi}
\calR^*(q,{-}\rmD_q\calE(t,q)).
\end{equation}
The energy-dissipation principle states that under reasonable
technical assumptions these relations are equivalent to a scalar
energy-dissipation balance. To motivate this we 
consider a lower semi-continuous convex function $\Psi:X\to \R_\infty$
on a reflexive Banach space $X$. Denote by $\Psi^*:X^*\to \R_\infty$
the Legendre-Fenchel dual, i.e.\ $\Psi(\xi)=\sup\set{\langle
  \xi,v\rangle -\Psi(v)}{v\in 
  X}$. Then, the Fenchel equivalences (see \cite{Fenc49CCF,
EkeTem76CAVP} or \cite[Thm 23.5]{Rock70CA}) state that
\[
\text{(i) } \xi \in \pl \Psi(v) \quad \Longleftrightarrow \quad
\text{(ii) } v \in \pl \Psi^*(\xi) \quad \Longleftrightarrow \quad
\text{(iii) } \Psi(v)+\Psi^*(\xi) = \langle \xi,v\rangle,
\] 
 where $\pl$ denotes the convex subdifferential. 
Indeed, by the definition of $\Psi^*$ we have the Fenchel-Young
inequality $\Psi(v)+\Psi^*(\xi)\geq \langle \xi,v\rangle $ for all
$v\in X$ and $\xi\in X^*$. Thus, in (iii) it would suffices to ask for
the inequality $\Psi(v)+\Psi^*(\xi)\leq \langle \xi,v\rangle $.
 
Applying this with $\Psi=\calR(q,\cdot)$, integration over time and
using the chain rule we see that $q$ solves \eqref{eq:ODE.GS} if
and only if $q$ satisfies the  \emph{energy-dissipation balance} 
\begin{equation}
  \label{eq:EDP1}
\begin{aligned} 
&\calE(T,q(T)) + \mfD(q) =  \calE(0,q(0))
- \int_0^T \rmD_t \calE(t,q(t)) \dd t,\\
&\text{where } \mfD(q):= \int_0^T \Big(\calR(q,\dot q) +
\calR^*\big(q,-\rmD_q\calE(t,q)\big)\Big) \dd t. 
\end{aligned}
\end{equation}
Indeed, using the chain rule $\frac{\rmd}{\rmd t} \calE(t,q(t)) =
\rmD_t\calE(t,q(t)) + \langle \rmD_q\calE(t,q(t)), \dot q(t)\rangle$
(the validity of which is the main technical assumption in the general
infinite-dimensional case) it is easy to go back from \eqref{eq:EDP1}
to \eqref{eq:ODE.GS}, as we deduce 
\[
\int_0^T   \Big(\calR(q,\dot q) +
\calR^*\big(q,-\rmD_q\calE(t,q)\big) - \langle\rmD_q\calE(t,q(t)),
\dot q(t)\rangle  \Big) \dd t = 0.
\]
As the integrand in non-negative by the Fenchel-Young inequality and
the integral is $0$, we conclude that the integrand is $0$ almost
everywhere, which means (iii) in the Fenchel equivalences. Thus (i)
and (ii) also hold almost everywhere, i.e.\ \eqref{eq:ODE.GS} holds. 
We refer to \cite{AmGiSa05GFMS,Miel16EGCG} for more details and exact
statements.

\subsection{Evolutionary \textGamma-convergence for gradient systems}
\label{su:EGC.GS}
We now consider families $(Q,\calE_\eps,\calR_\eps)$ of gradient
systems depending on a small parameter $\eps>0$. We are interested in
the limits $u_0$ of solutions as well as in suitable limiting gradient
systems $(Q,\calE_0,\calR_0)$.  
 
Hence, for $\eps\in [0,\eps_0]$ we consider the gradient-flow
equations 
\begin{equation}
  \label{eq:GF.eps}
  0=\pl_{\dot u} \calR_\eps(q_\eps, \dot q_\eps) + \rmD_q
  \calE_\eps(t,q_\eps), \quad q_\eps(0)=q^0_\eps. 
\end{equation}
Following \cite{Miel16EGCG} we say that the family $(Q,\calE_\eps,\calR_\eps)$ of gradient systems
\emph{E-converges} the gradient system $(Q,\calE_0,\calR_0)$, and
shortly write $(Q,\calE_\eps,\calR_\eps) \overset{\rmE}\to
(Q,\calE_0,\calR_0)$, if the following holds: If $q_\eps^0 \to q_0^0$
and $q_\eps:[0,T]\to Q$ are solutions of \eqref{eq:GF.eps} for
$\eps\in {]0,\eps_0[}$, then there
exist a subsequence $0<\eps_k \to 0$ and a solution $q_0:[0,T]\to Q$
for \eqref{eq:GF.eps} with $\eps=0$ such that 
\begin{equation}
  \label{eq:Eto}
  \forall\, t\in {]0,T]}:\quad q_{\eps_k}(t)\to q_0(t) \text{ and }
  \calE_{\eps_k} (t,q_{\eps_k}(t)) \to \calE_0(t,q_0(t)). 
\end{equation}
(A similar notion $\overset{\rmE}\weak$ can be defined by replacing
strong with weak convergence.) Note that the selection of subsequences
is only needed if the limiting underlying gradient systems does not have
uniqueness of solutions. In that case different subsequences may
converge for to different solutions of \eqref{eq:GF.eps}$_{\eps=0}$
with the same initial condition $q_0^0$.

A major drawback of this notion is that $\calR_0$ is not intrinsically
connected to the original gradient systems
$(Q,\calE_\eps,\calR_\eps)$. Indeed, if $(Q,\calE_0,\calR_0)$ and
$(Q,\calE_0,\wh\calR_0)$ generate the same gradient-flow equation
(i.e. $\pl_\xi \calR_0^*(q,{-}\rmD_q\calE_0(t,q))= \pl_\xi
\wh\calR_0^*(q,{-}\rmD_q\calE_0(t,q))$, see \eqref{eq:GS.exa} for
examples) and if $(Q,\calE_\eps,\calR_\eps) \overset{\rmE}\to
(Q,\calE_0,\calR_0)$, then we also have $(Q,\calE_\eps,\calR_\eps)
\overset{\rmE}\to (Q,\calE_0,\wh\calR_0)$. The notion of \emph{EDP
  convergence} is stricter and involves the effective dissipation
potential $\calR_\eps$  for $\eps\in {[0,\eps_0[}$ directly through
the dissipation functionals $\mfD_\eps$ defined via
\begin{equation}
\label{eq:mfD.eps}
\mfD_\eps(q(\cdot)):= \int_0^T \Big( \calR_\eps(q,\dot
q)+\calR_\eps^*\big(q,{-}\rmD_q\calE_\eps (t, q)\big) \Big) \dd t.  
\end{equation}
The following definition now asks \textGamma-convergence of
$\mfD_\eps$ to $\mfD_0$, and thus $\calR_\eps$ are intrinsically
involved. The new feature is that we ask much more than convergence of
these functionals along solutions $q_\eps$ converging to $q_0$. In
light of \cite{LMPR17MOGG} this seems to be essential, since the gradient
structures contain more information than the equations determining the
solutions. We refer to the discussion in Section \ref{se:Discuss}.

\begin{definition}[EDP-convergence, cf.\ \cite{LMPR17MOGG}] \label{def:EDPcvg}
\mbox{} \ The gradient systems
$(Q,\calE_\eps,\calR_\eps)_{]0,\eps_0]}$ are said to 
\emph{converge to the gradient system} $(Q,\calE_0,\calR_0)$ \emph{in the
sense of the energy-dissipation principle}, shortly 
``\emph{EDP-converge}'' or $(Q,\calE_\eps,\calR_\eps)_{]0,\eps_0[} 
\EDPto(Q,\calE_0,\calR_0)$, if the following conditions hold:
\begin{subequations}
\label{eq:CondEDPcvg}
\begin{align}
\label{eq:CondEDPcvg.a}
&(Q,\calE_\eps,\calR_\eps)  \Eto (Q,\calE_0,\calR_0),\\
\label{eq:CondEDPcvg.b} &\calE_\eps \Gto \calE_0, \quad \text{ and
} \quad \mfD_\eps \Gweak \mfD_0,  
\end{align}
\end{subequations}
where specific choice of the \textGamma-convergence 
$\Gweak$ in \eqref{eq:CondEDPcvg.b} needs to be specified in each
particular case. 
\end{definition}

Two remarks are in order. First, as we highlight in Section
\ref{se:Discuss}, the EDP-convergence does in general not imply that
the two contributions of the dissipation function (generated by
$\calR_\eps$ and $\calR_\eps^*$, respectively) converge
individually. Indeed, this may even be wrong when restricting to
solutions.

Second, it is one of the main results of this paper that the structure
of $\mfD_\eps$ may not be preserved by taking the
\textGamma-limit in general. Under suitable technical assumptions the 
techniques in \cite{Dalm93IGC} show that a \textGamma-limit $\mfD_0$ 
has the integral form $\mfD_0(q)=\int_0^T \calN_0(t,q,\dot q) \dd t$,
but $\calN_0$ may not have the form 
\[
\calN_0(t,q,\dot q)= \calR_0(q,\dot q) +
\calR^*_0(q,{-}\rmD_q\calE_0(t,q))
\]
for any $\calR_0$. 

In our wiggly-energy model as well as in many other applications 
we  have a time-dependent external loading $\ell:[0,T]\to Q^*$,
and we want to have a result that works uniformly in with respect to
$\ell$. Thus, we
look at driven gradient systems with
\[
\calE_\eps(t,q)= \calF_\eps(q) - \langle \ell(t),q\rangle \quad
\text{and} \quad \calF_\eps \Gto \calF_0. 
\]
Because of $\rmD_q\calE_\eps(t,q)=\rmD \calF_\eps(q) - \ell(t)$ and
the arbitrariness of $\ell$, we introduce the variable $\xi \in Q^*$
as a placeholder of variants for the restoring force
$-\rmD_q\calE_\eps$. Indeed, we use the decomposition
\begin{equation}
  \label{eq:Decomp.DE}
  -\rmD_q \calE_\eps(t,q) = \Xi_\eps(q)+\ell(t) - \Omega_\eps(q),
\end{equation}
where $\Xi_\eps$ is supposed to converge nicely to the desired limit
$\rmD \calF_0(q)$, while $\Omega_\eps(a) $ somehow 
converges to $0$. Thus, we can write $\mfD_\eps$ in the form 
\begin{align}
\nonumber
&\mfD_\eps(q) = \mfJ_\eps\big( q,-\rmD_q\calE_\eps(t,q){+}\Omega_\eps(q)\big), 
\quad \text{where }\\
&
\label{eq:def.mfJeps}
\mfJ_\eps(q,\xi)= \int_0^T\Big(\calR_\eps(q,\dot q)+ \calR_\eps^* \big( q,
\xi{-}\Omega_\eps(q)\big)\Big) \dd t.
\end{align}
As is observed in \cite{Visi13VFSS} it is important that $\dot q$ and
$\xi$ are in duality and that the convergences of $\dot q_\eps$ to
$\dot q_0$ and of $\xi_\eps$ to $\xi_0$ are such that the duality
pairing $(\dot q, \xi) \mapsto \int_0^T \langle \xi(t),\dot
q(t)\rangle \dd t$ is continuous. In most applications one uses 
\begin{equation}
  \label{eq:Cvg.q.xi}
  q_\eps \weak  q_0 \text{ in } \rmW^{1,p}(0,T;Q) \text{ (weakly)}\quad 
\text{and} \quad \xi_\eps \to \xi_0 \text{ in } \rmL^{p'}(0,T;Q^*) 
\text{ (strongly)}.
\end{equation} 
This explains why the decomposition \eqref{eq:Decomp.DE} is useful: we
obtain the strong convergence $\Xi_\eps(q_\eps(\cdot)) \to
\Xi_0(q_0(\cdot))$ and want to use $\Omega_\eps(q(\cdot)) \weak 0$ in a
suitable sense. 

Now, we may consider \textGamma-convergence for the functionals
$\mfJ_\eps$ with respect to the convergence in \eqref{eq:Cvg.q.xi},
denoted by ``$\weakstrong$''. Again, under suitable assumption the
theory in \cite{Dalm93IGC} predicts that a possible \textGamma-limit
takes the following form
\begin{equation}
  \label{eq:mfJ.cvg}
\mfJ_\eps \Gweakstrong \mfJ_0: (q,\xi)\mapsto \int_0^T \calM (q,\dot
q, \xi)\dd t, 
\end{equation}
where now $\calM (q,\cdot,\cdot):Q\ti Q^* \to [0,\infty]$ contains
the effective information on the dissipation for a given macroscopic
speed $v=\dot q\in Q$ and an effective macroscopic force $\xi \in
Q^*$. Even in the case $\Omega_\eps \equiv 0$ we see that the
convergence $\weakstrong$ from \eqref{eq:Cvg.q.xi} is the natural one
for studying the \textGamma-limit of $\mfJ_\eps$, since under suitable
coercivity assumptions one has 
\[
\calR_\eps(q, \cdot) \Gweak \calR_0(q,\cdot) \text{ in }Q \quad \Longleftrightarrow
\quad \calR_\eps^*(q,\cdot) \Gto \calR_0^*(q,\cdot) \text{ in } Q^*,
\]  
see \cite[p.\,271]{Atto84VCFO} and the survey
\cite[Sec.\,3.2]{Miel16EGCG}. 

As a remainder of the Young-Fenchel inequality
$\calR_\eps(q,v)+\calR_\eps^*(q,\xi) \geq \langle \xi,v\rangle$ one
can hope for the estimate 
\begin{equation}
  \label{eq:genYouFen}
  \forall \, q,v \in Q, \ \xi \in Q^*: \quad \calM (q,v,\xi)\geq
\langle \xi,v\rangle,
\end{equation}
however this has to be proved in each case using properties of
$\Omega_\eps$, see our Lemma \ref{le:calM0}(b) for the wiggly-energy
model. Then, the essential as in the energy-dissipation principle of
the previous subsection the limit evolution is given by
\begin{align*}
&\calM (q,\dot q, {-}\rmD_q \calE_0(t,q)) = -\langle
\rmD_q\calE_0(t,q), \dot q\rangle  \text{ \ or equivalently}\\
&\calE_0(T,q(t))+ \int_0^T\!\!\calM (q,\dot q, {-}\rmD_q \calE_0(t,q))\dd t
= \calE_0(0,q(0))-\int_0^T\langle \dot\ell(t),q\rangle \dd t, 
\end{align*}
where we assumed that $\calE_0(t,q)=\calF_0(q)-\langle \ell(t),q\rangle$
still satisfies a chain rule. 
While $\calM $ encodes information on the combined limit of
$(\calE_\eps,\calR_\eps)$, we can now go back looking at solutions
which necessarily stay in the so-called contact set
$\mathsf{C}_{\calM }(\cdot)$, namely   
\[
(\dot q(t),{-}\rmD_q \calE_0(t,q(t)) \in \mathsf C_{\calM }(q(t)) 
\text{ with } 
\mathsf C_\calM(q):= \bigset{(v,\xi)\in Q\ti Q^*}
                 { \calM(q,v,\xi)=\langle \xi,v \rangle }. 
\] 
This subset gives the admissible pairs $(v,\xi)$ of rates and forces
at a given state $q$, i.e.\ it defines a kinetic relation. 

Our definition of relaxed EDP-convergence now asks that this kinetic
relation is given in terms of a dissipation potential $\calR_\eff$. We
emphasize that using this approach the dissipation $\calR_\eff$ is
uniquely defined through the steps above, i.e.\ as in EDP-convergence
we find ``the'' effective dissipation potential, however in contrast
to EDP-convergence we are more flexible in term of the
\textGamma-limit $\mfD_0$ of $\mfD_\eps$, which may not have
$(\calR_0,\calR_0^*)$ form. That is also the reason why we use the
notation $\calR_\eff$, as there is no direct convergence of
$\calR_\eps$ to $\calR_\eff$, see the discussion in Section
\ref{se:Discuss}.

\begin{definition}[Relaxed EDP-convergence] \label{def:relEDP}
  We say that the family 
 $(Q,\calE_\eps,\calR_\eps)_{]0,\eps_0[}$ of gradient systems  
\emph{converges to the gradient system} $(Q,\calE_0,\calR_\eff)$
\emph{in the relaxed EDP sense}, and
shortly write $(Q,\calE_\eps,\calR_\eps)_{]0,\eps_0[} 
\relEDPto (Q,\calE_0,\calR_\eff)$, if the
following holds. 
\begin{subequations}
\label{eq:relEDP} 
\begin{align}
\label{eq:relEDP.a}&(Q,\calE_\eps,\calR_\eps)_{]0,\eps_0[} 
\Eto (Q,\calE_0,\calR_\eff), \\ 
\label{eq:relEDP.b}&\calE_\eps(t,q) =\calF_\eps(q)-\langle
\ell(t),q\rangle,\quad  \calF_\eps \Gto \calF_0,\\ 
\label{eq:relEDP.c}&\exists\, \Omega_\eps: \ \wt q_\eps\weak \wt q_0\text{ in
}\rmW^{1,p}(0,T;Q) \Longrightarrow \rmD_q\calF_\eps(\cdot, \wt q_\eps ) {-}
\Omega_\eps(\wt q_\eps)  \to \rmD_q \calF_0(\wt q_0),\\  
\label{eq:relEDP.d}& \mfJ_\eps \text{ defined in \eqref{eq:def.mfJeps}
  satisfies \eqref{eq:mfJ.cvg} with } \calM  \text{ satisfying
  \eqref{eq:genYouFen}} ,  \\
\label{eq:relEDP.e}&\exists\,\text{diss.\,pot.\,} \calR_\eff\ 
   \forall\, q\in Q: \
\mathsf C_{\calM }(q)=\bigset{(v,\xi)\in Q\ti Q^*}{ \xi \in
  \pl_v\calR_\eff(q,v)}.  
\end{align}
\end{subequations}
\end{definition}

The aim of this paper is to show that the theory sketched above can be
made rigorous for the wiggly-energy model. Thus, we have a first
non-trivial example that shows that relaxed EDP-convergence provides a
mechanically relevant concept for deriving effective gradient
structures where neither the Sandier-Serfaty theory
\cite{SanSer04GCGF} nor the EDP-convergence from \cite{LMPR17MOGG}
applies.

\subsection{Our model as gradient system and relaxed EDP-convergence}
\label{su:ModelGS}

For our wiggly-energy model, the gradient system is given by the state
space $\R$, the energy $\calE^\wigg_\eps:\R\ti[0,T]\to \R$ and a general
convex dissipation potential $\calR:\R\times\R\to\R$.  We choose
the following 
assumptions to keep the technicalities to a limit; however, it is easily
possible to generalize most assumptions except for the additive
structure of $\calE_\eps$ concerning the wiggly part $\ppp$. 
\begin{subequations}
\label{eq:AssumER}
 \begin{align}
 \label{eq:Assum.E1}
 &\calE^\wigg_\eps(t,u)=\Phi(u) - \ell(t)u +
   \eps \ppp(u,\inveps  u) \text{ with }\Phi \in \rmC^1(\R), \ \ell\in
   \rmC^1([0,T])\\
  \label{eq:Assum.E2}
 & \text{ and } \ppp \in \rmC^1(\R^2) \text{ with } \ppp(u,y{+}1)=\ppp(u,y) 
    \text{ for all }u,y\in \R;\\
\label{eq:Assum.R1}
  &\calR\in \rmC^1(\R^2), \quad \calR(u,v)\geq 0, \quad
   \calR(u,0)=0; \\
\label{eq:Assum.R2}
  &\forall\, u \in \R: \ \calR(u,\cdot ) \text{ is strictly convex}; \\
\nonumber
  & \exists\, p \in {]1,\infty[}\  \exists\, c_1,c_2>0\ \exists\
  \text{modulus of continuity }\omega \ \: \forall\,
  u,\, \wh u,\, v\in \R: \\
 \label{eq:Assum.R3}   
  &\qquad   c_1(|v|^p{-}1) \leq \calR(u,v) \leq c_2(1{+}|v|^p) \text{ and}\\
\label{eq:Assum.R4}
  &\qquad  |\calR(u,v) {-} \calR(\wh u,v)| \leq \omega(|u{-}\wh u|)
   (1{+}|v|^p). 
\end{align}
\end{subequations}
Assumption \eqref{eq:Assum.R3} implies that the dual dissipation potential
$\calR^*$ satisfies the estimate 
\begin{equation}
  \label{eq:p'growth}
  \forall\, u,\xi \in \R:\quad c_3 (|\xi|^{p'}{-}1) \leq \calR^*(u,\xi)
\leq c_4 (1{+}|\xi|^{p'}),
\end{equation}
where $p'=p/(p{-}1)$. Moreover, $\calR^*(u,\cdot)$ is continuously
differentiable and strictly convex. 
The last assumption \eqref{eq:Assum.R4} is a
uniform continuity statement that 
should be avoidable; however, it helps us settle some
technical issues which would otherwise destroy the chosen and
hopefully clear
\textGamma-convergence proof. Again, by using the Legendre-Fenchel
transform we find the corresponding uniform continuity statement for
$\calR^*$, namely 
\begin{equation}
  \label{eq:R*unif.cont}
  \forall\,
  u,\, \wh u,\, \xi\in \R:\quad  |\calR^*(u,\xi) {-} \calR^*(\wh
  u,\xi)| \leq C_p \omega(|u{-}\wh u|) 
   (1{+}|\xi|^{p'}), 
\end{equation}
where $C_p>1$ is a constant depending only on $p>1$.  

As a special case we consider power-law
potentials $\calR(u,v)= \frac{\nu(u)}p |v|^p$ giving
$\calR^*(u,\xi)=\frac{\mu(u)}{p'}|\xi|^{p'}$, where
$\mu(u)=\nu(u)^{1/(1-p)}$. So, the assumptions
\eqref{eq:Assum.R1}--\eqref{eq:Assum.R4} are satisfied if $\nu$ and
$1/\nu$ are positive, bounded and continuous.  


The gradient-flow equation has the usual form 
\begin{equation}
  \label{eq:ODE1} \pl_{\dot u}\calR(u,\dot u)=-
\rmD_u\calE^\wigg_\eps(t,u)= - 
\Phi'(u)+\ell(t) - \eps \pl_u \ppp(u,\inveps u) - \pl_y \ppp(u,\inveps
u ),
\end{equation}
where the wiggly part $\ppp:\R\ti \bbS^1 \to \R$ inserts the small inherent
length scale $\eps$ into the system via the periodicity variable 
$y=u/\eps$.

Following the abstract approach of Sections \ref{su:EDP.GS} and
\ref{su:EGC.GS} equation \eqref{eq:ODE1} is equivalent to the
energy-dissipation balance 
\begin{subequations}
\label{eq:EDB.wiggly}
\begin{align}
  \label{eq:EDB.wig.a}
 & \calE^\wigg_\eps(T,u(T)) +
  \mfJ_\eps^\wigg\big(u,{-}\rmD_q\calE^\wigg_\eps(\cdot,u){+}\Omega_\eps(u)\big)
  = \calE^\wigg_\eps(0,u(0)) - \int_0^T \!\dot\ell u \dd t,\\
\label{eq:EDB.wig.b}&
 \text{with } \Omega_\eps(u):=\pl_y\ppp(u,\inveps u) \text{ and}\\
\label{eq:EDB.wig.c}&
\mfJ_\eps^\wigg(u,\xi):= \int_0^T\!\Big( \calR(u(t),\dot u(t)) +
\calR^*\big(u(t),\xi(t) -\Omega_\eps(u(t))\big) \Big) \dd t.
\end{align}
\end{subequations}

The proof of relaxed EDP-convergence relies on the following 
technical result for the \textGamma-convergence of
$\mfJ_\eps^\wigg$. For this we define the limit dissipation functional  
\begin{align}
&\nonumber
 \mfJ_0^\wigg: \rmW^{1,p}(0,T)\ti \rmL^{p'}(0,T)\to [0,\infty] \text{
   via}
\\ 
&  \label{eq:mfJ0}
\mfJ_0^\wigg(u,\xi):=\int_0^T \calM (u,\dot u,\xi)\dd t \quad
\text{with }\\
  \label{eq:calM0}
  &\calM ^\wigg(u,v,\xi):=\inf_{z\in \rmW^{1,p}_v} \Big( \int_0^1\big[
  \calR(u,|v|\dot{z}(s))+     \calR^*\big(u,
    \xi - \pl_y \ppp(u,z(s))\big)\big] \dd s \Big) ,\\
 \nonumber &\text{where } \rmW^{1,p}_v=\bigset{ v\in \rmW^{1,p}(0,1)}
  {z(1)=z(0){+}\sign(v)}. 
\end{align}
Recalling the definition of weak-strong convergence
\eqref{eq:Cvg.q.xi} in $\rmW^{1,p}(0,T)\ti \rmL^{p'}(0,T)$, which
is denoted by $\weakstrong$, the following result holds.

\begin{theorem}[\textGamma-convergence of $\mfJ_\eps^\wigg$]  
\label{th:GCvg.Jwigg} If the gradient system
  $(\R,\calE^\wigg_\eps,\calR_\eps)$ satisfies the assumptions
  \eqref{eq:AssumER}, then $\mfJ_\eps^\wigg \Gweakstrong \mfJ_0^\wigg$.
\end{theorem}

As a first consequence we obtain a \textGamma-convergence result for
the dissipation function $\mfD_\eps^\wigg$ which takes the form 
\[
\mfD_\eps^\wigg(u)=\mfJ_\eps^\wigg(u,{-}\rmD_u\calE_\eps(\cdot,u){-}\Omega_\eps(u))
\text{ for } \eps\in {[0,\eps_0[},
\]
where we set $\Omega_0(u)=0$. 

\begin{corollary}[\textGamma-convergence of
  $\mfD_\eps^\wigg$]\label{co:GCvg.Dwigg} Taking the weak convergence
  $\weak$ in $\rmW^{1,p}(0,T)$ we have $\mfD_\eps^\wigg \Gweak \mfD_0^\wigg$.
\end{corollary}
\noindent
\begin{proof}
The liminf estimate for $\mfD_\eps^\wigg(u_\eps)$ with $u_\eps \weak
u_0$ in  $\rmW^{1,p}(0,T)$ follows easily from the liminf estimate for
$\mfJ_\eps^\wigg(u_\eps,\xi_\eps)$ if we use
\[
\xi_\eps= -\rmD_u \calE^\wigg_\eps(\cdot,u_\eps){+}\Omega_\eps(u_\eps) =
-\Phi'(u_\eps) +\ell -\eps \pl_u\ppp(u_\eps,\inveps u_\eps) \to \xi_0 =
-\Phi'(u_0)+\ell = - \rmD_u\calE_0(\cdot,u_0),
\]
where we used the compact embedding of $\rmW^{1,p}(0,T)$ into
$\rmC^0([0,T]) \subset \rmL^{p'}(0,T)$. 

For the limsup estimate we have to construct for each $\wh u_0$ a
recovery sequence $\wh u_\eps \weak \wh u_0$ in $\rmW^{1,p}(0,T)$ such
that $\mfD_\eps^\wigg(u_\eps)\to \mfD_0^\wigg(u_0)$. For this we set
$\wh\xi_0= -\rmD_u\calE_0(\cdot,\wh u_0)$ and use the recovery
sequence $(\wh u_\eps,\wh \xi_\eps) \weakstrong (\wh u_0,\wh \xi_0)$
such that $\mfJ_\eps^\wigg(\wh u_\eps,\wh \xi_\eps)\to
\mfJ_0^\wigg(\wh u_0,\wh\xi_0)$, whose existence is guaranteed by the
\textGamma-convergence of $\mfJ_\eps^\wigg$. 
Setting 
\[
\eta_\eps := -\rmD_u\calE_\eps (\cdot,\wh u_\eps)+\Omega_\eps(\wh
u_\eps) = -\Phi'(\wh u_\eps) +\ell - \eps \pl_u \ppp(\wh u_\eps,
\inveps \wh u_\eps)
\]
we find $\eta_\eps \to \wh\xi_0$ in $\rmL^{p'}(0,T)$, and Lemma
\ref{le:mfJ.Lip} yields $\mfJ^\wigg_\eps(\wh u_\eps,\wh\xi_\eps)-\mfJ_\eps^\wigg(\wh
u_\eps,\eta_\eps) \to 0$. Thus, we have 
\begin{align*}
&\mfD_\eps^\wigg(\wh u_\eps) -\mfD_0^\wigg(\wh u_0)\ = \ 
\mfJ_\eps(\wh u_\eps,\eta_\eps) - \mfJ_0^\wigg(\wh u_0,\xi_0)\\
&\qquad =\big(\mfJ_\eps(\wh u_\eps,\eta_\eps){-}\mfJ_\eps(\wh u_\eps,\wh\xi_\eps)\big) +
\big(\mfJ_\eps(\wh u_\eps,\wh\xi_\eps) {-} \mfJ_0^\wigg(\wh u_0,\xi_0)\big) \to 0+0.
\end{align*}
This is the desired limsup estimate. 
\end{proof} 

It remains to prove the equi-Lipschitz continuity of
$\mfJ_\eps^\wigg(u,\cdot)$ used in the above proof. 

\begin{lemma}\label{le:mfJ.Lip} If $(\R,\calE_\eps,\calR)$ satisfies
  \eqref{eq:AssumER}, then  there exists
  $C_*$ such that
\begin{align*}
&\forall\,\eps\in {[0,\eps_0[}, \; \xi,\eta\in \rmL^{p'}(0,T),\; 
u\in \rmW^{1,p}(0,T):\\ 
&\qquad \big|\mfJ_\eps^\wigg(u,\xi) {-} \mfJ_\eps^\wigg ( u,\eta)\big| \leq 
C_* \big(1{+}\|\xi\|_{\rmL^{p'}}{+}\|\eta\|_{\rmL^{p'}}\big)^{p'-1} 
  \| \xi{-}\eta \|_{\rmL^{p'}}.
\end{align*}
\end{lemma}
\begin{proof} Because $\calR^*$ is convex and has $p'$ growth
  (see \eqref{eq:p'growth}) there exists $C_*>0$ such that 
\[
\forall\, u,\xi,\eta \in \R: \quad \big|\calR^*(u,\xi)-
\calR^*(u,\eta)\big| \leq C_*(1{+}|\xi|{+}|\eta|)^{p'-1} |\xi{-}\eta| .
\]
Integration over $t\in [0,T]$ and applying H\"older's estimate gives
the desired result. 
\end{proof}

Our main result is now the relaxed EDP-convergence which follows from
the fact that the representation \eqref{eq:calM0} of $\calM $ can be
used to prove that $\calM (u,\cdot,\cdot):\R\ti \R\to {[0,\infty[}$
represents an subdifferential operator $v \mapsto
\pl_\xi\calR_\eff(u,v)$ for a uniquely defined effective dissipation
potential $\calR_\eff$. 

\begin{theorem}[Relaxed EDP-convergence]\label{th:RelaxedEDPcvg}
If the gradient systems $(\R,\calE^\wigg_\eps,\calR)$ satisfy the
assumptions \eqref{eq:AssumER} and if $\calM ^\wigg$ is defined as in
\eqref{eq:calM0}, then there exists an effective dissipation potential
$\calR_\eff$ such that \eqref{eq:relEDP.e} holds. 

Moreover, we have $(\R,\calE^\wigg_\eps,\calR)
\relEDPto (\R,\calE_0,\calR_\eff)$. 
\end{theorem} 
\begin{proof} The main parts for this proof are done in the following Sections
  \ref{se:Homog} and \ref{se:M0.prop}, and we refer to the corresponding
  results there. Nevertheless, we have the prerequisites to see the
  structure of the arguments already at this stage.  

As our energy $\calE_\eps^\wigg$ has the form $\calE_\eps^\wigg(t,u)=
\Phi(u)-\ell(t) + \ppp(u,\inveps u)$ and we set $\Omega_\eps(u)=\pl_y
\ppp(u,\inveps u)$ the conditions \eqref{eq:AssumER} easily give the
conditions \eqref{eq:relEDP.b} and \eqref{eq:relEDP.c}, where for the
second condition we use the compact embedding  
$\rmW^{1,p}(0,T) \Subset \rmC^0([0,T]) \subset \rmL^{p'}(0,T)$. 

Of course, the convergence $\mfJ_\eps^\eps \Gweak \mfJ_0^\wigg$ in
\eqref{eq:relEDP.d} is exactly what is stated in Theorem
\ref{th:GCvg.Jwigg} and proved in Section \ref{se:Homog}, whereas the
generalized Young-Fenchel estimate \eqref{eq:genYouFen} is
established in Lemma \ref{le:calM0}(b).

Proposition \ref{pr:R.eff} exactly provides the construction of
$\calR_\eff$ such that condition \eqref{eq:relEDP.e} holds.

Thus, it remains to establish the E-convergence $(\R,\calE_\eps,\calR)
\Eto (\R,\calE_0,\calR_\eff)$ (see \eqref{eq:Eto} for the definition) 
of condition \eqref{eq:relEDP.e}. For this we start with solutions
$u_\eps$ of \eqref{eq:ODE1} satisfying $u_\eps(0)\to u_0^0$ 
and exploit the standard arguments on
evolutionary $\Gamma$-convergence from \cite{SanSer04GCGF,
  Miel16EGCG}. As $u_\eps$ also satisfies the energy-dissipation
balance \eqref{eq:EDB.wig.a} we have the a priori estimate $\|
u_\eps\|_{\rmW^{1,p}(0,T)} \leq C$ and we find a subsequence with
$u_{\eps_k} \weak u_0 $ in $\rmW^{1,p}(0,T)$ which implies $u_{\eps_k}
\to u_0$ and hence $\calE_{\eps_k}(t,u_{\eps_k}(t)) \to
\calE_0(t,u_0(t))$ uniformly in $[0,T]$.  

Now we pass to the
limit $\eps_k \to 0$ in \eqref{eq:EDB.wig.a} and find 
\[
\calE_0(T,u_0(T))+ \mfJ_0^\wigg(u_0, {-}\rmD_u\calE_0(\cdot,u_0)) \leq
\calE_0(0,u_0^0 ) - \int_0^T\! \dot\ell u_0 \dd t,
\]
where we only used the liminf estimate from $\mfJ_\eps^\wigg \Gweak
\mfJ_0^\wigg$ and employed \eqref{eq:relEDP.c}. 
Now we argue as in the energy-dissipation principle (cf.\ the end of
Section \ref{su:EDP.GS}) by using the chain rule for $t \mapsto
\calE_0(t,u_0(t))$ and find $\calM ^\wigg\big(u_0,\dot
u_0,{-}\rmD_u\calE_0(t,u_0)\big)= -\rmD_u \calE_0(t,u_0)\dot u_0$. By
the definition of $\calR_\eff$ from \eqref{eq:relEDP.e} we
conclude that $0=\pl_{\dot u}\calR_\eff(u_0,\dot u_0) + \rmD_u
\calE_0(t,u_0)$ holds a.e.\ in $[0,T]$, i.e.\ $u_0$ is a solution of
the gradient system $(\R,\calE_0,\calR_\eff)$. 
\end{proof}

\section{The main homogenization result}
\label{se:Homog}

 This section contains the proof of Theorem \ref{th:GCvg.Jwigg}
which states $\mfJ_\eps \Gweakstrong \mfJ_0$, where from now on we
drop the superscript ``$\wigg$'' and always assume that the assumptions
\eqref{eq:AssumER} hold, as in the rest of the paper we only consider
the special case of our wiggly-energy model. The proof of the
technical homogenization result is obtained by extending the result of
\cite[Thm.\,3.1]{Brai02GCB}.

Before we start with  the proof of the homogenization result, we show
that the role of the variable $\xi\in \rmL^{p'}(0,T)$ is simply that
of a parameter, thus we are dealing with a parameterized
\textGamma-convergence as discussed in \cite{Miel11CDEB}.  
This comes from the fact that for $\xi$ we have strong
convergence and the functionals $\mfJ_\eps$ are equi-Lipschitz
continuous in $\xi$, as established in Lemma \ref{le:mfJ.Lip}. 
As indicated in \cite[Ex.\,3.1]{Miel11CDEB} we  see in our Example
\ref{ex:M.noncvx} that the functional $\xi \mapsto \mfJ_0(u,\cdot)$ is
not convex in general, despite of the convexity of $\mfJ_\eps(u,\cdot)$.
The following result shows that the Lipschitz
continuity in $\xi$ is preserved. We
refer to Section \ref{su:Disc.WW} for the case where the
\textGamma-limit of $\mfJ_\eps$ in the weak$\ti$weak topology which
gives a strictly lower limit that is indeed convex in $\xi$.  

\begin{lemma}[Freezing $\xi$]\label{le:freeze.xi}
(a) The weak$\ti$strong \textGamma-limit $\mfJ_0$ of $\mfJ_\eps$ exists if
and only if for all $\xi\in \rmL^{p'}(0,T)$ we have the weak
\textGamma-convergence $\mfJ_\eps(\cdot,\xi) \Gweak \mfJ_0(\cdot,\xi)$
in $\rmW^{1,p}(0,T)$.

(b) If the \textGamma-limit $\mfJ_0(\cdot,\xi_j)$ exists for
$\xi_1,\xi_2\in \rmL^{p'}(0,T)$, then for all $u\in
\rmW^{1,p}(0,T)$ we have
\begin{equation}
  \label{eq:mfJ0.equiLip}
\big| \mfJ_0(u,\xi_1)- \mfJ_0(u,\xi_2)\big| \leq 
C_* \big(1{+}\|\xi_1\|_{\rmL^{p'}}{+}\|\xi_2\|_{\rmL^{p'}}\big)^{p'-1} 
  \| \xi_1{-}\xi_2 \|_{\rmL^{p'}},
\end{equation}
where $C_*$ is from Lemma \ref{le:mfJ.Lip}. 

(c) If the weak \textGamma-limits $\mfJ_0(\cdot,\xi)$ exist for a
dense set in $\rmL^{p'}(0,T)$, then they exist for all $\xi \in
\rmL^{p'}(0,T)$.
\end{lemma}
\begin{proof}
\underline{Part (a).}  We proceed as in the proof of Corollary
\ref{co:GCvg.Dwigg}.  As $\xi_\eps \to \xi_0$ strongly, Lemma
\ref{le:mfJ.Lip} leads to 
\[
\big| \mfJ_\eps( u_\eps,\xi_\eps) - \mfJ_\eps(u_\eps,\xi_0)\big| \leq 
\wt C \| \xi_\eps {-}\xi_0\|_{\rmL^{p'}} \ \to \ 0,
\]
for $\eps \to 0$. Thus, it is easy to transfer the liminf estimate and
the construction of recovery sequences from
$\mfJ_\eps:\rmW^{1,p}(0,T)\ti \rmL^{p'}(0,T)\to \R$ to
$\mfJ_\eps(\cdot,\xi_0):\rmW^{1,p}(0,T) \to \R$ and vice versa. 

\underline{Part (b).} For the Lipschitz continuity we argue as
follows. For given $(u,\xi_j)$ 
we have a recovery sequence $(u^{(j)}_\eps,\xi_j)\weak (u,\xi_j)$ as
$\eps \to 0$, thus we have  
\begin{align*} 
\mfJ_0(u,\xi_1) - \mfJ_0(u,\xi_2) &\overset{\phantom{\text{Lem.\,\ref{le:mfJ.Lip}}}}{=} \lim_{\eps \to 0} \big(
\mfJ_\eps(u^{(1)}_\eps, \xi_1) {-} \mfJ_\eps( u^{(2)}_\eps,\xi_2)\big) \\
&\overset{\phantom{\text{Lem.\,\ref{le:mfJ.Lip}}}}{\leq} \liminf_{\eps\to 0} \big(
\mfJ_\eps(u^{(2)}_\eps, \xi_1) {-} \mfJ_\eps( u^{(2)}_\eps,\xi_2)\big)\\
&\overset{\text{Lem.\,\ref{le:mfJ.Lip}}}{\leq} \liminf_{\eps\to 0} C_* \big(
1+\|\xi_1\|_{\rmL^{p'}} + \| \xi_2\|_{\rmL^{p'}}\big)^{p'-1} \|
\xi_1{-}\xi_2 \|_{\rmL^{p'}}. 
\end{align*} 
Interchanging $\xi_1$ and $\xi_2$ we obtain the opposite result, whence 
\eqref{eq:mfJ0.equiLip} is established.

\underline{Part (c).}  Let $D\subset \rmL^{p'}(0,T)$ be the dense set
of $\xi$, for which $\mfJ_0(\cdot,\xi)$ exists. By part (b) this
function has a unique continuous extension $J:\rmW^{1,p}(0,T)\ti
\rmL^{p'}(0,T) \to \R$ that is still Lipschitz continuous in the
second variable. We have to show that this $J(\cdot,\xi)$ is indeed
the desired \textGamma-limit.
 
Given $\eta \in
\rmL^{p'}(0,T)\setminus D$ and $\delta>0$ we choose $\xi\in D$ with
$\|\xi-\eta\|_{\rmL^{p'}} \leq \delta$. For a given limit $u\in
\rmW^{1,p}(0,T)$ we first derive an approximate liminf estimate for
arbitrary $u_\eps \weak u$ via 
\begin{align*}
\liminf_{\eps \to 0} \mfJ_\eps(u_\eps,\eta)&\geq \liminf_{\eps\to 0}
\big(\mfJ_\eps(u_\eps, \xi) - \wt C \delta\big) \geq \mfJ_0 (u,\xi) -
\wt C \delta,   
\end{align*}
where $\wt C= C_*(1{+}\|\xi\|_{\rmL^{p'}}{+}\|\eta\|_{\rmL^{p'}})$.
Taking $\delta\to 0$ we have obtain the desired liminf estimate
$\liminf_{\eps \to 0} \mfJ_\eps(u_\eps,\eta) \geq \mfJ_0(u,\eta)$.

For the limsup estimate for $(\wh u,\eta)$ we have to construct a
recovery sequence $\wh u_\eps \weak \wh u$. For this we choose
$\xi^\delta \in D$ with $\| \xi^\delta - \eta\|_{\rmL^{p'}} <\delta$ and
then $\wt u^\delta_\eps\weak \wh u$ such that $\mfJ_\eps (\wt
u^\delta_\eps, \xi^\delta)\to \mfJ_0(\wh u,\xi^\delta)$ as $\eps \to 0$. 
By the equi-coercivity of
$\mfJ_\eps$ in $u$ (cf.\ \eqref{eq:Assum.R3}) all $\xi^\delta_\eps$
lie in a bounded and closed ball of $\rmW^{1,p}(0,T)$ where the weak
topology is metrizable. Hence we can extract a diagonal sequence $
\wh u_\eps = \wt u^{\delta(\eps)}_\eps \weak \wh u$ such that, using
Lemma \ref{le:mfJ.Lip} once again,  
$\mfJ_\eps (\wh u_\eps,\eta) \to \mfJ_0(\wh u,\eta)$, which is the desired
limsup estimate. 

Thus we have shown that $\mfJ_\eps(\cdot,\eta ) \Gweak
\mfJ_0(\cdot,\eta)$. 
\end{proof}

Our \textGamma-convergence result now concerns functionals of the form 
\begin{align}
\nonumber &
\mfJ_\eps(u,\xi)=\int_0^T \! N(\xi(t),u(t),\inveps u(t), \dot
u)\dd t\text{  with }\\
\label{eq:def.N}
& N(\xi,u,y,v):=\calR(u,v)+\calR^*(u,\xi{-}\pl_y\ppp(u,y)).
\end{align}
 Combining the uniform continuity estimates \eqref{eq:Assum.R4},
\eqref{eq:R*unif.cont}, the convexity and the upper bounds for $\calR$
and $\calR^*$ we easily obtain the following uniform continuity for
$N$:
\begin{equation}
  \label{eq:N.unif.cont}
\begin{aligned}  
  \exists\, C_N>0\ &\forall\, \xi_1,\xi_2,u_1,u_2,y, v_1,v_2\in \R:
  \quad \big| N(\xi_1,u_1,y,v_1)- N(\xi_2,u_2,y,v_2)\big| \\ 
 &\leq C_N\Big( \omega(|u_1{-}u_2|)
 \big(1{+}|v_1|^p{+}|v_2|^p{+} |\xi_1|^{p'} {+}|\xi_2|^{p'}\big) \\
&\qquad \quad +(1 {+}|v_1|^{p-1}{+}|v_2|^{p-1})|v_1{-}v_2| + (1
 {+}|\xi_1|^{p'-1}{+}|\xi_2|^{p'-1})| \xi_1{-}\xi_2|\Big) , 
\end{aligned}
\end{equation}
where $\omega$ is as in \eqref{eq:Assum.R4}. 

We follow the techniques in \cite[Thm.\,3.1]{Brai02GCB}, where the
case is treated that $N$ does not depend on $\xi$ and $u$.  The
generalization to the dependence on $t\mapsto \xi(t)$ with fixed $\xi$
in a dense subset $\rmC^0([0,T])$ of $\rmL^{p'}(0,T)$ and on $u=u_\eps
\weak u_0$ is handled by the uniform continuity assumption
\eqref{eq:Assum.R4}.

More importantly, we show that the limits of
``multi-cell problems'' can be replaced by a ``single-cell problem'',
which is contained in the following proposition. The essential
argument here is that we have a scalar problem for $y = z(s) = w(s)+Vs\in
\R$; in particular, the ordering properties of $\R$ together with the
1-periodicity in the variable $y$ allow us to use some simple
cut-and-paste rearrangements.

\begin{proposition}[Multi-cell versus single-cell  problem]\label{pr:N1cell} 
Consider a function $g \in \rmC(\R^2;{[0,\infty[})$ with 
\begin{subequations}
 \label{eq:gCond}
\begin{align}
 \label{eq:gCond.a}
&\forall\, v\in \R:\ g(\cdot,v) \text{ is 1-periodic},\qquad 
       \forall\, y \in \R:\ g(y,\cdot) \text{ is convex}, \\
 \label{eq:gCond.b}
&\exists\, p>1 \ \exists\, c_1,c_2>0\ \forall\, y,v\in \R:\quad
c_1\big( |v|^p -1 \big) \leq g(y,v) \leq c_2\big(1+ |v|^p\big),\\
 \label{eq:gCond.c}
&\forall\, y\in \R\ \forall\, v\in \R\setminus\{0\} : \quad g(y,v)>
g(y,0) \geq 0. 
\end{align}
\end{subequations}
Then, the following statements hold:
\\[0.2em]
(A) For all $V\in \R$ we have the identity 
\begin{subequations}
  \label{eq:HomN1cell}
\begin{align}  
\label{eq:HomN1cell.a}
G_\eff(V)&:=\lim_{L\to \infty}\  \inf\bigset{ \frac1L \int_0^L
    g\big(w(s){+}Vs, \dot w(s) {+}V\big) \dd s  }{w\in
    \rmW^{1,p}_\mafo{per}(0,T)} \\ 
\label{eq:HomN1cell.c}
&= \min\bigset{\int_0^1\! g\big(z(s), |V|\dot z(s)\big)\dd s }{
  z\in \rmW^{1,p}(0,1), \; z(1)=z(0){+} \sign(V) }. 
\end{align}
(B) Minimizers $z \in \rmW^{1,p}(0,1)$ in 
\eqref{eq:HomN1cell.c} exist  and are strictly monotone functions. 
\\[0.3em]
(C) For $V\neq0$ we have the alternative characterization 
\begin{align}
 \label{eq:HomN1cell.b}
  G_\eff(V)= \inf\bigset{ \int_{y=0}^1 g\big( y, \frac{V}{a(y)} \big)
    a(y)\dd y }{ a(y)>0, \ \int_0^1 a(y)\dd y=1},
\end{align}
\end{subequations}
and $V\mapsto G_\eff(V)$ is continuous and convex. 
\\[0.3em]
(D) If $g_1$ and $g_2$ are functions satisfying \eqref{eq:gCond} such
that 
\begin{equation}
  \label{eq:g1g2}
  \exists\, \delta_1,\delta_2 >0\ \forall\, y,v\in \R:\quad
  \big|g_1(y,v)-g_2(y,v)\big| \leq \delta_1 +\delta_2 |v|^p,
\end{equation} 
then the corresponding effective potentials $G_\eff^{(1)} $ and
$G_\eff^{(2)}$ satisfy the estimate 
\begin{equation}
 \label{eq:G1G2}
  \begin{aligned}   \forall\, v_1,v_2\in \R:\quad \big| G_\eff^{(1)}(v_1) -
    G_\eff^{(2)}(v_2)\big| \ &\leq\ \delta_1 +\frac{\delta_2}{c_1}\, 
     \big(c_1+ c_2 + c_2|v_1|^p\big)\\ 
    &\ \quad + \wh c
     \big(1{+}|v_1|^{p-1}{+}|v_2|^{p-1}\big) |v_1{-}v_2|, 
  \end{aligned} 
\end{equation}
where $\wh c$ only depends on $c_1,\,c_2$, and $p$ from
\eqref{eq:gCond.b}.
\end{proposition}
\begin{proof}
We define $G(L,V)$ to be the infimum in the right-hand side of
\eqref{eq:HomN1cell.a} and have to show $G(L,V)\to G_\eff(V)$ as $L\to
\infty$. For this we use the 1-periodicity of $g(\cdot, v)$. Moreover,
we use the coercivity of $g$ which guarantees the existence of
minimizers such that the infimum $G(L,V)$ is attained.  

We first treat the trivial case $V=0$ and then $V> 0$. The case 
$V<0$ is completely analogous to the case $V>0$. The main argument
for analyzing the minimizers is a simple cut-and-paste rearrangement
technique for 
the graph  $\bbG_z := \set{(s,z(s))\in \R^2}{s\in [0,L]}$. If we cut
this graph into finitely many pieces, we may translate these pieces
horizontally by arbitrary real numbers (using the fact that $g$ does not
depend on $s$) and may translate the pieces vertically by integer
values (using the 1-periodicity of $g(\cdot,v)$). If the result $\ol z$ is
again a graph of a continuous function, then $\ol z$ lies in
$\rmW^{1,p}(0,T)$ again and satisfies $\int_0^L g(z,\dot z)\dd s =
\int_0^L g(\ol z, \dot{\ol z}) \dd s$.

\underline{Step 1. The case $V=0$.}\\
We first observe that $G(L,0)= g_\mafo{min}:=\min\set{g(y,0)}{y\in
  \R}$, since $g(y,v)\geq g_\mafo{min}$ and we can choose $w\equiv
y_*$ with $g(y_*,0)=g_\mafo{min}$.  The minimizer $z$ for
\eqref{eq:HomN1cell} is given by $z\equiv y_*$.

\underline{Step 2. Monotonicity of $z: s\mapsto w(s){+}sV$.} Here we
consider general minimizers $w$ for $G(L,V)$ with $V>0$ and $LV\geq
1$. To show that $z$ is increasing, we assume that there exist $s_1$ and
$s_2$ with $0\leq s_1<s_2\leq L$ and $z(s_1)=z(s_2)$ such that
$z|_{[s_1,s_2]}$ is not constant. From this we produce a contradiction. 

With $y_*$ from Step 1 and using $LV\geq 1$ the intermediate-value
theorem provides $s_*\in [0,L]\setminus {]s_1,s_2[}$ such that
$z(s_*)=y_*$. We then have
\begin{equation}
  \label{eq:NonMonotone}
  \int_{s_1}^{s_2} g(z(s),\dot z(s)) \dd s > \int_{s_1}^{s_2}
g(z(s),0)\dd s \geq \int_{s_1}^{s_2} g(y_*,0)\dd s ,
\end{equation}
where the strict estimate ``$>$'' holds since $z$ is not constant on
this interval and $g$ satisfies \eqref{eq:gCond.c}.  We now define a
function $\wt z\in \rmW^{1,p}(0,L)$ by cutting out the non-monotone
interval ${]s_1,s_2[}$ and inserting a flat part where $\wt z$
takes the value $y_*$, see Figure \ref{fig:min.timerescale}. 
E.g.\
for the case $s_*\geq s_2$ we obtain
\[
\wt z(s)= \left\{ \ba{cl} z(s) &\text{for } s\in [0,L]\setminus{]s_1,s_*[},
  \\  z(s{+}s_2{-}s_1) &\text{for }s\in [s_1,s_1{+}s_*{-}s_2],\\ 
     y_*&\text{for } [s_1{+}s_*{-}s_2, s_*].\ea  \right. 
\]
The case $s_*\leq s_1$ is similar.
\begin{figure}
\centering 

\begin{tikzpicture}
\draw[very thick, color=black!50] (0,0)--(5,0)--(5,4)--(0,4)--(0,0);
\draw[thick,->] (0,0) --(5.5,0) node[right] {$s$};
\draw[thick,->] (0,0) --(0,4.5) node[above] {$z(s)$};

\node[below] at (0,0) {$0$};
\node[left] at (0,0) {$z(0)$};
\node[below] at (5,0) {$L$};
\node[left] at (0,4) {$z(0){+}LV$};
\node[left] at (0,3) {$y_*$};

\draw (1,2)--(1,-0.1) node[below] {$s_1$}; 
\draw (3,2)--(3,-0.1) node[below] {$s_2$}; 
\draw (-0.1,3) --(4,3)--(4,-0.1) node[below] {$s_*$};

\draw[ultra thick] (0,0)--(1,2)--(2,1)--(5,4);
\draw[ultra thick, color=red] (0,0)--(1,2);
\draw[ultra thick, color=blue] (3,2)--(5,4);

\end{tikzpicture}
\begin{tikzpicture}
\draw[very thick, color=black!50] (0,0)--(5,0)--(5,4)--(0,4)--(0,0);
\draw[thick,->] (0,0) --(5.5,0) node[right] {$s$};
\draw[thick,->] (0,0) --(0,4.5) node[above] {$\wt z(s)$};

\node[below] at (0,0) {$0$};
\node[left] at (0,0) {$z(0)$};
\node[below] at (5,0) {$L$};
\node[left] at (0,4) {$z(0){+}LV$};
\node[left] at (0,3) {$y_*$};

\draw (1,2)--(1,-0.1) node[below] {$s_1$}; 
\draw (-0.1,3) --(4,3)--(4,-0.1) node[below] {$s_*$};

\draw[very thick,dashed, color=black!30] (0,0)--(1,2)--(2,1)--(5,4);
\draw[ultra thick, color=red] (0,0)--(1,2);
\draw[ultra thick, color=blue] (1,2)--(2,3);
\draw[ultra thick, color=orange] (2,3)--(4,3) node[pos=0.5, below]
{\color{black}$ \underbrace{\mbox{}\hspace{1.8cm}\mbox{}}_{s_2{-}s_1}$};
\draw[ultra thick, color=blue] (4,3)--(5,4);

\end{tikzpicture}

\caption{The new function $\wt z$ (right side) is constructed
  from the not monotonically increasing function $z$ (left side) by removing
  the non-monotone part on $[s_1,s_2]$ and inserting a flat part of the
  same length $s_2{-}s_1$ with value $\wt z(s)=y_* \in \mafo{armgin}\,
  g(\cdot,0)$.}
\label{fig:min.timerescale}
\end{figure}
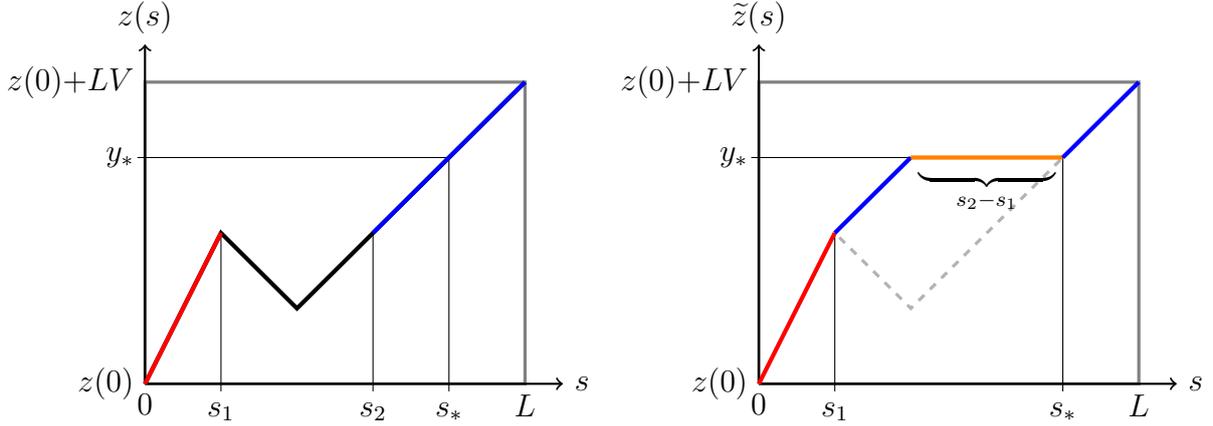
By construction we have $\wt z \in \rmW^{1,p}(0,L)$ and $\wt z(L)=\wt
z(0){+} LV$. Hence, $\wt z$ is a competitor for the minimization
problem $G(L,V)$. Because \eqref{eq:NonMonotone} implies 
$\int_0^L  g(z,\dot z)\dd s > \int_0^L g(\wt z, \dot{\wt z}) \dd s $
we see that $z$ cannot be a minimizer, which is the desired
contradiction. Thus, statement (B) is shown.  

\underline{Step 3. Claim: $\forall\,V> 0 \ \forall\, k\in \N$ with $k/V\geq 1$ we have $G(k/V,V)=G(1/V,V)$.} 
\\
We start from a minimizer $w_V$ for $G(1/V,V)$ and use the
1-periodicity of $g(\cdot,v)$. Extending $w_V$ periodically to $w^k_V
\in \rmW^{1,p}_\mafo{per}(0,k/V)$ we can insert it as competitor for
$G(k/V,V)$ and conclude $G( k/V,V)\leq G(1/V, V)$.  

For the opposite estimate consider a fixed $k\geq 2$ and 
take a minimizer $w \in
\rmW^{1,p}_\mafo{per}(0,k/V)$ for $G( k/V,V)$. We extend $w$
periodically to all of $\R$ and  define $z: \R\ni s\mapsto w(s)+sV$. 
and 
\[
T:= \set{s_2{-}s_1 }{ s_1,s_2\in \R,\ z(s_2)=z(s_1)+1} \quad
\text{ and } \tau_* := \min T.
\]
The set $T$ is non-empty as $z(k/V)=z(0)+k$.
By Step 2 $z$ is monotone, hence $\tau_*>0$ by using periodicity which
gives compactness. Choosing $s_j$ with $z(s_j)=z(0)+j$ for
$j=1,...,k{-}1$ and setting $s_0=0$ and $s_k=k/V$, we have $
k/V= \sum_{j=1}^k (s_j{-}s_{j-1})$. Thus, at least one $s_j{-}s_{j-1}$
is less or equal $1/V$, which implies $\tau_*\leq 1/V$. 

By shifting $z$ horizontally, we may assume $z(\tau_*)=z(0){+}1$. 
If $\tau_*=1/V$ we have $z(1/V)=z(0)+1$ such that $w:s\mapsto z(s)-Vs$
satisfies $w(0)=w(1/V)=w(k/V)$.  Hence, $w|_{[0,1/V]} $ is a
competitor for $G(1/V,V)$, and $w_{[1/V,k/V]}$ is a competitor for
$G((k{-}1)/V,V)$ (after shifting $s$ to $s-1/V$). Hence, we obtain 
\begin{equation}
  \label{eq:Gk.1.k-1}
\begin{aligned}
\frac kV G(k/V,V) &= \int_0^{1/V} g(w{+}Vs,\dot w{+}V)\dd s +
\int_{1/V}^{k/V} g(w{+}Vs,\dot w{+}V)\dd s\\ 
&\geq \frac1V G(1/V,V) +
\frac{k{-}1} V G((k{-}1)/V,V).
\end{aligned}
\end{equation}

We want to show the same lower bound for the case $\tau_*<1/V$. This
is done by a cut-and-paste rearrangement. We decompose $[0,k/V]$
into at most 5 parts via $0<t_1<t_2<t_3<t_4\leq k/V$.  We set
$t_2:=\tau_*<t_3:=1/V$ and choose $t_4 > 1/V$ such that
$z(t_4)=z(0)+j_*$ with $j_*\geq 2$ and $z(t_4{-}t_3)\geq z(0)+j_*-1$.
Now the intermediate-value theorem applied to the difference of
$z|_{0,\tau_*}$ and $\ol z:[0,\tau_*]\ni s \mapsto
z(t_4{-}t_3{+}s)-j_*+1$ gives at least one zero $t_1 \in [0,\tau_*]$
as $z(0)\leq \ol z (0) = z(t_4{-}t_3) -j_*+1$ and
$\ol z (t_3) =z(\tau_*)\leq z(t_3) $ by
monotonicity. 

We define the rearrangement $\wh z$ concatenation of
vertically shifted versions of $z$ on the intervals $[0,t_1]$,
$[t_3,t_4]$, $[t_2,t_3]$, $[t_1,t_2]$, and $[t_4,k/V]$, namely 
\[
 \wh z(s)= \left\{ \ba{cl} 
   z(s) &\text{for } s \in [0,t_1]\cup [t_4,k/V], \\  
  z(s{+}t_4{-}t_3)-j_*+1& \text{for } s \in [t_1,t'_2], \\
  z(s{+}t_2{-}t_3) &\text{for }s\in [t'_2,t'_3],\\
  z(s{+}t_2{-}t_4)+j_*-1&\text{for }s\in [t'_3,t_4], 
\ea \right. 
\] 
where $t'_2= t_3$ and $t'_3=t_4-t_2+t_1$.
See Figure \ref{fig:min.rearrange} for an illustration. 
\begin{figure}
\centering 
\begin{tikzpicture}[scale=0.8]
 \draw[color=black!30, thick] (0,0) -- (6,0) --(6,9) --(0,9)--(0,0);
 \draw[color=black!30, thick] (2,0) -- (4,0) --(4,9) --(2,9)--(2,0);
 \draw[color=black!30, thick] (0,3) -- (6,3) --(6,6) --(0,6)--(0,3);
 \draw[color=black!30, thick] (0,0) -- (6,9);
 \draw[->] (0,0) -- (6.3,0) node[right] {$s$};
 \draw[->] (0,0) -- (0,9.3) node[above] {$z(s)$};

 \draw[ultra thick] (0,0)--(2,4)--(3.5,5)--(4.5,6)--(6,9);
 \draw[ultra thick,color=red] (1,2)--(1.5,3);
 \draw[ultra thick,color=blue] (1.5,3)--(2,4)--(3.5,5);
 \draw[ultra thick,color=green] (3.5,5)--(4.5,6);

 \draw (1,2.2) -- (1,-0.2) node[below] {$t_1$};
 \draw (1.5,3.2)--(1.5,-0.2) node[below] {$t_2$};
 \draw (3.5,5.2)--(3.5,-0.2) node[below] {$\wh t_3$};
 \draw (4.5,6.2)--(4.5,-0.2) node[below] {$t_4$};
 \node[below] at (2,0) {$\frac 1V$}; 
 \node[below] at (4,0) {$\frac 2V$};
 \node[below] at (6,0) {$\frac 3V$};
 \node[left] at (0,3) {$1$}; 
 \node[left] at (0,6) {$2$}; 
 \node[left] at (0,9) {$3$}; 

\end{tikzpicture}    
\begin{tikzpicture}[scale=0.8]
 \draw[color=black!30, thick] (0,0) -- (6,0) --(6,9) --(0,9)--(0,0);
 \draw[color=black!30, thick] (2,0) -- (4,0) --(4,9) --(2,9)--(2,0);
 \draw[color=black!30, thick] (0,3) -- (6,3) --(6,6) --(0,6)--(0,3);
 \draw[color=black!30, thick] (0,0) -- (6,9);
 \draw[->] (0,0) -- (6.3,0) node[right] {$s$};
 \draw[->] (0,0) -- (0,9.3) node[above] {$\wh z(s)$};

 \draw[ultra thick, dashed, color=black!30] 
           (0,0)--(2,4)--(3.5,5)--(4.5,6)--(6,9);
 \draw[ultra thick] (0,0)--(1,2);
 \draw[ultra thick,color=green] (1,2)--(2,3);
 \draw[ultra thick,color=blue] (2,3)--(2.5,4)--(4,5);
 \draw[ultra thick,color=red] (4,5)--(4.5,6);
 \draw[ultra thick] (4.5,6)--(6,9);

 \draw (1,6.1) -- (1,-0.2) node[below] {$t_1$};
 \draw (2,6.1)--(2,-0.0);
 \node[above] at (1.5,6)
           {$\overbrace{\mbox{}\hspace*{0.9cm}\mbox{}}^{t_4{-}\wh t_3}$}; 
 \draw (4,7.6)--(4,-0.0);
 \node[above] at (3,0.0)
           {$\underbrace{\mbox{}\hspace*{1.8cm}\mbox{}}_{\wh t_3{-}t_2}$};
 \draw (4.5,7.6)--(4.5,-0.2) node[below] {$t_4$};
 \node[above] at (4.25,7.6)
           {$\overbrace{\mbox{}\hspace*{0.4cm}\mbox{}}^{t_2{-}t_1}$};
 \node[below] at (2,0) {$\frac 1V$}; 
 \node[below] at (4,0) {$\frac 2V$};
 \node[below] at (6,0) {$\frac 3V$};
 \node[left] at (0,3) {$1$}; 
 \node[left] at (0,6) {$2$}; 
 \node[left] at (0,9) {$3$}; 
 \fill (2,3) circle (3pt);
\end{tikzpicture}
\caption{Rearrangement of $z$ leads to $\wh z$, which intersect the
  diagonal $s\mapsto z(0)+Vs$ at $s=1/V$ (filled circle). With $\wh t_3 = t_4-t_3+t_1$, the parts of
  the graph associated with $[t_1,t_2]$ and $[\wh t_3,t_4]$ are
  interchanged by vertical integer-valued shifting and horizontal
  adjustment to make the function continuous.} 
\label{fig:min.rearrange}
\end{figure}
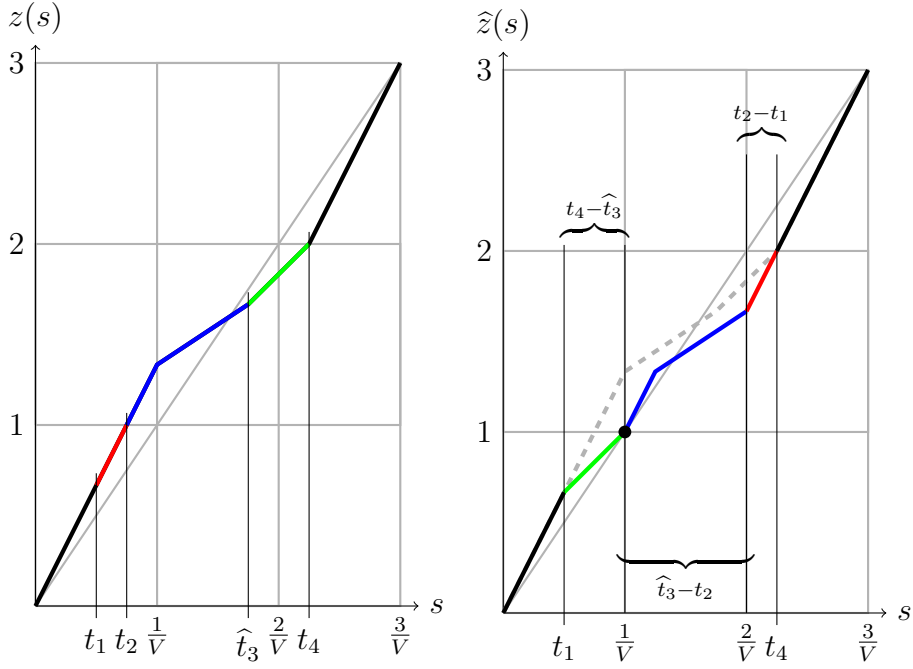

 By construction $z$ and $\wh z$ are
minimizers for $G(k/V)$, but $\wh z$ additionally satisfies $\wh
z(1/V)=\wh z(0)+1$, as in the case $\tau_*=1/V$. By induction we 
 find $ G(k/V,V)\geq  G(1/V,V)$. Since the opposite
estimate was shown above, we conclude $G(k/V,V)=G(1/V,V)$. 
\medskip

\underline{Step 4. Limit $G(L,V)\to G(1/V,V)$ for $L\to \infty$.}
\\
We already know the values at $G(k/V,V)=G(1/V,V),$ and now estimate
the function for $L\in {]k/V,(k{+}1)/V[}$. Using
$g^*_V=\max\set{g(u,V)}{u\in \R}$ and taking the minimizer
$z_L$ for $G(L,V)$ we extend $z_L \in \rmW^{1,p}(0,L)$ to $\wt z\in
\rmW^{1,p}(0,(k{+}1)/V)$ via $\wt z(s)=z(0)+sV$ for $s>L$, then  
\begin{align*}
L \,G(L,V)&=\int_0^{L} g(z_L,\dot z_L)\dd s \geq \int_0^{(k+1)/V} g(\wt z ,
\dot{\wt z}) \dd s -g^*_V \big( \frac{k{+}1}V-L\big) \\
&\geq \frac{k{+}1}V G((k{+}1)/V,V) -g^*_V/V \ \geq L \, G(1/V,V)
- g^*_V /V.
\end{align*} 
This implies $\liminf_{L\to \infty} G(L,V)\geq G(1/V,V)$. The opposite
inequality follows by taking the minimizer $z_{k/V}$ and extending it
affinely to a competitor for $G(L,V)$. This results in 
$\frac kV G(1/V,V)= \frac kV G(k/V,V) \geq LG(L,V) - g^*_V/V$ and
$\limsup_{L\to \infty} G(L,V) \leq G(1/V,V)$ follows, and $G(L,V) \to
G(1/V,V)$ is established. 

To establish the identity \eqref{eq:HomN1cell} we simply observe that the
minimizers $z$ of \eqref{eq:HomN1cell.c} and the minimizers $w$ of
$G(1/V,V)$  are related by $z(s)= w(|V|s)+\sign(V) \,s$. Thus, part
(A) is established. 
 
\underline{Step 5. Convexity of $G_\eff$.}\\
Obviously monotone functions $s\mapsto z(s)$ as competitors in
\eqref{eq:HomN1cell.c} can be approximated by strictly monotone
functions in $\rmW^{1,p}$. For these functions we can invert $y=z(s)$
to obtain $s=\sigma(y)$. Thus for $a(y)=sign(V)\sigma'(y)$ we have $a(y)>0$
and $\int_0^1 a(y) \dd y =1$. Thus, transforming the integral in
\eqref{eq:HomN1cell.c} gives the desired formula
\eqref{eq:HomN1cell.b}. 

The convexity of $g(y,\cdot)$ implies the convexity of $(v,a)\mapsto
g(y,v/a)a=:h(y,a,v)$. With this we set $\calH(a,v)=\int_0^1
h(y,a(y),v)\dd y$, which is still convex in $(a,v)$.
Thus, for $\theta\in {]0,1[}$ and $v_0,v_1\in
\R$ we choose for $\eps>0$ functions $a_0$ and $a_1$ such that 
$\calH(a_j,v_j)\leq G_\eff(v_j)+\eps$ for $j=0$ and $1$.  
For $v_\theta= (1{-}\theta)v_0+\theta v_1$ we obtain
\begin{align*}
  G_\eff(v_\theta)&=\inf\bigset{ \calH(a,v_\theta)}{\int_0^1 a(y) \dd
    y =1} \ \leq \calH\big((1{-}\theta)a_0{+}\theta
  a_1,v_\theta\big)\\
& \overset{\calH\text{ cvx}}\leq (1{-}\theta) \calH(a_0,v_0) + 
\theta \calH(a_1,v_1) \leq (1{-}\theta) G_\eff(v_0)+ \theta
G_\eff(v_1) + \eps. 
\end{align*}
As $\eps>0$ was arbitrary the desired convexity is established. 

\underline{Step 6. Continuous dependence of $G_\eff$ from $g$.}\\
To obtain \eqref{eq:G1G2} we first consider the case $v_1=v_2=V$ and
denote by $z_j$ any minimizers for $G_j(1/V,V)$. By comparing with
$\wh z(s)=\sign(V) s$ we first obtain the upper bound
\[
G_\eff^{(j)}(V)=G_j(1/V,V)=\int_0^1 g_j(z_j,|V|\dot z_j)\dd s \leq
\int_0^1 g_j(s,|V|)\dd s \leq c_2(1{+}|V|^p).
\]
Second, using the lower bound for $g_j$ we find 
\[
G_\eff^{(j)}(V)=\int_0^1 g_j(z_j,|V|\dot z_j)\dd s\geq
c_1|V|^p\int_0^1|\dot z_j|^p \dd s - c_1,
\]
which gives the a priori estimate  $c_1|V|^p\int_0^1|\dot z_j|^p \dd
s\leq c_1+c_2+ c_2|V|^p$. 
Now we compare the two effective potentials as follows 
\begin{align*}
&G_\eff^{(2)}(V){}-G_\eff^{(1)}(V)=\int_0^1 \big(g_2(z_2,|V|\dot z_2)-
g_1(z_1,|V|\dot z_1) \big) \dd s \\
&\quad \leq \int_0^1 \big(g_2(z_1,|V|\dot z_1)-
g_1(z_1,|V|\dot z_1) \big) \dd s \ \leq \  \int_0^1 \big(\delta_1+
\delta_2 |V|^p|\dot z_1|^p\big) \dd s\\  
&\quad  = \ \delta_1 + 
 \delta_2 |V|^p \int_0^1|\dot z_1|^p \dd s
 \  \leq \ \delta_1  +\frac{\delta_2}{c_1} \big(c_1+c_2 + c_2|V|^p\big).
\end{align*}
By interchanging 1 and 2, we obtain the same bound for
$G_\eff^{(1)}(V)-G_\eff^{(2)}(V)$ and \eqref{eq:G1G2} is established
for $v_1=v_2=V$. 

By the triangle inequality it suffices to estimate
$G_\eff^{(1)}(v_1)-G_\eff^{(1)}(v_2)$. For this we can use that
$G_\eff^{(1)}$ is convex according to part (C) and satisfies the
bounds $0\leq G_\eff^{(1)}(V)\leq c_2(1{+}|V|^p)$. Thus, 
\[
|G_\eff^{(1)}(v_1)-G_\eff^{(1)}(v_2)| \leq \wh
c(1{+}|v_1|^{p-1}{+}|v_2|^{p-1}) | v_1{-}v_2|
\]
follows by standard convexity theory. Hence, part (D) is established
as well.
\end{proof}

\begin{remark}[Non-uniqueness without monotonicity]
  \label{re:UniqStrictM} Minimizers in \eqref{eq:HomN1cell.c}
are neither unique nor strictly monotone for functionals based on
$g(y,v)= \max\{|v|,v^2\}$. For $V=1/2$ we have the minimizers
$z(s)=s/2$ as well as $z(s)=\min\{s,1/2\}$. So, our assumption on
strict convexity is indeed important. 
\end{remark} 

As a consequence of Proposition \eqref{pr:N1cell}(C) we obtain a very
useful uniform continuity for the effective contact potential $\calM
$. For this, we recall that $\calM (U,V,\Xi)$ (cf.\ \eqref{eq:calM0})
is obtained by setting $g^{U,\Xi} (y,v) = N(\Xi,U,y,v)$, then $\calM
(U,V,\Xi)=G^{U,\Xi}_\eff(V)$. Exploiting the continuity of $N$ (see
\eqref{eq:N.unif.cont}), we obtain the following result.

\begin{corollary}[Continuity of $\calM $]\label{co:M0.cont}
If $N$ (see \eqref{eq:def.N}) satisfies  \eqref{eq:N.unif.cont}, then 
there exists $C_\calM>0$ such that $\calM $ (see \eqref{eq:I.rmM0}) satisfies
\begin{equation}
  \label{eq:M0.contin}
\begin{aligned}  
\forall\, v_j,\xi_j&\in \R: \quad 
\big|\calM (u_1,v_1,\xi_1) - \calM (u_2,v_2,\xi_2)\big|  \\ 
& \leq 
C_\calM \Big(\omega(|u_1{-}u_2|)
 \big( 1 {+}|v_1|^p {+}|v_2|^p {+}|\xi_1|^{p'} {+}|\xi_2|^{p'}
 \big)  \\ 
&\quad \qquad 
 + \big(1 {+} |v_1|^{p-1} {+}|v_2|^{p-1}\big)|v_1{-}v_2| + 
 \big( 1 {+} |\xi_1|^{p'-1} {+}|\xi_2|^{p'-1} \big) |\xi_1{-}\xi_2|\Big),
\end{aligned}
\end{equation}
where $\omega$ is from \eqref{eq:Assum.R4}.
\end{corollary} 
\begin{proof} We simply apply part (D) of Proposition \ref{pr:N1cell} 
with $g_j(y,v)=N(\xi_j,u_j,y,v)$. Then, inserting
\eqref{eq:N.unif.cont} into \eqref{eq:g1g2} allows us to conclude 
\eqref{eq:G1G2}, which is indeed the desired estimate
\eqref{eq:M0.contin}, because
$\calM(u_j,v_j,\xi_j)=G^{u_j,\xi_j}_\eff(v_j)$. 
\end{proof}

We have now prepared all the tools for first showing the desired
liminf estimate and then the limsup estimate by constructing suitable
recovery sequences. Both results are suitable generalizations of
\cite[Thm.\,3.1]{Brai02GCB}.  (Recall that we dropped the superscript
$^\text{wig}$ which was used in Section \ref{s:Model}.)

\begin{proposition}[The liminf estimate] \label{prop-liminf} Let
  $\mfJ_\eps, \;\mfJ_0:\rmW^{1,p}(0,T)^2 \to \R$ be defined as in
  \eqref{eq:EDB.wig.c} and \eqref{eq:mfJ0}, respectively. Then,
 \[
(u_\eps,\xi_\eps)\weakstrong (u_0,\xi_0)\text{ \  in  } 
\rmW^{1,p}(0,T)\ti \rmL^{p'}(0,T) \quad \Longrightarrow \quad 
\mfJ_0(u_0,\xi_0) \leq \liminf\limits_{\eps\searrow 0} \mfJ_\eps(u_\eps,\xi_\eps) .
\]
\end{proposition}
 \begin{proof} 
By Lemma \ref{le:freeze.xi} we know that it suffices to consider
$\xi_\eps = \xi$ with $\xi\in \rmC^0([0,T))$. We keep this choice
fixed for the rest of the proof. Moreover, we keep $u_0\in
\rmW^{1,p}(0,T)\subset \rmC^0([0,T])$ fixed. 

The main idea is to use continuity in time of $\xi$ and $u_0$ as
well as the uniform convergence $\|u_\eps{-}u_0\|_{\rmL^\infty}\to 0$
to approximate 
\[
N(\xi(t),u_\eps(t),\inveps u_\eps(t),\dot
u_\eps(t)) \quad \text{by} \quad  
N(\xi(t_j),u_0(t_j),\inveps u_\eps(t),\dot
u_\eps(t))
\]
on suitable subintervals $[t_{j-1},t_j] \subset [0,T]$. 
By \eqref{eq:N.unif.cont}
for every $\delta>0$ we find $\eta>0$ with
\begin{subequations}
\label{eq:N.estim}
\begin{equation}
  \label{eq:N.est.a}
  |\xi{-}\wh\xi|+ |u{-}\wh u| <\eta \quad \Longrightarrow \quad 
\big|N(\xi,u,y,v)-N(\wh\xi, \wh u, y, v)\big| \leq \delta\big( 1
+|\xi|^{p'} + |v|^p\big).
\end{equation}
We now fix an arbitrary $\delta>0$, which finally can be made as small
as we like.  

We define a partition $0=t_0<t_1< \cdots < t_{n-1}<t_n=T$ such that 
\begin{equation}
  \label{eq:N.est.b}
  |\xi(t){-}\xi(t_j)| <\eta/3 \text{ and } |u_0(t){-}u_0(t_j)|<
\eta/3 \quad\text{for } t\in [t_{j-1},t_j] \text{ and }j=1,...,n.
\end{equation}
\end{subequations}
Moreover, we choose $\eps_1>0$ such that $\|u_\eps {-}
u_0\|_{\rmL^{\infty}} <\eta/3 $ for $\eps\in {]0,\eps_1[}$. 

Then, we can estimate $\mfJ_\eps(u_\eps,\xi)$ from below as follows
\begin{align*}
\mfJ_\eps(u_\eps, \xi)&=\sum_{j=1}^n \int_{t_{j-1}}^{t_j}
N(\xi(t),u_\eps(t), \inveps u_\eps(t),\dot u_\eps(t)) \dd t \\
& \geq \sum_{j=1}^n \int_{t_{j-1}}^{t_j} \Big( N(\xi(t_j), u_0(t_j),
\inveps u_\eps(t),\dot u_\eps(t)) - \delta (1+\|\xi\|_\infty^{p'}+
|\dot u_\eps(t)|^p\Big) \dd t .
\end{align*}  
Because $u_\eps \weak u_0$ we have $\| \dot u_\eps\|_{\rmL^p}^p\leq
C_{\dot U}<\infty$, and hence can pass to the liminf for $\eps\searrow
0$ by using \cite[Thm.\,3.1]{Brai02GCB} for each of the summands
$j=1,...,n$ separately:
\[
\liminf_{\eps \to 0} \mfJ_\eps(u_\eps, \xi) \geq 
\sum_{j=1}^n \int_{t_{j-1}}^{t_j} \calM (u_0(t_j),\dot u_0(t) ,
\xi(t_j))\dd t -\delta T\big(1{+}\|\xi\|_\infty^{p'}{+} C_{\dot U}
\big) 
\] 
Here we used that $g^{u,\xi}(y,v)=N(\xi,u,y,v)$ in Proposition
\ref{pr:N1cell} giving $G^{u,\xi}_\eff(V)=\calM (u,V,\xi)$. 
Employing the uniform continuity of $\calM $ established in
\eqref{eq:M0.contin} yields 
\[
 \big| \calM (u_0(t),V,\xi(t)) -
 \calM (u_0(t_j),V,\xi(t_j))\big|   
\leq C \delta\,(1+ |V|^p).
\]
Thus, we can further estimate from below as follows
\begin{align*}
\liminf_{\eps \to 0} \mfJ_\eps(u_\eps, \xi) &\geq 
\sum_{j=1}^n \int_{t_{j-1}}^{t_j}\!\!\! \big( \calM (u_0(t),\dot u_0(t) ,
\xi(t)) {-} C\delta (1{+}|\dot u_0(t)|^p) \big)\dd t -  
\delta T\big(1{+}\|\xi\|_\infty^{p'}{+} C_{\dot U} \big) \\
&= \mfJ_0(u_0,\xi) - \delta \wh C. 
\end{align*}
As $\delta>0$ can be chosen arbitrarily small, the desired liminf
estimate is established. 
\end{proof}

The final limsup estimate is obtained by providing recovery sequences
for piecewise affine functions $\wh u$ and piecewise constant
functions $\wh\xi$ and exploiting a
standard density argument. So we can use that $V=\dot{\wh u}(t)$ and $\Xi=\wh
\xi(t)$ are constant in a macroscopic subinterval, but the
construction of recovery sequences is still 
complicated as $t \mapsto \wh u(t)$ is not constant. So locally on the
scale $O(\eps)$ we approximate via $\wh u_\eps(t) \approx \wh u(t_*)+ \eps
Z(t_*,\inveps (t{-}t_*))$, where $Z(t_*,\cdot)$ is 
obtained from the minimizers $z\in \rmW^{1,p}(0,1)$ for $\calM (\wh u(t_*),\wh
u'(t_*),\wh\xi(t_*))$ (cf.\ \eqref{eq:calM0}). We  keep such an
approximation on intervals of length $\eps^{1/2}$ and adjust $\wh
u(t_*)$ then on the neighboring intervals. 

Indeed, for given $(U,V,\Xi)
\in \R^3$ we take a minimizer $z_{U,V,\Xi}\in \rmW^{1,p}(0,1)$, where for $V\neq
0$ we may assume $z(0)=0$ without loss
of generality. For $V\neq 0$ we define the
\emph{shape functions} $Z_{U,V,\Xi}:\R\to \R$ via 
\begin{align}
\label{eq:ShapeFcn}
Z_{U,V,\Xi}(t) := z_{U,V,\Xi}(|V|t) \text{ for }t\in
[0,\tfrac1{|V|}],\quad Z_{U,V,\Xi}(t{+}\tfrac kV)= Z_{U,V,\Xi}(t){+}k.
\end{align} 
Note that the definition of $Z_{U,V,\Xi}$ is such that $\R\ni t\mapsto
Z_{U,V,\Xi}(t)-Vt$ is periodic with period $1/|V|$.

\begin{proposition}[The limsup estimate, recovery
  sequences] \label{prop-limsup}  
For all pairs $(\wh u, \wh \xi)\in\rmW^{1,p}(0,T)  \ti \rmL^{p'}(0,T)$ 
there exists a recovery sequence $\wh u_\eps \weak \wh u$ in
$\rmW^{1,p}(0,T)$ such that for all  $\wh\xi_\eps \to \wh \xi$ in
$\rmL^{p'}(0,T)$ we have 
$\mfJ_\eps(\wh u_\eps,\wh\xi_\eps)\rightarrow\mfJ_0(\wh u,\wh\xi)$. 
\end{proposition}
\begin{proof}
 \underline{Step 1: Continuity of $\mfJ_0$.} Using 
the uniform continuity of $\calM $ established in
\eqref{eq:M0.contin}, we easily obtain that $\mfJ_0: \rmW^{1,p}(0,T)\ti
\rmL^{p'}(0,T) \to \R$ is continuous in the norm topology. Thus, by
standard arguments of \textGamma-convergence it suffices to provide
the construction of a recovery sequences for $(\wh u,\wh \xi)$ on a subset of 
$\rmW^{1,p}(0,T)\ti\rmL^{p'}(0,T)$
that is dense in the norm topology. Then, the same diagonal argument
as in the proof of Lemma \ref{le:freeze.xi}(c) can be applied. 

\underline{Step 2: Restriction to a dense subset $D \subset \rmW^{1,p}(0,T)\ti
\rmL^{p'}(0,T)$.} We define $D$ as follows. We consider dyadic
partitions $\set{t_{j,N}:= kT/2^N }{ k=0,...,2^N}$ of $[0,T]$ and
assume that pairs $( \wh u,\wh\xi)$ in $D$ are such that $\dot{\wh u}$ and
$\wh \xi$ are constant on the intervals
${]t_{j-1,N},t_{j,N}[}$. Moreover, we assume that the slopes $
V_{j,N}= \dot{\wh u}(t)$ for $t \in {]t_{j-1,N},t_{j,N}[}$ are non-zero. 
By standard arguments we see that $D$ is dense in $\rmW^{1,p}(0,T)\ti
\rmL^{p'}(0,T)$.

As all $\mfJ_\eps$ and $\mfJ_0$ are integral functionals it is
now sufficient to give the recovery construction of a $(\wh
u,\wh\xi)\in D$ on  one subinterval $[t_{j-1,N},t_{j,N}]$. For $\wh u$
we take care that the values at both ends remain unchanged, so
that joining the different constructions stays in $\rmW^{1,p}(0,T)$. 

\underline{Step 3: Recovery construction.} To simplify the notation we
write $[a,b]=[t_{j-1,N},t_{j,N}]$, $V=\frac1{b-a}(\wh u(b){-}\wh
u(a))$, and  $\wh\xi(t)=\Xi$. We use the shape functions
$Z_{U,V,\Xi}$ introduced in \eqref{eq:ShapeFcn} for the fixed values
$V$ and $\Xi$, but still need to adjust $U$ accordingly. This is done
on the intermediate scale $\eps^{1/2}$, i.e.\ we divide $[a,b]$ in 
\[
n_\eps:= \big\lfloor \frac{b{-}a} {\eps^{1/2}} \big\rfloor \quad
\text{(floor function)},
\]
subintervals of equal length via $a^\eps_k:= a+ k(b{-}a)/n_\eps$. 
Letting $U^\eps_k=\wh u(a^\eps_k)$ for $k=0,1,...,n_\eps$ we assume for simplicity
$U^\eps_k\in\eps\Z$ and we define the
approximation $\wh u_\eps:[a^\eps_{k-1},a^\eps_k]\to \R$ via 
\[
\wh u_\eps(t)= \left\{\ba{cl} 
 U^\eps_{k-1} + \eps Z_{U_k,V,\Xi}\big(\inveps (t-a^\eps_{k-1})\big)  & \text{for }
a^\eps_{k-1} \leq t \leq x^\eps_k ,\\ 
   U^\eps_{k}+V(t{-}a^\eps_k)=\wh u(t)&  \text{for } x^\eps_k \leq t\leq a^\eps_k, 
\ea \right.
\]
where $x^\eps_k := a^\eps_{k-1} +\frac{\eps}{|V|}\big\lfloor \frac{|V|(a^\eps_k{-}a^\eps_{k-1})}{\eps} \big\rfloor$. 
The number of used periods of the shape function $Z_{U,V,\xi}$ behaves
like $1/(\eps^{1/2}|V|)\to \infty$ and covers
$[a_{k-1}^\eps,x_k^\eps]$, which is most of the interval
$[a_{k-1}^\eps,a_k^\eps]$, while the remaining part
$[x^\eps_k,a^\eps_k]$ with $\wh u_\eps=\wh u$ has at most of length
$\eps |V|$.  Using $Z_{u,V,\Xi}(m/V)=m$ for all $m \in \Z$ we see that
$\wh u_\eps$ lies in $\rmW^{1,p}(a_{k-1}^\eps,a_k^\eps)$. Moreover, it
coincides with $\wh u$ on the points $a^\eps_k$ and thus we also have
$\wh u_\eps \in \rmW^{1,p}(a,b)$.

Because of the monotonicity of $Z_{U,V,\Xi}$ and $Z_{u,V,\Xi}(m/V)=m$
we have the obvious estimate $|Z_{U,V,\Xi}(t)-Vt|\leq 1$ which implies
$\|\wh u_\eps - \wh u\|_{\rmL^\infty}\leq \eps$. As we show below we
   have $\mfJ_\eps(\wh u_\eps,\wh\xi)\leq C$ for all $\eps\in
  {]0,1[}$. Hence the equi-coercivity of $\mfJ_\eps$ (cf.\
  \eqref{eq:Assum.R3}) yields $\| \dot{\wh u}_\eps\| \leq C$. Together
  with the uniform convergence, this implies $\wh u_\eps \weak \wh u$
  in $\rmW^{1,p}(0,T)$.  
 
\underline{Step 4: Limsup estimate}. 
We need to estimate the limsup of $\mfJ_\eps(\wh u_\eps,\wh \xi)$ from
above by $\wh J_0(\wh u,\wh \xi)$. Of course it suffices to do this in
the finitely many subintervals $[a,b]=[t_{j-1,N},t_{j,N}]$. 
We first observe that $\wh u$ is bounded and hence takes values in
$[-R,R]$ for a suitable $R$. Defining the piecewise constant
approximation $\ol u_\eps(t)=\wh u(a^\eps_k)$ for $t\in {[a^\eps_k,
  a^\eps_{k+1}[}$ our construction gives 
\[
 \|\wh u_\eps - \wh u\|_{\rmL^\infty}\leq \eps \quad \text{and} \quad 
 \|\ol u_\eps - \wh u\|_{\rmL^\infty} \leq 2\eps^{1/2}.
\]
Thus, using that $\mfJ_\eps$ is defined in terms of $N$ we have 
\begin{align*}
\mfJ_\eps(\wh u_\eps,\wh \xi)&= 
\int_a^b N(\Xi,\wh u_\eps(t),\inveps \wh u_\eps(t), \dot{\wh
  u}_\eps(t))\dd t \\
&=\sum_{k=1}^{n_\eps} \Big(\int_{a^\eps_{k-1}}^{x^\eps_k} N(\Xi,\wh
  u_\eps(t),\inveps \wh u_\eps(t), \dot{\wh   u}_\eps(t))\dd t + 
 \int_{x^\eps_k}^{a^\eps_k} N(\Xi,\wh
  u(t),\inveps \wh u(t), V)\dd t.
\end{align*}
We can now estimate from above by replacing $\wh u_\eps$ by the
interpolant $\ol u_\eps$ and can then use that $\wh u_\eps$ restricted
to $[a^\eps_{k-1},x^\eps_k]$ is exactly given by the optimal shape
functions $Z_{U^\eps_{k-1},V,\Xi}$\,. Using the uniform continuity
\eqref{eq:N.unif.cont} and $U^\eps_{k-1}\in\eps\Z$, we obtain the upper bounds 
\begin{align*}
\mfJ_\eps(\wh u_\eps,\wh \xi)&\leq \sum_{k=1}^{n_\eps} \Big(\!
\int\limits_{a^\eps_{k-1}}^{x^\eps_k}\!\!\!\! \big( N(\Xi,\ol u_\eps(t),\inveps \wh
u_\eps(t), \dot{\wh   u}_\eps(t)) {+} C\omega(\|\wh u_\eps{-}\ol
u_\eps\|_\infty)( 1{+}|\dot{\wh u}_\eps|^p) \big) \dd t +
C(|a^\eps_k{-}x^\eps_k|)\Big) 
\\
&= \sum_{k=1}^{n_\eps} \Big( (x^\eps_k{-}a^\eps_{k-1})
\big(\calM (U_{k-1}^\eps , V,\Xi) + C_V \omega(3\eps^{1/2}) \big) +
C\eps/V\Big),    
\end{align*}
where we used that $\dot{\wh u}_\eps(t)= \dot Z_{U^\eps,V,\Xi}$ is
bounded uniformly in $\rmL^p$ via $C_V=C(1{+}|V|^p)$, see Step 5 in
the proof of Proposition \ref{pr:N1cell}. 

Now, replacing the factor $(x^\eps_k{-}a^\eps_{k-1})$ by
$(a^\eps_k{-}a^\eps_{k-1})$, which is an error of $O(\eps)$ we find 
\begin{align*}
\limsup_{\eps\to 0} \mfJ_\eps(\wh u_\eps,\wh \xi)&\leq
\limsup_{\eps\to 0 } \int_a^b \calM (\ol u_\eps(t),V,\xi) = \int_a^b
\calM  \big(\wh u(t),\dot{\wh u}(t),\wh\xi(t)\big) = \mfJ_0(\wh u,\wh\xi),
\end{align*}
where we again used the continuity \eqref{eq:M0.contin} for $\calM $
and $\ol u_\eps \to \wh u$ in $\rmL^\infty(a,b)$. 
\end{proof}

In summary, we have now finished the proof of the main homogenization
results in Theorem \ref{th:GCvg.Jwigg}, which states the
\textGamma-convergence of $\mfJ_\eps$ to $\mfJ_0$ in the
weak$\ti$strong topology of $\rmW^{1,p}(0,T)\to \rmL^{p'}(0,T)$. 
The necessary liminf estimate is given in 
Proposition \ref{prop-liminf} and the existence of recovery sequences in
Proposition \ref{prop-limsup}.

\section{Properties of the effective contact potential $\calM $}
\label{se:M0.prop}

In this section we discuss the properties of $\calM $ that can be
derived directly from its definition in terms of the value function of
a minimization problem, see \eqref{eq:calM0}. In the rest of this
section, we drop the dependence on the variable $u$, because it is
simply playing the role of a fixed parameter. 

Moreover, we shortly
write $\mfp(y)= \pl_y \ppp(u,y)$, such that $\mfp:\R\to \R$ is an
arbitrary continuous and 1-periodic function with average $0$, viz.\
$\int_0^1 \mfp(y) \dd y=0$. We use the abbreviations
\[
\ol\mfp:= \max \set{\mfp(y)}{ y\in \R} \quad \text{and} \quad
\ul\mfp:= \min \set{\mfp(y)}{ y\in \R}.
\]  

\subsection{$\calM$ and the effective dissipation potential $\calR_\eff$}
\label{su:BasicEffective}

The first result concerns elementary properties that follow directly
from the fact that $\calM $ is defined in terms of the dual sum
$\calR+\calR^*$. 

\begin{lemma}[Basic properties of $\calM $]
\label{le:calM0}\mbox{}\\ 
(a) For all $v,\xi$ we have $\calM (v,\xi)\geq v\xi$. 
\\[0.3em]
(b)  For all $\xi\in \R$ we have 
\[
\calM (0,\xi)=\min\limits_{\pi\in[\ul\mfp,\ol\mfp]}
\calR^{*}(\xi{-}\pi) \quad \text{and } \calM (v,\xi)\geq
\calM (0,\xi) \text{ for all }v.
\] 
(c) If $\calR(-v)=\calR(v)$ for all $v$, then also
$\calM (-v,\xi)=\calM (v,\xi)$ for all $v,\xi\in \R$. If
additionally, $\mfp(y)=-\mfp(y_*{-}y)$ for some $y_*$ and all $y$,
then also $\calM (v,-\xi)=\calM (v,\xi)$.  
\end{lemma}
\begin{proof}
\underline{Part (a).}  For a minimizer $z$ for $\calM (v,\xi)$, we
simply apply the Young-Fenchel inequality  under the integration in the definition of $\calM $ and use that
$\mfp$ has average $0$: 
\begin{align*}
\calM (v,\xi)=\int_0^1 \big(\calR(|v|\dot
z){+}\calR^*(\xi{-}\mfp(z))\big) \dd s \geq \int_0^1 |v|\dot
z(s) (\xi{-}\mfp(z(s)))\big) \dd s = |v|(z(1){-}z(0)) \xi.
\end{align*}
Because of $z(1)=z(0)+\sign(v)$ we obtain the desired result. 

\underline{Part (b).} The result for $v=0$ is trivial, as we can
choose a constant minimizer $z(s)=z_*$. When comparing $v=0$ and
$v\neq 0$ we take a minimizer for $z_{v,\xi}$ and estimate
\[
\calM (v,\xi)=\int_0^1 \big(\calR(|v|\dot
z_{v,\xi}){+}\calR^*(\xi{-}\mfp(z_{v,\xi}))\big) \dd s
\geq \int_0^1 \min_{\pi\in [\ul\mfp,\ol\mfp]} \calR^*(\xi{-}\pi) \dd s
  = \calM (0,\xi). 
\]

\underline{Part (c).} The first symmetry follows since minimizers
$z_{v,\xi}$ give minimizers $z_{-v,\xi}:s \mapsto z_{v,\xi}(1{-}s)$
and vice versa. For the second symmetry we consider $z_{v,-\xi}:s
\mapsto -z_{v,\xi}(s)$.   
\end{proof}

 The next result concerns the most important point for our
effective contact potential $\calM $, namely the analysis of the
contact set 
\[
\mathsf C_{\calM }:=\bigset{(v,\xi)}{ \calM (v,\xi)=v\xi}
\]
We show that this set is the graph of the subdifferential of a
unique effective dissipation potential $\calR_\eff$.

\begin{proposition}[Effective dissipation potential] \label{pr:R.eff}
There is a unique dissipation potential $\calR_\eff:\R\to \R$ such
that 
\begin{equation}
  \label{eq:sfC.graphR}
  \mathsf C_{\calM }= \mafo{graph}(\pl \calR_\eff) = \bigset{(v,\xi)}{
    \xi\in \pl\calR_\eff(v)} = \bigset{(v,\xi)}{
    \calR_\eff(v)+\calR_\eff^*(\xi)=v\xi}. 
\end{equation}
If $\calR$ is strictly convex (and hence $\calR^*$ differentiable),
then the potential $\calR_\eff$ is characterized by the fact that
$\pl\calR^*_\eff(\xi)$ is the harmonic mean of the functions $[0,1]\in
y \mapsto \pl\calR^*(\xi{-}\mfp(y))$, viz.\ 
\[
\pl\calR_\eff^*(\xi)= \left\{ \ba{cl} 0& \text{for }\xi\in [\ul\mfp,\ol\mfp], \\
   K(\xi)&\text{for } \xi<\ul\mfp \text{ or } \xi>\ol\mfp, \ea \right.
\text{ where } 
K(\xi)=\bigg( \int_{0}^1 \frac {\dd y}{\pl\calR^*(\xi{-}\mfp(y))} 
 \bigg)^{-1}. 
\]
\end{proposition}
\begin{proof} 
In the proof of Lemma \ref{le:calM0}(a) we have seen that
$\calM (v,\xi)$ can only hold with equality if the minimizer
$z_{v,\xi}$ satisfies 
\[
\calR(|v|\dot z_{v,\xi}( s ))+
\calR^{*}(\xi-\mfp(z_{v,\xi}( s ))) = 
|v|\dot z_{v,\xi}( s )\,\big(\xi{-}\mfp(z_{v,\xi}( s ))\big) \text{ for a.a.\
} s\in [0,1]. 
\] 
By the Fenchel equivalences $z=z_{v,\xi}$ has to satisfy the
differential inclusion   
\begin{equation} \label{eq:Cont.Eq}
|v| \dot z(s) \in \pl\calR^* \big(\xi- 
\mfp (z(s)) \big), \quad z(1)=z(0)+\sign v. 
\end{equation}
If $\pl\calR^*$ is continuous, then we can solve the equation via
separation of the variables $z$ and $s$, and the boundary condition
gives 
\[
1= \int_{s=0}^1 \dd s = \int_{s=0}^1 \frac{|v|\dot z(s) \dd
  z}{\pl\calR^*\big(\xi{-}\mfp(z(s))\big)}= |v|\sign(v) \int_{y=0}^1
\frac{\dd y}{\pl\calR^*\big(\xi{-}\mfp(y)\big)}= \frac{v}{K(\xi)}.
\]  
Thus, the formula for $K$ is established. We observe that $\xi\mapsto
K(\xi)$ is monotone and $\xi K(\xi)\geq 0$. 
Hence, $\calR^*_\eff(\xi)=\int_0^\xi
K(\eta)\dd \eta$ gives the desired dual effective dissipation
potential. Defining $\calR_\eff $ by Legendre transform, the Fenchel
equivalences provide 
the desired relation between $\mathsf C_{\calM }$ and the graph of
$\calR_\eff$. 
\end{proof}

The explicit formula for $\pl\calR^*_\eff$ clearly shows how the
effective dissipation potential depends on the wiggly
$\mfp(y)=\pl_y\ppp(u,y)$. In particular, we obtain the sticking region
$\xi \in [\ul\mfp,\ol\mfp]$, where one has $v=0$. 
The special case $\calR(v)=\frac1{2\mu} v^2$ and $\mfp(y)=\hat a \sin(2\pi y)$
from \cite{Jame96HPT,AbChJa96KMWE} can be calculated explicitly, and we obtain 
\[
\pl\calR_\eff^*(\xi)=\mu \sign(\xi) \sqrt{\xi^2{-}\hat a{}^2}
\text{ for }\xi^2\geq \hat a^2 \quad \text{and} \quad \pl\calR_\eff^*(\xi)=0
  \text{ for } \xi^2\leq \hat a^2.
\]

\subsection{Expansions for $\calM $}
\label{su:ExpandM0}

We now want to study the behavior of $\calM (v,\xi)$ for small $v$,
which emphasizes the sticking phenomenon induced by the wiggly
energy landscape. To simplify the argument we assume that $\calR$
behave like a power near $v=0$. In fact, we restrict to the case $v>0$ 
by assuming
\begin{equation}
  \label{eq:calR.power}
  \calR(v)= \frac r\alpha v^\alpha + O(v^{\alpha+\delta}) \text{ for
  }v\searrow 0,
\end{equation}
where $\alpha>1$ and $r,\delta>0$. The proof involves an argument of
Modica-Mortola type (cf.\ \cite{ModMor77EGC} and
\cite[Ch.\,6]{Brai02GCB}) as for small velocities the minimizers $z$ for
$\calM $ are mostly near minimizers for $y\mapsto
\calR^*(\xi-\mfp(y))$ but have a transition layer of width $|v|$ to
make a jump of size 1.

\begin{lemma}[Expansion of $\calM $ for $v\approx
  0$]\label{le:Exp.M0}
Assume that in addition to all previous assumptions we also have
\eqref{eq:calR.power}, then for $v>0$ we have 
\begin{equation}
  \label{eq:calM.Expand}
  \calM (v,\xi)=M_0(\xi) \ + \
v \,M_1(\xi) + o(v) \text{ for }v\searrow 0,
\end{equation}
with $M_0(\xi)= \min_{\pi\in [\ul\mfp,\ol\mfp]} \calR^*(\xi{-}\pi)$
and $M_1(\xi)= \int_{y=0}^1
\Psi\big(\calR^*(\xi{-}\mfp(y)){-}M_0(\xi)\big) \dd y  $, where
$\Psi:{[0,\infty[}\to {[0,\infty[}$ is the inverse function of $\calR^*
:{[0,\infty[}\to {[0,\infty[}$. 

In particular, for $\xi\in  [\ul\mfp,\ol\mfp]$ we have $M_0(\xi)=0$
and if additionally $\calR$ is symmetric, then
$M_1(\xi)=\int_0^1|\xi{-}\mfp(y)|\dd y$.
\end{lemma} 
\begin{proof} We fix $\xi$ and choose $y_* \in \mafo{argmin}
  \calR^*(\xi{-}\mfp(\cdot))$. 
We rewrite $\calM (v,\xi)$ in the form 
\[
\calM (v,\xi)= \calM (0,\xi) + v M_1(v,\xi) \text{ with }
M_1(\xi,v)= \min_{z(1)=z(0)+1}  \int_0^1\frac1v\big(\calR(v\dot z) +
 G_\xi(z(y)) \big) \dd s ,
\] 
where $G_\xi(z)= \calR^*(\xi{-}\mfp(z))-\calR^*(\xi{-}\mfp(y_*))\geq
0$.   

Setting $s=v\tau$ and $w(\tau)=z(v\tau)$ we see that $w$ has to
minimize $\int_0^{1/v} \big(\calR(w'(\tau)+ G_\xi(w(\tau)\big) \dd
\tau$ under the constraint $w(1/v)=w(0)+1$. Indeed, by periodicity of
$\mfp$ in $y$ we may assume $w(0)=y_*$, so we are in the classical
Modica-Mortola setting of phase transitions. 

Our assumption \eqref{eq:calR.power} guarantees that $\calR^*$ is
strictly increasing for $\xi>0$, hence we can write $G_\xi(z)=
\calR^*(H_\xi(z))$ with $H_\xi(z)=\Psi(G_\xi(z))$. Now, the methods in  
\cite[Ch.\,6]{Brai02GCB} give the convergence
$M_1(v,\xi) \to M_1(0,\xi)$ with 
\[
M_1(0,\xi) =\min_{w(-\infty)=y_* \atop w(\infty)=y_*+1} 
 \int_{\tau\in \R}
 \big[\calR(w'(\tau)){+}\calR^*\big(H_\xi(w(\tau))\big) \big] \dd \tau
 = \int_{y_*}^{y_*+1} H_\xi(z) \dd z. 
\]
Because of the periodicity of $\mfp$ this is the desired formula for
$M_1$. 

The last statement follows if we use $\calR^*(-\xi)=\calR^*(\xi)$
which gives $\Psi(\calR^*(\eta))=|\eta|$. 
\end{proof} 

The formula for $M_1(\xi)$ can be made more explicit in the case of a
homogeneous potential $\calR(v)=\frac\nu p |v|^p$.  We have
$\calR^*(\eta)=  \frac1{p'}\nu^{-1/(p-1)} |\xi|^{p'}$ and  
$\Psi(\sigma)=  \nu^{1/p} (p'\sigma)^{1/p'}$.  

We finally look at the rate-independent limit that was already studied
in \cite{Miel12ERID}. The relevant time rescaling is obtained by
\[
\text{replacing } \  \calR \ \text{ by } \ \calR_\delta: v \mapsto
\frac1\delta \calR(\delta v),
\]
where $\delta$ is positive parameter that tends to $0$ in the
rate-independent limit, cf.\ \cite{EfeMie06RILS,MiRoSa09MSJR}. 

This scaling obviously gives $\calR^*_\delta(\xi)=\frac1\delta
\calR^*(\xi)$, so that the associated rescaled effective contact
potential is $\calM_\delta(v,\xi)=\frac1\delta \calM(\delta v,\xi)$. 
The above result provide the following convergence. We obtain indeed the same
result as in \cite[Prop.\,3.1]{Miel12ERID}, where a joint limit was
taken (i.e.\ $\delta_\eps\searrow 0$ with $\eps\searrow0$) while our
result is a double limit, where first $\eps\to 0$ and then $\delta\to
0$.  

\begin{corollary}[Rate-independent limit]\label{co:RateIndep}
Under the above assumptions including \eqref{eq:calR.power} and
$\calR(-v)=\calR(v)$  we have 
\[
\calM_\delta(v,\xi)\xrightarrow{\delta\to0} \calM_{\rmR\rmI}(v,\xi)= \left\{ \ba{cl} 
|v| M_1(\xi) &\text{for } \xi\in [\ul\mfp,\ol\mfp],\\ 
 \infty& \text{for } \xi \not\in [\ul\mfp,\ol\mfp], \ea
\right.
\text{ with }M_1(\xi)=\int_0^1 |\xi{-}\mfp(y)| \dd y. 
\]
\end{corollary}
\begin{proof} \underline{Case $\xi \not\in [\ul\mfp,\ol\mfp]$.} We have
  $\calM_\delta(v,\xi)\geq \calM_\delta(0,\xi)=\frac1\delta
  M_0(\xi)$. Because of $M_0(\xi)>0$ for this case we are done. 

\underline{Case $\xi\in [\ul\mfp,\ol\mfp]$.} We now have $M_0(\xi)=0$,
and Lemma \ref{le:Exp.M0} gives the result.  
\end{proof}

Finally we discuss the kinetic relation $v=\pl\calR^*_\eff(\xi)$ for
$\xi$ slightly outside the sticking region $[\ul\mfp , \ol\mfp]$ and
for very large $\xi$. For simplicity we restrict to the quadratic
case. 

\begin{lemma}[Expansion of kinetic relation]\label{le:Exp.KinRel} 
 Assume $\calR(v)=\frac12 v^2$ and let 
 $\mfp$ have a unique maximizer $z_*$ 
 such that  
 $\mfp(z) = \ol\mfp -c_{*}|z{-}z_{*}|^{\alpha} +
 O(|z{-}z_{*}|^{\gamma})$  with $c_*>0$, \  $1<\alpha<\infty$,
 and  $\gamma > 2\alpha -1$. Then, 
 \[
\calK(\xi) =
 c_{*}^{1/\alpha} S_\alpha^{-1}\max\{0,\xi{-}\ol\mfp\}^{(\alpha-1)/\alpha}
 + o(|\xi{-}\ol\mfp|^{(\alpha-1)/\alpha}) \quad\text{ for }\quad \xi
 \to \ol\mfp 
\]
with $S_\alpha = 2\sum_{0}^{\infty}(-1)^n \big(\frac{1}{\alpha n
  +1}+\frac{1}{\alpha(n+1)-1}\big)$. In the case $\alpha = 2$ we have
$S_2=\pi$ and $\calK(\xi) =
\sqrt{c_{*}}\pi^{-1}\big(\max\{0,\xi{-}\ol\mfp\}\big){}^{1/2}
+o(|\xi{-}\ol\mfp|^{1/2})$.  

For general $\mfp$ we obtain $\calK(\xi)
- \xi \rightarrow 0$ as $|\xi| \rightarrow \infty$
\end{lemma}
\begin{proof}
The computation is performed only on $[z_*,1]$ since we are able to conclude by symmetry. 
We define $h(z) = \mfp(z+z_*) - \ol\mfp + c_{*}|z|^{\alpha}$. With $\delta>0$ fixed we observe
\[
\int^{z^{*}+\delta}_{z^{*}} \frac{1}{\eps+\ol\mfp-\mfp(z)}\dd z =
\int^{\delta}_{0} \frac{1}{\eps+c_{*}z^\alpha} \dd z +
\int^{\delta}_{0}\frac{1}{\frac{(\eps+c_{*}z^\alpha)^2}{h(z)}+\eps+c_{*}z^\alpha}
\dd z.
\]

We want to argue only for the leading order term. Since $\gamma >
2\alpha - 1$ we have
\[0 \leq
\int^{\delta}_{0}\frac{1}{\frac{(\eps+c_{*}z^\alpha)^2}{h(z)}+\eps+c_{*}z^\alpha}
\dd z \leq \int^{\delta}_{0}\frac{h(z)}{(\eps+c_{*}z^\alpha)^2}
\dd z \leq \int^{\delta}_{0} c_{*}^{-1}h(z) z^{-2\alpha} \dd
z\rightarrow 0
\]
as $\delta\searrow0$. Let $\delta_\eps = (\eps/c_{*})^{1/\alpha}$. For
$h < x$ we use the geometric series
\begin{equation}
 \frac{1}{x+h} = \sum\limits_{n=0}^{\infty}
 (-1)^n\frac{h^n}{x^{n+1}}.\label{eq:tailor} 
\end{equation}
With this we  compute
\begin{align*}
 \int^{\delta_\eps}_{0} \frac{1}{\eps+c_{*}z^\alpha}\dd z &\overset{\eqref{eq:tailor}}{=} 
 \int^{\delta_\eps}_{0}\sum\limits_{n=0}^{\infty}(-1)^n\eps^{-(n+1)} (c_{*}z^{\alpha})^n \dd z\\
 &=\sum\limits_{n=0}^{\infty}(-1)^n\eps^{-(n+1)} \frac{c_{*}^n}{\alpha
   n +1}\delta_\eps^{\alpha n +1}\ 
 = \
 c_{*}^{-1/\alpha}\eps^{\frac{1}{\alpha}-1}\sum\limits_{n=0}^{\infty}(-1)^n
 \frac{1}{\alpha n +1} .
\end{align*}
 For the remaining interval we obtain 
\begin{align*}
 \int^{\delta}_{\delta_\eps} &\frac{1}{\eps+c_{*}z^\alpha}\dd z \overset{\eqref{eq:tailor}}{=} 
 \int^{\delta}_{\delta_\eps}\sum\limits_{n=0}^{\infty}(-1)^n\frac{\eps^{n}}{(c_{*}z^{\alpha})^{n+1}} \dd z\\
 &=\sum\limits_{n=0}^{\infty}(-1)^n \frac{\eps^n}{c_{*}^{n+1}}\frac{1}{\alpha(n+1)-1}
 \left(
 \left(\frac{c_{*}}{\eps}\right)^{n+1}\left(\frac{\eps}{c_{*}}\right)^{\frac{1}{\alpha}}
 -\frac{\delta}{\delta^{\alpha(n+1)}}
 \right)\\
 &=c_{*}^{-1/\alpha}\eps^{\frac{1}{\alpha}-1}\left(\sum\limits_{n=0}^{\infty}(-1)^n
   \frac{1}{\alpha (n +1)-1}\right) -\sum\limits_{n=0}^{\infty}(-1)^n
 \frac{1}{\alpha (n
   +1)-1}\left(\frac{\eps}{c_{*}\delta^\alpha}\right)^n 
  \frac{\delta^{1-\alpha}}{c_{*}} \,.
\end{align*}
We set $\delta=\eta_\eps$ such that $\eps = o(\eta_\eps)$
  and obtain
  $\int^{1}_{z_*+\eta_\eps} \frac{1}{\eps+\ol\mfp(u)-\mfp(z)}\dd z =
  o(\eps^{\frac{1}{\alpha}-1})$. This leads to 
\[
\calK(\xi) = c_{*}^{1/\alpha}
S_\alpha^{-1}(\max\{0,\xi{-}\ol\mfp\})^{1-\frac{1}{\alpha}}_{+} +
o(|\xi{-}\ol\mfp\}|^{1-\frac{1}{\alpha}})
\] 
with $S_\alpha$ as given above. 
For general $\mfp$ the limit $|\xi| \rightarrow \infty$ yields
\begin{align*}
 \calK(\xi)-\xi &= \xi\frac{1-\int_0^1\big(1-\frac{\mfp(z)}{\xi}\big)^{-1}\dd z}{\int_0^1\big(1-\frac{\mfp(z)}{\xi}\big)^{-1}\dd z}\\
 &=
 \frac{\int_0^1\mfp(z)\frac{\xi}{\mfp(z)}\big(1-\big(1-\frac{\mfp(z)}{\xi}\big)^{-1}\big)\dd
   z}{\int_0^1\big(1-\frac{\mfp(z)}{\xi}\big)^{-1}\dd z}\ 
 \longrightarrow \ \frac{\int_0^1-\mfp(z) \dd z}{\int_0^1 1 \dd z}=0.
\end{align*}
 This is the desired result.
\end{proof}

 Finally, we look at the case that the maximum of $\mfp$ is attained by
a linear approach, i.e.\ the limiting case $\alpha=1$ that is excluded
in the previous lemma. 

\begin{lemma}\label{le:1}  Assume $\calR(v)=\frac12 v^2$ and let
  $\mfp$ have a unique maximum such that $\mfp(z) = \ol\mfp
  -c_{*}|z {-} z_{*}| + O(|z{-}z_{*}|^\gamma)$ with $\gamma>1$, then
 \[
\calK(\xi) =
\frac{c_{*}}{2}\bigg(\log \big((\xi-\ol\mfp)^{-1}_{+}\big)\bigg)^{-1} +
o\bigg(\big(\log \big((\xi-\ol\mfp)^{-1}_{+}\big)\big)^{-1}\bigg) \quad\text{ as }\quad \xi\searrow \ol\mfp.
\]
\end{lemma}
\begin{proof}
  As in the proof of the previous lemma the computation is performed
  only on $[z_*,1]$ since we are able to conclude by symmetry.
  We define $h(z) = \mfp(z+z_*) - \ol\mfp + c_{*}|z|$. With
  $\delta>0$ fixed we observe
\[
\int^{z^{*}+\delta}_{z^{*}} \frac{1}{\eps+\ol\mfp-\mfp(z)}\dd z =
\int^{\delta}_{0}  \frac{1}{\eps+c_{*}z} \dd z +
\int^{\delta}_{0}\frac{1}{\frac{(\eps+c_{*}z)^2}{h(z)}+\eps+c_{*}z}
\dd z.
\]
We want to argue only for the leading order term. We have
\[
0 \leq
\int^{\delta}_{0}\frac{1}{\frac{(\eps+c_{*}z)^2}{h(z)}+\eps+c_{*}z}
\dd z \leq \int^{\delta}_{0}\frac{h(z)}{(\eps+c_{*}z)^2} \dd z
\leq \int^{\delta}_{0} c_{*}^{-1}h(z) z^{-2} \dd z\rightarrow 0
\]
as $\delta\searrow0$. For the remaining term we compute 
\[
\int^{\delta}_{0}  \frac{1}{\eps+c_{*}z} \dd z =
c_{*}^{-1}\big(\log(\inveps )+\log(\eps+c_{*}\delta)\big).
\]
We set $\delta=\eta_\eps$ such that $\eps =
  o(\eta_\eps)$. Then we have
  $\int^{1}_{z_*+\eta_\eps}\frac{1}{\eps+\ol\mfp-\mfp(z)}\dd z =
  o\big(\log(\inveps)\big)$.
\end{proof}

The following remark shows that $\pl\calR^*_\eff$ need not be
continuous.

\begin{remark} \label{re:kink} For $\mfp(z) = \ol\mfp
  -c_{*}|z{-}z_{*}|^{\alpha} + O(|z{-}z_{*}|^{\gamma})$ with $c_*>0$
  and $0<\alpha<1$ the integrand $z\mapsto (\xi {-} \mfp(z))^{-1}$
  remains integrable for $\xi \searrow \ol\mfp$, such that
  $\pl\calR^*_\eff(\xi)\to \sigma_*>0$. Hence, $\calR^*_\eff$ is
  Lipschitz continuous, but not differentiable, and $\pl\calR^*_\eff$
  is multi-valued, namely $\pl\calR^*_\eff(\ol\mfp)=[0,\sigma_*]$.
\end{remark}

\subsection{Lower and upper bounds on $\calR_\eff$}
\label{su:LowUppBouReff}

Here we provide a few bounds on $\calR_\eps$ and its Legendre dual
$\calR_\eff$ in terms of $\calR$, $\calR^*$, $\ul\mfp$, and
$\ol\mfp$. Throughout we restrict to the case $v\geq 0$ (and hence
$\xi\geq 0$, but similar
results hold for $v\leq 0$. 

The first result simply uses the fact that the harmonic mean can be
estimated from above and below by the maximum and the minimum,
respectively. 

\begin{proposition}[Bounds for $\calR_\eff$]\label{pr:Bound.Reff}
We always have the estimates
\begin{subequations}
\label{eq:Bou.Reff}
\begin{align}
\label{eq:Bou.R}
& \forall\, v\geq 0:\quad B^\mafo{low}_{\calR,\mfp}(v)  \leq \calR_\eff(v) \leq 
 \ol\mfp v + \calR(v),  \\
\label{eq:Bou.R*}
&\forall\, \xi\geq \ol\mfp: \quad \calR^*(\xi{-}\ol\mfp) \leq
\calR^*_\eff(\xi) \leq
\calR^*(\xi{-}\ul\mfp)-\calR^*(\ol\mfp{-}\ul\mfp)
\end{align}
\end{subequations}
where $B^\mafo{low}_{\calR,\mfp}(v)=\ol\mfp v $ for $v\in
[0,\pl\calR^*(\ol\mfp{-}\ul\mfp)]$ and $B^\mafo{low}_{\calR,\mfp}(v)=
\ul\mfp v{+}\calR(v){+}\calR^*(\ol\mfp{-}\ul\mfp)$ otherwise.%
\end{proposition}
\begin{proof}
  From $\ul\mfp \leq \mfp(y)\leq \ol\mfp$ we immediately obtain $
  \pl\calR^*(\xi{-}\ol\mfp) \leq \pl\calR^*_\eff(\xi) \leq
  \pl\calR^*(\xi{-}\ul\mfp) $ for all $\xi\geq \ol\mfp$.  Using
  $\calR_\eff^*(\xi)=0$ for $\xi\in [0,\ol\mfp]$ integration of
  these inequalities gives \eqref{eq:Bou.R*}.

For taking the Legendre transform, which is anti-monotone, in
\eqref{eq:Bou.R*} we have to extend the lower and upper bounds for
$\calR_\eff^*$ by $0$ on the interval $[0,\ol\mfp]$; then we obtain
\eqref{eq:Bou.R}. 
\end{proof}

Under additional assumptions these simple bounds can be improved. The
following result applies in particular to the case $\calR(v)=\frac
rp|v|^p$ with $p>1$, because 
${]0,\infty[}\ni v \mapsto 1/\pl\calR(v)= \frac1r v^{1-p}$ is convex. 

\begin{proposition}[Improved bound for $\calR_\eff$]\label{pr:ImpBouR}
Assume that the mapping ${]0,\infty[}\ni v \mapsto 1/\pl\calR(v)$ is
convex, then we have $\forall\, \xi\geq 0: \ \calR_\eff^*(\xi)\leq
\max\{0,\calR^*(\xi){-}\calR^*(\ol\mfp)\} $ or equivalently
\[
\calR_\eff(v) \geq \left\{ \ba{cl} \ol\mfp v& \text{for }v\in
  [0,\pl\calR^*(\ol\mfp)], \\ \calR^*(\ol\mfp) + \calR(v)& \text{for }
  v \geq \pl\calR^*(\ol\mfp). \ea \right. 
\]
\end{proposition} 
\begin{proof} Using the convexity of $1/\pl\calR^*$ we can apply
  Jensen's inequality and use $\int_0^1 \mfp(y)\dd y=0$. Thus, we
  obtain $\pl\calR^*_\eff(\xi)\leq \pl\calR^*(\xi)$ for all $\xi\geq
  \ol\mfp$. 

Integration gives the upper bound for $\calR^*_\eff$, and Legendre
transforms leads to the lower bound for $\calR_\eff$. 
\end{proof} 

In the case of the last result we obtain the simple bounds
$\calR_\eff^* \leq \calR^*$ and $\calR_\eff\geq \calR$. We expect that
these simple estimates hold in more general cases.\medskip

In the case of a $p$-homogeneous potential $\calR(v)=\frac rp|v|^p$ 
the dissipation $\pl\calR(v)v$ equals $p$ times the dissipation
potential, which is Euler's formula for homogeneous functions. For the
effective dissipation $\calR_\eff$ this homogeneity is destroyed, but
we still have a one-sided bound. 

Because $\pl\calR_\eff^*$ is defined as the harmonic mean of
$\pl\calR^*(\xi{-}\mfp(\cdot))$ we know that
$\pl\calR_\eff^*:{]\ol\mfp,\infty[}\to {[0,\infty[}$ is as smooth as
$\pl\calR^*$ and that $\pl\calR^*_\eff(\xi)=0$ for $\xi\in
{[0,\ol\mfp[}$. In general, there might be a kink at $\xi=\ol\mfp$,
see Remark \ref{re:kink}. For simplicity of the presentation we
restrict the following result to the case that $\calR^*_\eff$ is
differentiable.

\begin{proposition}[$p$-homogeneous case]\label{pr:pHom}
Assume that $\calR(v)=\frac rp|v|^p$ with $p>1$ and $r>0$ and
that $\calR^*_\eff$ is differentiable. Then we have 
\begin{equation}
  \label{eq:pHomReff}
  \pl\calR_\eff(v) v = \alpha(v) \calR_\eff(v)
\end{equation}
with a continuous function $\alpha:\R \to [1,p]$ satisfying
$\alpha(0)=1$ and $\alpha(v)\to p$ for $|v|\to \infty$. 
\end{proposition}
\begin{proof} Our proof uses the corresponding dual statement
  $\pl\calR^*_\eff(\xi) \xi = \beta(\xi) \calR^*_\eff(\xi)$ for $\xi
  \not\in [\ul\mfp , \ol\mfp]$.  It is  enough to consider the case $\xi
  >\ol\mfp$ as $\xi<\ul\mfp$ works analogously.
We relate $\alpha(v)$ and $\beta(\xi)$ for $\xi=\pl\calR_\eff(v)$
via
\[
\alpha(v)\calR_\eff(v)=\pl\calR_\eff(v)v=
\calR_\eff(v){+}\calR_\eff^*(\xi)=
\pl\calR^*_\eff(\xi)\xi=\beta(\xi)\calR^*_\eff(\xi).
\]
Hence, we have  $(\alpha(v){-}1)\,
(\beta(\xi){-}1)=1$, and it suffices to show that
$\beta:{]\ol\mfp,\infty[}\to {]p',\infty[}$ is continuous with
$\beta(\xi)\to \infty$ for $\xi\searrow \ol\mfp$ and $\beta(\xi)\to
p'$ for $\xi \to \infty$.   

From the convexity and differentiability of $\calR^*_\eff$ we
conclude that $\xi\mapsto \pl\calR^*_\eff(\xi)$ is even continuous. 
Thus, for $\xi>\ol\mfp$ the monotonicity of $\pl\calR^*_\eff$ gives
\[
\calR^*_\eff(\xi)= \int_{\ol\mfp}^\xi \pl\calR^*_\eff(\eta) \dd \eta
\leq (\xi{-}\ol\mfp) \pl\calR^*_\eff(\xi). 
\]
Hence for $\xi > \ol\mfp$ we have $\beta(\xi)=\xi\pl\calR^*_\eff(\xi)
/\calR^*_\eff(\xi) \geq \xi/(\xi{-}\ol\mfp)\to \infty$ for
$\xi\searrow \ol\mfp$ easily follows.  Moreover,
$\pl\calR^*_\eff(\xi)- \pl\calR^*(\xi)\to 0 $ for $\xi \to \infty$
implies $\beta(\xi)\to p'$. 

Thus, it remains to show $\beta(\xi)>p'$. For this it is sufficient to
show $\calH(\xi):=p'\calR_\eff^*(\xi)-\pl\calR_\eff^*(\xi)\xi <0$. The
continuity of $\pl\calR^*_\eff$ yields $\calH(\ol\mfp)=0$, and thus
the result follows from $\calH'(\xi)< 0$ for $\xi>\ol\mfp$. Using the
explicit form of $\pl\calR^*(\eta) = r_*\eta^{p'-1}$ for $\eta>0$ and
the definition of $\pl\calR_\eff^*$ in terms of the harmonic mean we
find
\begin{align*}
\calH'(\xi)&= (p' {-} 1) \pl\calR_\eff^*(\xi) -
\frac{\xi \int_0^1 (p' {-} 1) (\xi{-}\mfp)^{-p'} \dd y}{ 
     \big(\int_0^1  (\xi{-}\mfp)^{1-p'} \dd y\big)^2 }\\
& = (p' {-} 1)
\pl\calR_\eff^*(\xi) \Big( 1 - \frac{\int_0^1 h\dd y \int_0^1
  h^{-p'}\dd y }{\int_0^11\dd y\;\int_0^1 h^{1-p'}\dd y }\Big),
\end{align*}
where we set $h(y)=\xi-\mfp(y)>0$ (because of $\xi>\ol\mfp)$ and used
$\xi = \int_0^1 h(y)\dd y$ (because $\mfp$ has average $0$).   

We now estimate the denominator of the fraction in the right-hand side
by the numerator using suitable version of H\"older's inequality:
\begin{align*}
\int_0^1 1\dd y &= \int_0^1 h^{p'/(p'+1)} h^{-p'/(p'+1)}\dd y < 
\| h^{p'/(p'+1)}\|_{\rmL^{(p'+1)/p'}} \|
h^{-p'/(p'+1)}\|_{\rmL^{p'}}\\
& =\big(\int_0^1 h\dd y \big)^{p'/(p'+1)} \big(\int_0^1 h^{-p'} \dd
   y\big)^{1/(p'+1)},\\
\int_0^1 h^{1-p'}\dd y &= \int_0^1 h^{1/(p'+1)} h^{-p'{}^2(p'+1)}\dd y
<   \big(\int_0^1 h\dd y \big)^{1/(p'+1)} \big(\int_0^1 h^{-p'} \dd
   y\big)^{p'/(p'+1)}. 
\end{align*}
Here we have have strict inequality as $y\mapsto
h(y)=\xi{-}\mfp(y)$ is non-constant.  Multiplying these two estimates
we have established $\calH'(\xi)<0$, and the proof is complete.%
\end{proof}

\subsection{Convexity properties of $\calM$} 
\label{su:calM.convex}

In light of the Fitzpatrick functions considered in
\cite{Visi13VFSS,Visi15?SSFE,Visi17?EGCW} (see also Section
\ref{su:EGCWeakT}) and for the question about bipotentials in the
sense of \cite{BuDeVa08ECBG,BuDeVa08?NMCM} (see also Section
\ref{su:Bipot}) it is natural to ask what type of convexity properties
the function $(v,\xi)\mapsto \calM(v,\xi)$ has.

We first observe that $\calM$ cannot be convex in both variables,
if $\ppp$ is nontrivial. This
follows easily from the expansion $\calM(v,\xi)=M_0(\xi)+ vM_1(\xi) +
o(v)_{ v\searrow 0}$ obtained in Lemma \ref{le:Exp.M0}. As
$M_0(\xi)=0$ for $\xi\in [\ul\mfp,\ol\mfp]$ we see that for those
$\xi$ we have 
\[
\rmD^2 \calM(v,\xi) = \bma{cc} 0 & M'_1(\xi)\\ M'_1(\xi)& v
M''_1(\xi)\ema + o(1) \text{ for } v\searrow 0.
\]
This contradicts convexity because $\det\rmD^2
\calM(v,\xi)=-M'_1(\xi)^2 +o(1)_{v\searrow 0}<0$.

The next result states that $\calM(\cdot,\xi)$ is always convex.

\begin{proposition}[Convexity of $\calM(\cdot,\xi)$] \label{pr:BiPot} 
For all $\xi\in \R$ the function $\calM (\cdot,\xi):\R\to \R$ is convex.
\end{proposition}
\begin{proof} 
  This convexity was already established in Proposition
  \ref{pr:N1cell}(C). For completeness we give a second and
  independent proof.

To show convexity of  $\calM (u,\cdot,\xi)$ we recall that Theorem
\ref{th:GCvg.Jwigg} states that  
$\mfJ_0:(u,\xi)\mapsto \int_0^T \calM(u,\dot u,\xi)\dd t$ is the
\textGamma-limit of $\mfJ_\eps$ in the weak$\ti$strong topology of
$\rmW^{1,p}(0,T)\ti \rmL^{p'}(0,T)$. The standard theory of
\textGamma-convergence \cite{Dalm93IGC,Brai02GCB} 
now implies that $\mfJ_0$ is lower
semicontinuous. In particular $v \mapsto \int_0^T \calM(u,v,\xi)\dd t
$ must be weakly lower semicontinuous in $\rmL^{p}(0,T)$, which
implies that $\calM(u,\cdot,\xi)$ must be convex. 
\end{proof}

 We now turn to the question of convexity of $\xi \mapsto
\calM(v,\xi)$ for fixed $v\in \R$. For this, we start from the
functionals
\[
\calN_{v,\xi}(z):=\int_{s=0}^1 \Big(\calR(|v|\dot z(s))+
\calR^*\big(\xi{-}\mfp(z(s))\big) \Big) \dd s, 
\]
then $\calM(v,\xi) = \inf\set{\calN_{v,\xi}(z(\cdot)) }{z \in
  \rmW^{1,p}_v}$. 

To study the convexity of $\calM (v,\cdot)$ 
we derive a characterization, which is the basis of the subsequent
analysis. The main idea is to invert for the minimizer $z_{v,\xi}$ of 
$\calN_{v,\xi}$ the relation $y=z_{v,\xi}(s)$ into $s=S_{v,\xi}(y)$, which
transforms the nonlinear function $y \mapsto \mfp(y)$ into a
non-constant coefficient. The new functional is then convex in the
unknown functions $S:y\mapsto S(y)$.  For this we define the two 
convex functions $\psi_+,\,\psi_-:\R\to [0,\infty]$ via
\[
\psi_\pm:\rho\mapsto \left\{ \ba{cl} |\rho|\calR(1/\rho)&\text{for
  } \pm \rho>0,\\ \infty &\text{for }\pm\rho<0, \ea\right.
\]
where the value at $\rho=0$ is fixed by lower semicontinuity. 
For simplicity, we consider subsequently the case $v>0$ only and
write $\psi=\psi_+$. The case $v<0$ can be done
similarly by using $\psi_-$. By \eqref{eq:Assum.R1} we have
$\pl\calR(0)=0$ which implies $\psi(\rho)\to 0$ for $
\rho\to \infty$. For $\rho\approx$ we have $\psi(\rho)\geq
c_1\rho^{1-p} -c_1\rho$, i.e.\ $\psi$ blows up at $\rho=0$. 

We now recall the representation of $\calM(v,\xi)$ introduced in
Proposition \ref{pr:N1cell}, see \eqref{eq:HomN1cell.b}, which is the
basis for our subsequent convexity discussion. Defining the functional
\[
\calT_{v,\xi}(a):= \int_{0}^1\!\!
\Big(\calR\big(\frac{v}{a(y)}\big) +
\calR^*( \xi{-}\mfp(y)) \Big) a(y) \dd y = \int_{0}^1 \!\!
\Big(v\,\psi\big(\frac{a(y)}v \big) +
a(y)\calR^*( \xi{-}\mfp(y)) \Big)  \dd y 
\]
we can express $\calM(v,\xi)$ for $v>0$ in the form
\begin{equation}
  \label{eq:M.altern}
\calM (v,\xi) = \inf\bigset{ \calT_{v,\xi}( a)}{ a > 0 \text{ and } \int_{y=0}^1
  a(y) \dd y =1}.   
\end{equation}
It is not difficult to show $\calT_{v,\xi}$ admits a minimizer
$a=A_{v,\xi}$, which is unique by the strict convexity of
$\calT_{v,\xi}$. Moreover \eqref{eq:Assum.R3} implies $\psi(\rho)\geq
c\rho^{1-p}$ for small $\rho$, so $A_{v,\xi}$ is bounded from below by
a positive constant. The point now is that the minimizer $A_{v,\xi}$
can be obtained almost explicitly, since the Euler--Lagrange equations
are given by
\begin{equation}
  \label{eq:EL.calT}
  \psi'(a(z)/v) + \calR^*(\xi{-}\mfp(z))= \HH,
\end{equation}
where the constant Lagrange multiplier $\HH$ associated with the constraint
$\int_0^1 a\dd z =1$ has to be chosen as a function of $(v,\xi)$ such
that $a$ satisfies the constraint, namely $\HH=\HHhh(v,\xi)$. 

For this we use the Legendre transform $\psi^*:{]{-}\infty,0]}\to
{[0,\infty]}  $ of $\psi=\psi_+$ given by
\[
\psi^*(\sigma)=\infty \text{ for } \sigma>0 \quad \text{and} \quad 
\psi^*(\sigma)=\psi_*(\sigma):= \sup\set{\sigma s
  - \psi(s)}{s>0} \text{ for } \sigma<0. 
\]
With this we have
\begin{equation}
  \label{eq:A.v.xi}
  a=A_{v,\xi}(z)= v\, \psi'_*\big( \HHhh(v,\xi) {-} 
\calR^*(\xi{-}\mfp(z))\big).
\end{equation}
Thus, the value $\HH=\HHhh(v,\xi)$ is determined by solving
\[
1 = v \int_0^1 \psi'_*\big( \HH {-} G(\xi,z)\big) \dd z \quad \text{with
  } G(\xi,z):= \calR^*(\xi{-}\mfp(z)).
\]
Note that $\psi_*(\sigma)$ is only defined for $\sigma= \HH-G(\xi,z)\leq
0$. Thus, we always assume 
\[
\HH < \inf\set{G(\xi,z) }{ z\in [0,1]}.
\]
Because of $G(\xi,z)\geq 0$ the case $\HH<0$ is always admissible, while
$\HH\geq 0$ can only be allowed when $\xi $ lies outside
$[\ul\mfp, \ol\mfp]$. 

It is now advantageous to introduce the functional 
\[
\calW(\xi,\HH):= \int_{0}^1\!\! w(\xi,\HH,z) \dd z \text{ with }
w(\xi,\HH,z):= \psi_*\big(\HH-G(\xi,z)\big).
\]
The following formulas for the partial derivatives of $\calW$ are
immediate when after interchanging integration with respect to $z\in
[0,1]$ and differentiations. 
\begin{align*}
&\calW_\HH(\xi,\HH)=\int_0^1 \psi'_*(\HH{-}G)\dd z>0, \qquad 
 \calW_\xi(\xi,\HH)= -\int_0^1\psi'_*(\HH{-}G) G_\xi \dd z ,\\
& \calW_{\HH\HH}(\xi,\HH)=\int_0^1 \psi''_*(\HH{-}G)\dd z>0, \quad 
 \calW_{\xi \HH}(\xi,\HH)= -\int_0^1 \psi''_*(\HH{-}G)G_\xi\dd z>0, \\
&\calW_{\xi \xi}(\xi,\HH)=\int_0^1 \Big(  \psi''_*(\HH{-}G) G^2_\xi - 
\psi'_*(\HH{-}G) G_{\xi\xi} \Big)\dd z. 
\end{align*}
Thus, $\HH=\HHhh(v,\xi)$ is obtained by solving the equation 
$1 = v \calW_\HH(\xi,\HH)$. 

\begin{remark}[Involution property]\slshape
\label{re:Involution} In fact, we may evaluate $\calW$ for $\HH=0$
explicity, since
\[
w(\xi,0,y)= - |\xi{-}\mfp(y))| \text{ for }(\xi,y)\in \R\times [0,1].
\]
For this we use the relation $\psi_*(-\calR^*(\eta))=-\eta$ for all
$\eta\in\R$, which holds under the additional evenness assumption
$\calR^*(-\eta)=\calR^*(\eta)$ (see \cite[Sec.\,4.2,
eqn.\,(4.9)]{LiMiSa14?OETP} for a proof).  
Hence, we obtain $\calW(\xi,0)= -\int_0^1 | \xi {-} \mfp(y)| \dd y $,
which immediately implies that $\calW(\cdot,0)$ is concave. Moreover,
for $\xi \not\in
\mathrm{range}(\mfp)=[\ul\mfp,\ol\mfp]$ we obtain
$\calW(\xi,0)= -|\xi|$ because of $\int_0^1 \mfp(z) \dd z=0$.   

Note that $\HH=0$ corresponds via \eqref{eq:EL.calT} and the definition
of $\psi$ and $a=S'_{v,\xi}=1/Z'_{v,\xi}$ to the equation $
\calR(vz') - vz'\calR'(vz')+\calR^*(\xi{-}\mfp(z))=0$. Using Fenchel's
equivalence this implies the pointwise contact relation
\[
\calR(vz'(s)) + \calR^*(\xi{-}\mfp(z(s))) =
vz'(s)\big(\xi{-}\mfp(z(s))\big) 
\]  
as established for $(v,\xi) \in \mathsf C_\calM$, see \eqref{eq:Cont.Eq}.
\end{remark}

The following identities are useful in the sequel.

\begin{lemma}[Identities connecting $\calW$ and $\calM$]
\label{le:M.vs.calW} \mbox{}

(A) $\calM (v,\xi)= \big(\HH - v\calW(\xi,\HH)\big)|_{\HH=\HHhh(v,\xi)}$;

(B) $\HHhh_v(v,\xi) = -\calW^2_\HH/\calW_{\HH\HH} \big|_{\HH=\HHhh(v,\xi)} 
\ \text{ and } \ 
\HHhh_\xi(v,\xi) = -\calW_{\xi \HH}/\calW_{\HH\HH} \big|_{\HH=\HHhh(v,\xi)}$;

(C) $\calM_v(v,\xi)=-\calW(\xi,\HHhh(v,\xi)) $, \ 
    $\calM_\xi(v,\xi)=- v \calW_\xi(\xi,\HHhh(v,\xi))$; 

(D)  $\calM_{vv}(v,\xi) = \calW_\HH^2/\calW_{\HH\HH}\big|_{\HH=\HHhh(v,\xi)} >0$, 
  \  $\calM_{v\xi}(v,\xi)=- \calW_\xi + \calW_\HH \calW_{\xi
  \HH}/\calW_{\HH\HH} \big|_{\HH=\HHhh(v,\xi)}$,\\
 \mbox{}\qquad \quad  $\calM_{\xi\xi}(v,\xi)=\frac{v}{\calW_{\HH\HH}}
 \big( \calW_{\xi \HH}^2 -\calW_{\HH\HH} \calW_{\xi\xi}\big)\big|_{\HH=\HHhh(v,\xi)}$. 
\end{lemma}
\begin{proof}
ad (A): Fenchel-equivalence means that $s=\psi'_*(\sigma) $ holds 
 if and only if $\psi(s)+\psi_*(\sigma)=s\sigma$. Thus, we have 
\[
\psi\big(\psi'_*(\sigma)\big) = \sigma \psi'_*(\sigma) - \psi_*(\sigma),
\]
We use this for $\sigma= \HH-G$ when inserting
the minimizer $a=A_{v,\xi}$ from
\eqref{eq:A.v.xi} into $\calT$ to obtain 
\begin{align*}
& \calM (v,\xi)\ = \ \calT(v,\xi;A_{v,\xi}) \ = \  \int_0^1 \!\!\Big(
v\psi(\psi'_*(\sigma(z))) + v\psi'_*(\sigma(z)) G(\xi,z)\Big) \dd z\\
&= v \int_0^1\!\!\Big( (\HH{-}G) \psi'_*(\HH{-}G) - \psi_*(\HH{-}G)
+G\psi'_*(\HH{-}G)  \Big) \dd z  \
 =\ \big( \HH- v \calW(\xi,\HH)\big)\big|_{\HH=\HHhh(v,\xi)}. 
\end{align*}
For the first derivatives of $\calM$ we use the implicit function
theorem on $1= v \calW_\HH(\xi,\HHhh(v,\xi))$ and obtain (B).
Now using the relations (B) and (C) the chain rule provides the
relations (D). 
\end{proof}

As $\calW_{\HH\HH}$ is positive, the convexity of $v \mapsto \calM
(v,\xi)$ follows for arbitrary $\xi\in \R$. For the convexity of $\xi
\mapsto \calM (v,\xi)$ we need to show that
\begin{equation}
  \label{eq:W.2ndD}
\calW_{\xi \HH}(\xi,\HH)^2 \geq \calW_{\HH\HH}(\xi,\HH)\calW_{\xi\xi}(\xi,\HH)   
\end{equation}
for all relevant $\xi$ and $\HH$. We see that this is not always
the case. However, we have a positive result if $\calR$ is
$p$-homogeneous, because in this case also $\psi_*$ is of power-law type and a 
nontrivial cancellation takes place.

\begin{theorem}[Convexity of $\calM (v,\cdot)$]\label{th:M.cvx.xi}
Assume $\calR(v)=r|v|^p$ for $p>1$ and $r>0$. Then for all $v\in \R$
the function $\calM (v,\cdot):\R\to \R$ is convex. 
\end{theorem}
\begin{proof}
It is sufficient to show \eqref{eq:W.2ndD}. 
To this end we note that the assumptions imply
\[
\calR^*(\eta)=r_*\eta^{1/a} \ \text{ and }
\psi_*(\sigma)=-f_*(-\sigma)^a 
\]
where $a=1-1/p\in {]0,1[}$.  By the homogeneity of \eqref{eq:W.2ndD}
we may assume $r_*=f_*=1$ for simplicity.  We establish the
desired inequality in two steps, one for $\HH\leq 0$ and one for
$0<\HH\leq \min G(\xi,\cdot)$ with quite different arguments.

\underline{Step 1:} $\calW_{\xi\xi}(\xi,\HH)\leq 0$ for $\HH\leq 0$. 
\\
We use $\calW_{\xi\xi}(\xi,\HH)=\int_0^1 w_{\xi\xi}(\xi,\HH,z) \dd z$
with $w_{\xi\xi}(\xi,\HH,z) = \psi''_*(\HH{-}G) G_\xi^2 -
\psi'_*(\HH{-}G)G_{\xi\xi}$. 
The power-law structure of $\calR^*$ easily gives the identity 
\[
(1{-}a)G_\xi^2 = G G_{\xi\xi}= \HH G_{\xi\xi}-(\HH{-}G) G_{\xi\xi}.
\]
Similarly, the power-law structure of $\psi_*$ gives 
\[
(1{-}a) \psi'_*(\HH{-}G) = (G{-}\HH)\psi''_*(\HH{-}G).
\]
Using these two relations we can simplify $w_{\xi\xi}$ and find 
\begin{equation}
  \label{eq:w.xi.xi}
  w_{\xi\xi}(\xi,\HH,z)= \psi''_*(\HH{-}G)\frac{G_{\xi\xi}}{1{-}a}\Big( 
G - \big(G{-}\HH\big) \Big) =\psi''_*(\HH{-}G)\frac{G_{\xi\xi}}{1{-}a} \HH
  .
\end{equation}
With $a<1$, $\psi''_*, G_{\xi\xi}\geq0$ we conclude $w_{\xi\xi}\leq
0$, and by integration of a non-positive function we obtain
$\calW_{\xi\xi} \leq 0$, and \eqref{eq:W.2ndD} trivially holds
because of $\calW_{\HH\HH}\geq 0$.

\underline{Step 2.} For $\HH>0$ we establish the estimate by showing
\begin{equation}
  \label{eq:W.xi.E}
\text{(a) \ } | \calW_{\xi \HH}| \geq \frac{\HH^{1-a}}{a} \calW_{\HH\HH} 
\quad \text{and} \quad \text{ (b) \ }
| \calW_{\xi \HH}|\geq \frac{a}{\HH^{1-a}} \calW_{\xi\xi}.
\end{equation}
The major observation for $\HH>0$ is that $G(\xi,z)=\calR^*(\xi{-}\mfp(z))
= |\xi{-}\mfp(z)|^{1/a} \geq \HH>0$ implies 
\[
|G_\xi(\xi,z)|=  \frac1a |\xi-\mfp(z)|^{(1-a)/a} \geq \HH^{1-a}/a. 
\] 
In particular, the continuous function $z\mapsto G_\xi(\xi,z)$ cannot
change the sign. Thus, we conclude
\begin{align}
 \label{eq:W.mix.Sign}
|\calW_{\xi \HH}|&= \Big| \int_0^1 \!\!G_\xi \psi''_*(\HH{-}G)\dd z \Big| = 
 \int_0^1 \! \big|G_\xi\big| \psi''_*(\HH{-}G)\dd z \\
 \nonumber 
& \geq 
  \int_0^1 \frac{ \HH^{1-a}}{a} \psi''_*(\HH{-}G)\dd z= \frac{ \HH^{1-a}}{a}
  \calW_{\HH\HH}>0.
\end{align}
Thus, \eqref{eq:W.xi.E}(a) is established. 

For part (b) we can use relation \eqref{eq:w.xi.xi}, which obviously
also holds for $0<\HH\leq \min G(\xi,\cdot)$. With  
\[
|G_{\xi}(\xi,z)|=\frac{1}{a} |\xi{-}\mfp(z)|^{(1-a)/a} = \frac a{1{-}a}
|\xi{-}\mfp(z)|\, G_{\xi\xi}(\xi,z) \geq \frac{a \HH^a}{1{-}a} \,G_{\xi\xi}(\xi,z)
\]
we find $|w_{\xi \HH}| =|G_\xi |\psi''_*(\HH{-}G) \geq \frac{a \HH^a}{1{-}a}
\, \psi''_*(\HH{-}G) G_{\xi\xi}(\xi,z)= a\HH^{a-1} w_{\xi\xi}$. Again
using \eqref{eq:W.mix.Sign} we can integrate this estimate, which
yields \eqref{eq:W.xi.E}(b).  

Multiplying the two estimates in  \eqref{eq:W.xi.E} finishes the proof
of \eqref{eq:W.2ndD} in the case $\HH>0$. Exploiting the last relation
in assertion (D) of Lemma \ref{le:M.vs.calW} provides the desired
convexity of $\xi \mapsto \calM (v,\xi)$.  
\end{proof}  

We conclude this subsection by showing that for general dissipation
potentials $\calR^*$ we cannot expect to have convexity for
$\calM (v,\cdot)$. A counterexample can be constructed by
exploiting part (D) in Lemma \ref{le:M.vs.calW} for an even function
$\calW(\cdot,\HH)$, then in addition to the obvious relation 
$\calW_{\HH\HH}>0$ we have $\calW_{\xi \HH}(0,\HH)=0$ and hence it suffices to
show $\calW_{\xi\xi}(0,\HH)>0$ for some $\HH$. Based on \eqref{eq:w.xi.xi}
it suffices to choose $G(\xi,z)=\calR^*(\xi{-}\mfp(z))$ having a small
second derivative $G_{\xi\xi}$ while $G_\xi$ is large.

\begin{example}[$\calM (v,\cdot)$ may be nonconvex]\label{ex:M.noncvx}\slshape
For a simple counterexample we consider the case that $\mfp(z)=\pm 2$
for $z\in [0,1/2]$ and $z\in 
{]1/2,1[}$ respectively. Continuity can be restored in very small
layers that don't destroy the non-convexity generated below.
 
Moreover, we only consider $|\xi|\leq 1$, since non-convexity occurs
near $\xi=0$. Thus, the relevant values of $\eta=\xi-\mfp(z)$ satisfy 
$|\eta|=|\xi{-}\mfp(z)|\in [1,3]$. 

The dual dissipation potential is chosen as 
\[
\calR^*(\eta) = \left\{ \ba{cl}  \eta^2& 
                   \text{for }|\eta|\leq 1,\\ 
2|\eta|{-}1    &\text{for } |\eta| \in [1,3],\\
21-8\sqrt{7{-}|\eta|\,} &\text{for } |\eta| \in [3,6],\\ 
\text{convex extension}&\text{for }|\eta|\geq 6.   \ea
\right.
\] 
We find the potential $\psi_*$ in the form
\[
\psi_*(\sigma) = \left\{ \ba{cl} \infty & \text{for } \sigma>0, \\
 -\sqrt{-\sigma} &\text{for } \sigma\in [-1,0], \\
  (\sigma{-}1)/2& \text{for }  \sigma\in [-5,-1], \\ 
  (\sigma^2{+}42\sigma {-}7)/64& \text{for }  \sigma\in [-13,-5], \\ 
 \text{convex extension} & \text{for }\sigma \leq -13.\ea \right.
\]

Thus, we can express the function $w(\xi,\HH,z)$ explicitly in a
certain range of $(\xi,\HH)$, because
integration over $z\in [0,1]$ only leads to two different values
$\mfp(z)=\pm 2$:
\begin{align*}
\calW(\xi,\HH) &= \frac12\Big( \psi_*\big(\HH - 2|\xi{+}2|+1\big)  +
               \psi_*\big(\HH - 2|\xi{-}2|+1 \big) \Big)\\
 & = \frac12\Big( \psi_*\big(\HH -3 -2\xi\big)  +
               \psi_*\big(\HH - 3 +2\xi  \big) \Big)\ = \ 
\frac{\HH^2+36 \HH -124}{64} + \frac{\xi^2}{16},
\end{align*}
where we used $|\xi|\leq 1$ for the first identity and $\HH\in
[-8,-4]$.  Thus, we have $\calW_{\xi \HH}\equiv 0$, $\calW_{\HH\HH}=
1/32>0$ and $\calW_{\xi\xi}=1/8 >0$, which implies $\calW_{\xi \HH}^2 -
\calW_{\HH\HH} \calW_{\xi\xi} \equiv -1/256$ for $|\xi|\leq 1$ and $\HH\in
[-8,-4]$. 

We can even solve $v\calW_\HH(\xi,\HH)=1$ and calculate $\calM (v,\xi)$
explicitly according to Lemma \ref{le:M.vs.calW}(A). First we find $
\HH=\HHhh(v,\xi)= {32}/v - 18 $ and obtain
\[
\calM (v,\xi)= \frac{16}v - 18 + v\big(7-\frac{\xi^2}{16}\big) 
\quad \text{for }(v,\xi) \in \big[ 
\frac{32}{14}, \frac{32}{10}\big]\ti[-1,1].    
\]
Thus, the concavity of $\calM (v,\cdot)$ on $[-1,1]$ is seen explicitly
because of $v\geq 32/14$.
\end{example}

\subsection{Bipotential-property of the limiting dissipation}
\label{su:Bipot}

In this section we consider the question whether the functional
\[
(v,\xi) \mapsto \calM (v,\xi)
\]
defined in \eqref{eq:calM0} is a \emph{bipotential} in the sense of 
\cite{BuDeVa08ECBG,BuDeVa08?NMCM}, see also \cite[Sec.\,3.1]{MiRoSa12BVSV}
and \cite[Sec.\,3.1]{MiRoSa16BVSI}, where they are
also called \emph{contact potentials}. For a reflexive 
Banach space $X$ with dual space
$X^*$ a function $B:X\ti X^*\to \R_\infty$ is called
\emph{bipotential} if it satisfies the following three conditions:
\begin{subequations}
\label{eq:Bipot}
\begin{align}
\label{eq:Bipot.a}& \forall\, v\in X\ \forall\, \xi\in X^*: \quad 
B(v,\cdot):X^*\to \R_\infty \text{ and } B(\cdot,\xi):X\to \R_\infty
\text{ are convex}, \\
\label{eq:Bipot.b}& \forall\, v\in X\ \forall\, \xi\in X^*: \quad
B(v,\xi)\geq \langle \xi,v\rangle,\\
\label{eq:Bipot.c}& 
\forall\, \wh v\in X\ \forall\, \wh\xi\in X^*:\quad \wh\xi \in \pl_v B( \wh v, \wh\xi) \ \Longleftrightarrow \ 
 \wh v \in \pl_\xi B( \wh v, \wh\xi) \ \Longleftrightarrow 
\ B( \wh v, \wh\xi)=\langle   \wh\xi, \wh v\rangle. 
\end{align}
\end{subequations}

 Under quite general assumptions one can show that effective
contact potentials $\calM(q,\cdot,\cdot):\bfQ\ti\bfQ^*\to \R $
satisfy the convexity of $\calM(q,\cdot,\xi)$ and the estimate
$\calM(q,v,\xi)\geq \langle \xi, v\rangle $. Hence, we can expect the
weaker property
\[
\calM(q,v,\xi) = \langle\xi, v\rangle  \quad \Longleftrightarrow \quad
\xi \in \pl_v \calM(v,\xi).
\]
(For this it is sufficient that for all $\xi$ there is at least one $v$
such that $\calM(q,v,\xi)=\xi v$.) In that case we can use the
energy-dissipation principle starting from the derived
energy-dissipation balance 
\[
\calE_0(T,q(T)) + \int_0^T\calM\big(q,\dot
q,{-}\rmD\calE_0(t,q)\big) \dd t \leq
\calE(0,q(0))+\int_0^T\pl_t\calE_0(t,q) \dd t
\]
(by involving a suitable chain-rule inequality) to obtain the
subdifferential inclusion 
\[
 0 \in \pl_v \calM(q,\dot q,{-}\rmD\calE_0(t,q)) + \rmD\calE_0(t,q).
\]
The disadvantage of such a formulation is that $\rmD\calE_0$ appears
twice and the dependence on $\pl_v\calM(v,\xi)$ on $\xi$ is difficult
to control in general cases. If $\calM$ is even a bipotential, one
also has the inverted equation 
\[
\dot q \in \pl_\xi\calM(q,\dot q,{-}\rmD\calE_0(t,q)),
\]
where now $\dot q$ shows up twice. These forms are not easy to handle,
but they allow for new applications, e.g.\ in the mechanics of
friction or soil mechanics, see \cite{BuDeVa08ECBG, BuDeVa08?NMCM,
  BulDes17SBEN}.
 
It is exactly the key ingredient of our notion of relaxed
EDP-convergence, that we asked that our effective contact
potential $\calM$ is such that the conditions in \eqref{eq:Bipot.c}
are in fact equivalent to the corresponding relations for $\calM_\eff:
(v,\xi)=\calR_\eff(v)+ \calR^*_\eff(\xi)$.  Nevertheless it is
interesting to check whether $\calM$ is indeed a bipotential. 

In the previous subsection we have analyzed the question of separate
convexity for $\calM$, i.e.\ convexity of $v\mapsto \calM(v,\xi)$ and
$\xi \mapsto \calM(v,\xi)$. We have seen that the first convexity
always holds, while the second is false in general. So we cannot
expect $\calM$ to be a bipotential without assuming further
properties. The following result shows that in the case
$\calR(v)=\frac rp|v|^p$ we have indeed a
bipotential.

\begin{theorem}[Bipotential property] \label{th:Bipotential} 
Assume that $\calR(v)=\frac rp|v|^p$ with $r>0$ and $p>1$, then for
all 1-periodic $\mfp\in \rmC^0(\R)$ with average $0$, the effective
contact potential $\calM$ is a bipotential, i.e.\ \eqref{eq:Bipot}
holds. 
\end{theorem} 
\begin{proof} \underline{Step 1: First two conditions}
Obviously, the conditions \eqref{eq:Bipot.a} and \eqref{eq:Bipot.b}
are satisfied for $B=\calM$, see Proposition \ref{pr:BiPot}, Theorem \ref{th:M.cvx.xi}, and  Lemma
\ref{le:calM0}(a). 

\underline{Step 2: Condition 3 ``$\Rightarrow$''.} 
It remains to establish third condition
\eqref{eq:Bipot.c}, which reads here 
\begin{equation}
  \label{eq:M.3Cond}
  \xi\in \pl_v \calM (v,\xi) \ \Longleftrightarrow \   
  \calM (v,\xi)=\xi v \ \Longleftrightarrow \ v\in \pl_\xi  \calM (v,\xi).
\end{equation}
Of course, \eqref{eq:Bipot.b} for $B=\calM$ immediately gives the implication 
\[
\calM(v,\xi)=v\xi \ \Longrightarrow \ \big( \xi\in \pl_v \calM (v,\xi)
\text{ and } v\in \pl_\xi  \calM (v,\xi)\big).
\]
It remains to show that the two outer relations in \eqref{eq:M.3Cond}
are equivalent to the middle relation.  \medskip

\underline{Step 3: The case $v=0$.}  In this case the first condition
in \eqref{eq:M.3Cond} reads $\xi\in \pl_v\calM(0,\xi)$. Our expansion 
\eqref{eq:calM.Expand} gives $\pl_v M(0,\xi)= [-M_1(\xi),M_1(\xi)]$.
Because of $\calR^*(\eta)=c_*|\eta|^{p'}$ we have 
\[
M_1(\xi)=\left\{\ba{cl}
    \int_0^1( |\xi{-}\mfp(y)|^{p'}{-}|\xi{-}\ol\mfp|^{p'})^{1/p'}
    \dd y &\text{ for }\xi\geq \ol\mfp,\\ 
   \int_0^1 |\xi{-}\mfp(y)| \dd y &\text{ for }\xi\in [\ul\mfp,\ol\mfp],\\ 
\int_0^1( |\xi{-}\mfp(y)|^{p'}{-}|\xi{-}\ul\mfp|^{p'})^{1/p'}
    \dd y &\text{ for }\xi\leq \ul\mfp. \ea\right. 
\]

For $\xi>\ol\mfp$ we have $\xi{-}\mfp(y)> 
\big((\xi{-}\mfp(y))^{p'}{-} a^{p'}\big)^{1/p'}$ for all $ a\in
{]0,\xi{-}\ol\mfp]}$. This yields
\begin{align*}
\xi & = \int_0^1 (\xi{-}\mfp(y))\dd y \gneqq \int_0^1
\big((\xi{-}\mfp(y))^{p'}-(\xi{-}\ol\mfp)^{p'}\big)^{1/p'}\dd y = M_1(\xi).
\end{align*}
Thus, $\xi>\ol\mfp$ cannot occur in the case $\xi\in
\pl_v\calM(0,\xi)$. Similarly, $\xi <\ul\mfp$ is impossible. 
In the remaining case $\xi\in [\ul\mfp,\ol\mfp]$ we have
$\calM(0,\xi)=M_0(\xi)=0$, i.e.\ in \eqref{eq:M.3Cond} with $v=0$
the first condition implies the middle condition. 

Now we start from the third condition $0 \in \pl_\xi
\calM(0,\xi)$. From $\calM(0,\xi)=M_0(\xi)=\min_{\pi\in
  [\ul\mfp,\ol\mfp]} \calR^*(\xi{-}\pi)$ we obtain $\xi \in
[\ul\mfp,\ol\mfp]$ and thus the middle condition
$\calM(0,\xi)=M_0(\xi)=0\xi$ again holds. 

\underline{Step 4: The case $v\neq 0$.} It suffices to consider the
case $v>0$, as $v<0$ can be treated similarly.  For a simpler analysis
we transform this to the variables $(\xi,\HH)$. According to the
formulas in Lemma \ref{le:M.vs.calW} we have to show the equivalence
(where $v=\calV(\xi,\HH)=1/\calW_\HH(\xi,\HH)$)
\begin{equation}
  \label{eq:calW.3Cond}
  \xi=-\calW(\xi,\HH) \ \Longleftrightarrow \   
  \HH-v \calW(\xi,\HH) = \xi  v \ \Longleftrightarrow \ 
  v = -v \calW_\xi(\xi,\HH).  
\end{equation} 
It is obvious that all the relations hold for $\HH=0$ and $\xi \geq
\ol\mfp$. 

Concerning the first relation in \eqref{eq:calW.3Cond}, the strict
monotonicity of $\psi_*$ implies $-\calW(\xi,\HH)>-\calW(\xi,0)=
\int_0^1 |\xi{-}\mfp(z)|\dd x \geq |\xi|$. So there cannot be
solutions with $\HH<0$. The solution set for $\HH=0$ is clearly given by
$\set{\xi}{ \xi\geq \ol\mfp}$.  For $\HH>0$ the function $\calW$ is only
defined for $0<\HH\leq \calR^*(\xi{-}\ol\mfp) $, where we used
$\calW\leq 0$ such that $\xi\geq 0 $ for solution of the left
relation. Again using the strict monotonicity of $\psi_*$ we conclude
$-\calW(\xi,\HH)<\calW(\xi,0) =\xi$ because of $ \xi\geq \ol\mfp$, which
follows from $\HH>0$. Hence, the solution set of the left relations is $
\set{(\xi,\HH)}{ \HH=0,\ \xi\geq \ol\mfp}$. Clearly, on this set the
middle relation holds.
 
We now study the solution set of the right relation in
\eqref{eq:calW.3Cond}, which simplifies to $\calW_\xi(\xi,\HH)=-1$
because of our assumption $v>0$. Obviously, we have
$\calW_\xi(\xi,0)=-1$ for $\xi>\ol\mfp$ 
and $\calW_\xi(\xi,0)=+1$ for $\xi<\ul\mfp$. 
As in the proof of
Theorem \ref{th:M.cvx.xi} we have  
\[
 \calW_{\xi \HH}(\xi,\HH)=-\int_0^1
 \psi''_*\big(\HH{-}\calR^*(\xi{-}\mfp(z))\big)G_\xi(\xi,z)
\dd z, \text{ where }G(\xi,z):= \calR^*(\xi{-}\mfp(z) ) .
\]
For $\xi\not\in [\ul\mfp,\ol\mfp]$ the sign of
$G_\xi(\xi,z)$ equals that of $\xi$, hence we conclude  
$\xi \calW_{\xi \HH}(\xi,\HH) <0$ for $\xi\not\in
[\ul\mfp,\ol\mfp]$. This implies 
\[
\HH>0 \ \Longrightarrow \ \left\{\ba{cl}
\calW_\xi(\xi,\HH)<\calW(\xi,0)=-1& \text{for }\xi>\ol\mfp, \\
\calW_\xi(\xi,\HH)>\calW(\xi,0)=1& \text{for }\xi<\ul\mfp.
\ea\right.
\]
Moreover, for $\HH<0$ and $\xi>\ol\mfp$ we find
$\calW_\xi(\xi,\HH)>\calW(\xi,0)=-1$. Now, restricting to the case of a
power-type dissipation potential $\calR(v)=c|v|^p$ (with $p>1$ as in
Theorem we have \ref{th:M.cvx.xi}) we have $ \calW_{\xi\xi}(\xi,\HH)<0$
for $\xi\in \R$ and $\HH<0$. Thus, for $\xi \leq
\ol\mfp$ and $\HH<0$ we obtain  the estimate 
\[
\calW_\xi(\xi,\HH)\geq \calW_\xi( \ol\mfp{+}1 , \HH)\geq
\calW_\xi( \ol\mfp{+}1 , 0) = -1.
\]
Altogether we conclude that the solution set of the right relation in
\eqref{eq:calW.3Cond} is exactly the same for the left relation. 
\end{proof}

We emphasize that the restriction to the power-law potentials $\calR$
is a sufficient condition for the property that $\calM$ is a
bipotential. However, this is certainly not necessary. We essentially
need the two nontrivial conditions (i) that $\xi\mapsto
\calM (v,\xi)$ is convex for all $v$ and (ii) $\xi \mapsto \calW(\xi,\HH)$ is
concave for all $\HH<0$.  

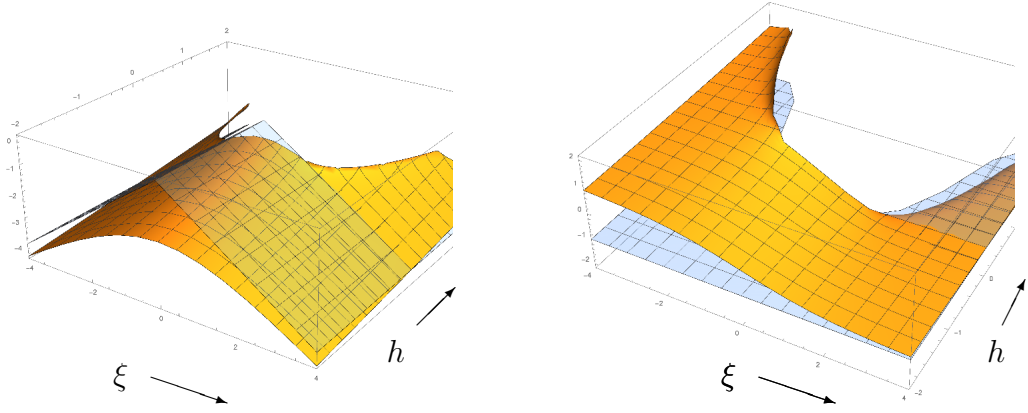
\begin{figure}
\centerline{\includegraphics[width=0.4\textwidth]{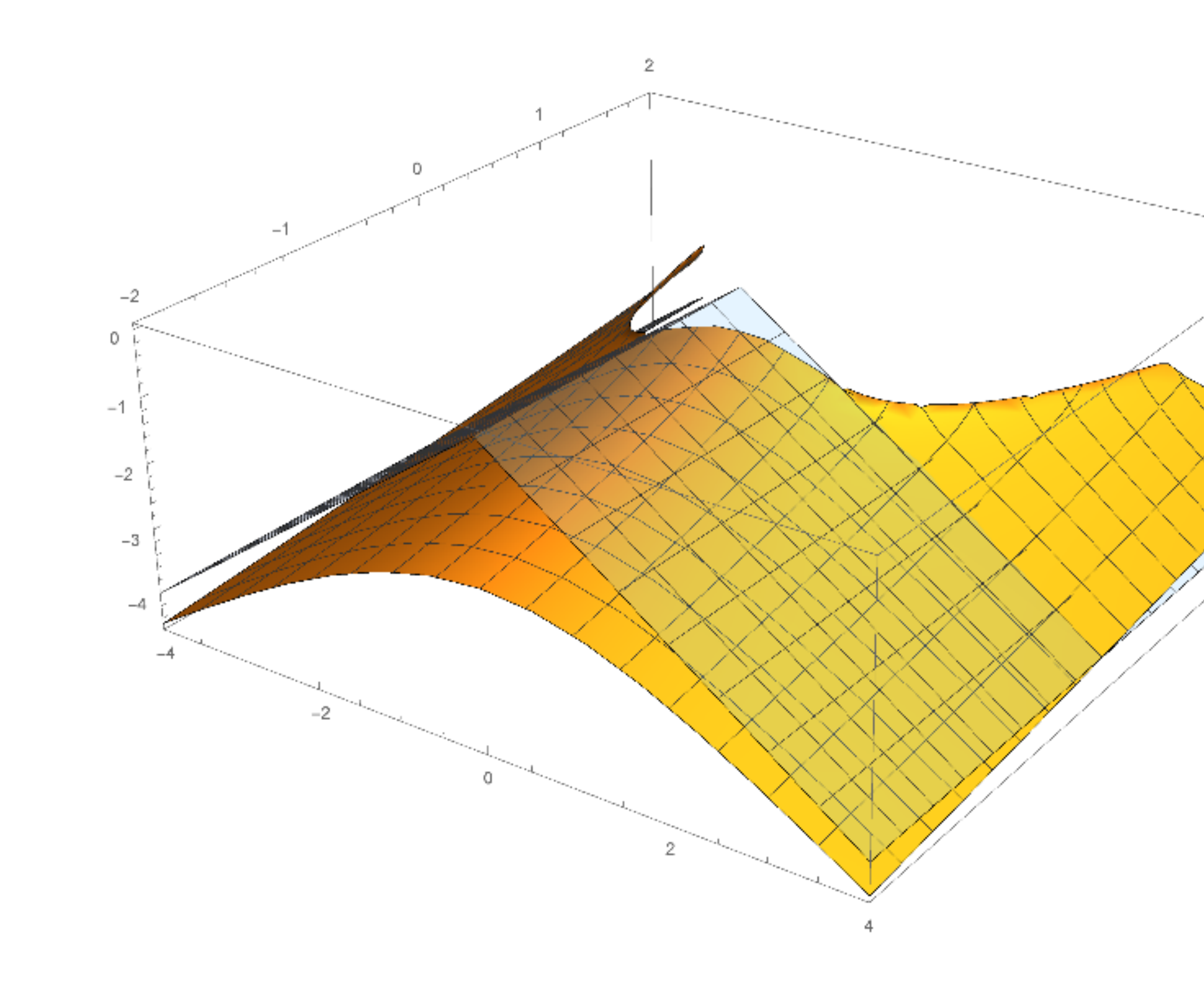}
\qquad 
\includegraphics[width=0.4\textwidth]{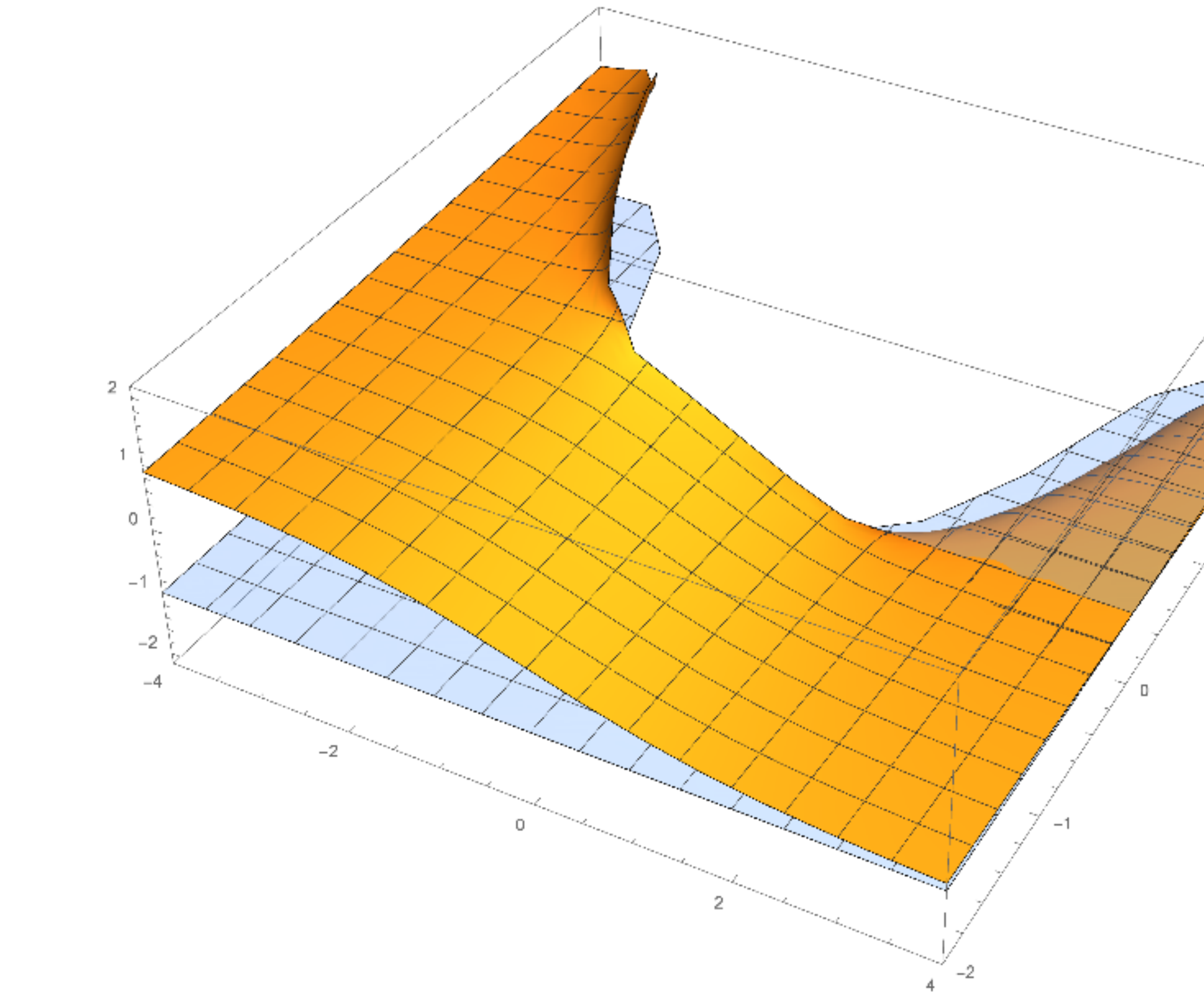}
\unitlength1cm
\begin{picture}(0,0)
\put(-12,0.2){$\xi$} \put(-11.5,0.3){\vector(3,-1){1}}
\put(-4,0.2){$\xi$} \put(-3.5,0.3){\vector(3,-1){1}}
\put(-0.5,0.5){$\HH$} \put(-0.3,1){\vector(1,2){0.3}}
\put(-8.4,0.5){$\HH$} \put(-8,1){\vector(1,1){0.5}}
\put(-4,0.2){$\xi$} \put(-3.5,0.3){\vector(3,-1){1}}
\end{picture}}
\caption{Left: The function $(\xi,\HH)\mapsto \calW(\xi,\HH)$ compared
  with $-|\xi|$. Right: The function $\calW_\xi$ compared with
  $-1$. In both cases the intersection occurs for $\HH=0$ and $\xi\geq
  \ol\mfp=1$. }
\label{fig:Bipot}
\end{figure}

\section{Discussion}
\label{se:Discuss}

Here we provide some discussion points concerning the notions of
evolutionary \textGamma-convergence.  But
first in Section \ref{su:Disc.WW} we highlight that it is important to study the
\textGamma-convergence of $\calJ_\eps$ in the weak$\ti$strong
topology, since using the weak$\ti$weak topology results in a smaller
dissipation function $\calM_\rmw$ that is obviously useless, as it
does not longer satisfies the estimate $\calM_\rmw(u,v,\xi)\geq v\xi$.
In Section \ref{su:EGCWeakT} we recall the notion of
\emph{evolutionary \textGamma-convergence of weak-type} introduced in
\cite{Visi15?SSFE}. Also there, it is strongly highlighted that the
topology for \textGamma-convergence needs to be strong enough to make
the bilinear mapping $(v,\xi) \mapsto \int_0^T \langle
\xi(t),v(t)\rangle \dd t $ continuous. The last subsections 
highlight the
difference between EDP-convergence and relaxed EDP-convergence.

\subsection{\textGamma-limit in weak$\ti$weak topology}
\label{su:Disc.WW}

We now consider $\mfJ_\eps$ on $\rmW^{1,p}(0,T)\ti \rmL^{p'}(0,T)$
equipped with the weak$\ti$weak topology, which is the natural topology for the
family $\mfJ_\eps$ in the sense that is exactly the coarsest topology
in which we have equi-coercivity (i.e.\  $\mfJ_\eps(u_\eps,\xi_\eps) \leq C_1
$ implies $  \|u_\eps\|_{\rmW^{1,p}} + \|\xi_\eps\|_{\rmL^{p'}} \leq
C_2$). Fortunately, in our wiggly-energy model we have a better
convergence for $\xi_\eps$ because of the relation $\xi_\eps
=-\rmD_u\calE_\eps(\cdot,u_\eps) + \Omega_\eps(u_\eps)$, which gave
strong convergence.   

Here we want to highlight that taking the
\textGamma-limit in the weak$\ti$weak topology leads to a functional 
\[
\mfJ_\rmw: (u,\xi)\mapsto \int_0^T \calM_\rmw(u,\dot u,\xi)\dd t
\]
 that is too small. Indeed, using the same techniques as in Section
 \ref{se:Homog}  it can be shown that the
  \textGamma-limit with respect to this weaker topology is given by
 \[
\calM_\rmw(u,v,\xi) = \min_{z \in \rmW^{1,p}_v(0,1)}\bigg\{\int_0^{1}
 \calR\big(u,|v|\dot{z}(s)\big)\dd s
 +\calR^{*}\Big(u, \xi -\int_0^{1} \partial_y\ppp(u,z(s))\dd
 s\Big)\; \bigg\}.
\] 
We clearly obtain $\calM_\rmw\leq \calM $ with $\calM $ from  
\eqref{eq:I.rmM0}. Note that $\calM_\rmw(u,v,\xi)$ is jointly convex
in $(v,\xi)$,  so it must be smaller that $\calM(u,v,\xi)$ in cases where
the latter is not convex in $\xi$.

While convexity may be considered as a nice add-on, the lower bound
$\calM_\rmw(u,v,\xi)\geq v\xi$ is essential for the energy-dissipation
principle to go back from the energy-dissipation estimate to the
subdifferential inclusion. However, $\calM_\rmw$ does no
longer satisfy this important lower bound. To see this, we consider
the example $\calR(u,\dot u)=\frac12\dot 
u^2$ and $\ppp(u,y)=a|y|$ for $|y|\leq \frac12$ and then periodically
extended. Assuming $a,v>0$ and inserting the piecewise interpolant of
the points $z(0)=0$, $z(\frac34)= \frac12$, and $z(1)=1$ into the
minimization problem defining $\calM_\rmw$ a simple calculation yields
the upper bound $\calM_\rmw(v,\xi)\leq \frac23 v^2 + \frac12(\xi-\frac
a2)^2$. Hence, we obtain $\calM_\rmw(\frac a2,\frac a2)\leq \frac16 a^2$
which is strictly smaller than $v\xi= \frac14 a^2$.

\subsection{Evolutionary $\Gamma$-convergence of \emph{weak-type}}
\label{su:EGCWeakT}

The definition of EDP-convergence and in particular that of relaxed
EDP-convergence is relatively close to the notion of
\emph{evolutionary $\Gamma$-convergence of the weak-type} introduced
in \cite{Visi13VFSS,Visi15?SSFE,Visi17?EGCW}. There the class of
monotone flows in the form  
\begin{equation}
  \label{eq:MonoFlow}
  \dot q + \bfA(q) \ni \ell(t)
\end{equation}
are studied, where $\bfA$ is a maximal monotone operator on an
evolution triple $\bfQ\subset \bfH\sim \bfH^* \subset \bfQ^*$. 
The operator $\bfA$
can be represented in the sense of Fitzpatrick by a function
$\calG:\bfQ\ti \bfQ^*\to \R$ as follows:
\begin{align*}
&\calG \text{ is convex and } \calG(q,\xi)\geq \langle \xi,q\rangle
\text{ for all }(q,\xi)\in \bfQ\ti \bfQ^*,\\
& \xi \in \bfA(q) \quad \Longleftrightarrow \quad (q,\xi) \in \mathsf
C_\calG=\bigset{(\eta,v)}{ \calG(v,\eta)=\langle \eta,v\rangle}. 
\end{align*}
The energy-dissipation principle is replaced by an \emph{extended
Brezis-Ekeland-Nayroles principle}, namely
\[
\frac12\|q(T)\|_\bfH^2 + \mathfrak G(q,\ell) = \frac12\| q(0)\|_\bfH^2, \text{
  where } 
\mathfrak G(q,\ell):=\int_0^T \big(\calG(q,\ell{-}\dot q)- \langle
\ell,q\rangle\big) \dd t. 
\]
For families of monotone flows and associated representation functions 
$\calG_\eps$ one can then study ``static \textGamma-convergence'' for
the functionals $\mathfrak G_\eps$. The applicability of this theory
to monotone operators certainly generalizes aspects of our general
EDP-convergence in Section \ref{su:EGC.GS}, however it is also more
restrictive as these monotone flows are only singly nonlinear, which
means for gradient systems $(\bfQ,\calE,\calR)$ that $\calR(q,v)$ cannot
depend on $q\in \bfQ$ and that either $\calE$ or $\calR$ are quadratic. 

More general classes of pseudo-monotone operators are considered with a
further extension of the Brezis--Ekeland--Nayroles principle in 
\cite{Visi17?SCSP}.

\subsection{Mosco convergence implies EDP-convergence}
\label{su:Mosco}

A simple abstract framework for EDP-convergence can be developed in cases
where we have 
\[
\calE_\eps \Gweak \calE_0 \quad \text{and} \quad \calR_\eps \Gweak \calR_0.
\]
However, these two convergences are certainly not sufficient for
EDP-convergence, as they are satisfied in our wiggly-energy model with
$\calR_0=\calR$, but $(\R,\calE_0,\calR)$ is certainly not the correct
limit. 

A general abstract theory was developed in
\cite[Th.\,4.8]{MiRoSa13NADN}, see also
\cite[Sec.\,3.3.2]{Miel16EGCG} for a simplified case and
discussion. It relies on the more restrictive notion of Mosco
convergence $\calF_\eps \Mto \calF_0$ on a Banach space $\bfQ$, which
means $\calF_\eps \Gto \calF_0$ and $\calF_\eps \Gweak \calF_0$. 
 
The setup starts from a reflexive Banach space $\bfQ$ and a densely
and compactly embedded energy space $\bfZ\Subset \bfQ$. The energies
$\calE_\eps: \bfQ\to \R_\infty:=\R\cup\{\infty\}$ are assumed to be
equi-coercive in $\bfZ$ and satisfy $\calE_\eps \Gweak \calE_0$ in
$\bfZ$, which is equivalent to $\calE \Mto \calE_0$ in $\bfQ$. 

The dissipation potentials $\calR_\eps:\bfZ\ti \bfQ\to
[0,\infty]$ satisfy $p$-equicoercivity with $p>1$:
\[
\exists\, c_1,C_2,C_3>0\ \forall\,\eps\in [0,1]\ \forall\, q\in \bfZ\
\forall\, v\in \bfQ:   
\quad c_1\|v\|_\bfQ^p-C_3 \leq \calR_\eps(q,v) \leq C_2\|v\|_\bfQ^p
+C_3.
\]
The convergence of $\calR_\eps$ to $\calR_0$ is the following Mosco
convergence:
\begin{equation}
  \label{eq:calR.Mto}
  \text{If } q_\eps \weak q_0 \text{ in } \bfZ, \ \text{ then }
\calR_\eps(q_\eps, \cdot) \Mto \calR_0(q_0,\cdot) \text{ in } \bfQ. 
\end{equation}
Still these conditions are not enough for EDP-convergence (as they
hold in our wiggly-energy model), so the crucial additional condition
in \cite[Thm.\,4.8]{MiRoSa13NADN} is the closedness of the
subdifferentials of the family $(\calE_\eps)_{\eps\in [0,1]}$, i.e.\ 
\[
 \left\{ \ba{c} 
    q_\eps \to q_0, \  \xi_\eps \weak \xi_* \text{ in }\bfQ^*, 
   \\  \xi_\eps \in \pl\calE_\eps(q_\eps) , \ 
     \calE_\eps(q_\eps) to E_0 
\ea \right\} \ \Longrightarrow \  
 \xi_*\in \pl\calE_0(q_0) \text{ and }E_0=\calE_0(q_0).
\]
This can be achieved if one has equi-$\lambda$-convexity, i.e.\ there
exists $\lambda_*\in \R$, such that all functions $q \mapsto
\calE_\eps(q)+\lambda_*\| q\|_\bfQ^2$ are convex.  

If all these conditions (together with some other standard conditions)
hold, then one obtains $(\bfQ,\calE_\eps,\calR_\eps) \EDPto
(\bfQ,\calE_0,\calR_0)$. Indeed, in  \cite[Thm.\,4.8]{MiRoSa13NADN}
EDP-convergence is not mentioned, however, the proof of evolutionary 
\textGamma-convergence is done in a way which exactly shows all
ingredients of EDP-convergence. 

This is in contrast to the typical Sandier-Serfaty approach
\cite{SanSer04GCGF,Serf11GCGF}, where only estimates along the precise
solutions of the gradient flows are needed.

\subsection{EDP-convergence versus relaxed EDP-convergence}
\label{su:EDPrelaxEDP}
 
More advanced cases of EDP-convergence are discussed in
\cite{LMPR17MOGG}. We recall that EDP-convergence distinguishes from
relaxed EDP-convergence that the limiting dissipation functional
$\mfD_0$ is given in terms of $\calM$ having the form 
\[
\calM(q,v,\xi)= \calM_\eff(q,v,\xi):=\calR_\eff(q,v)+ \calR^*_\eff(q,\xi).
\]
In the general case this identity is not true, and it is interesting
to ask whether we have an estimate of the form $\calM\geq \calM_\eff$,
since only this estimate is needed to show evolutionary
\textGamma-convergence. Yet, for our wiggly-energy model 
Proposition 
\ref{pr:Bound.Reff} yields the opposite estimate, namely
\[
\calM_\eff(0,\xi)=\calR_\eff^*(\xi) > \calR^*(\xi{-}\ol\mfp) =
M_0(\xi) = \calM(0,\xi) 
\text{ for all }\xi >\ol\mfp.
\]
Moreover, for $0<v\ll 1$ and $\xi\in [\ul\mfp,\ol\mfp]$ we have
$\calM(v,\xi)= v M_1(\xi)+ o(v)_{v\to 0}$ with 
$M_1(\xi)= \int_0^1|\xi{-}\mfp(y)|\dd y $, see Lemma
\ref{le:Exp.M0}. For $\xi\in {[0,\ol\mfp[}$ we have $M_1(\xi)<
M_1(\ol\mfp) =\ol\mfp$, so we again have $\calM_\eff(v,\xi)= v \ol\mfp
+o(v)> \calM(v,\xi)=v M_1(\xi)+o(v)$. 

We feel that this is the typical feature of relaxed EDP-convergence,
and conjecture that $\calM(v,\xi)\leq \calM_\eff $ and that equality
holds only in the case of true EDP-convergence. Of course, the
difference of $\calM_\eff - \calM$ always vanishes on the contact set
$\mathsf C_\calM$, which highlights that the representation of the
operator $v \mapsto \pl\calR_\eff(v)$ can well be given in terms of a
function $\calM$ that is smaller that $\calM_\eff$.  We illustrate this
by looking at the following special case.

\begin{example}\label{ex:Meff-M}\slshape
We want to justify our conjecture by an example calculation for the
case $\calR(v)=\frac12v^2$ and $\mfp$ taking the two values $\pm 1$
with weight $1/2$. 

For $\xi>1$ we find $\rmD\calR^*_\eff(\xi)=\xi-1/\xi$ and hence
$\calR_\eff$ and $\calR^*_\eff$ have the form 
\[
\calR_\eff(v)= \frac14 \big(v^2 + |v|\sqrt{v^2{+}4}\big)+\mafo{Asinh}(|v|/2) \
\text{ and } \ 
\calR_\eff:\xi \mapsto \frac12 \max\{0,\xi^2 {-}1\} -
\max\{0,\log \xi\}.
\] 
We can also evaluate $\calW$ explicitly and find
$\calW(\xi,h)=-\frac12\big(\sqrt{(\xi{-}1)^2-2h}+\sqrt{(\xi{-}1)^2-2h}
\big) $, and Lemma \ref{le:M.vs.calW} gives $\ol\calM(\xi,h)=h-V(\xi,h)\calW(\xi,h)=
h+\sqrt{(\xi{-}1)^2-2h} \sqrt{(\xi{-}1)^2-2h}$, where
$v=V(\xi,h):=1/\calW_h(\xi,h)>0$. Thus, we may check our conjecture
for positive $v$, because $v=V(\xi,h)$ tends to $0$ for $h\to -\infty$
and to $\infty$ for $h\nearrow \frac12 \min\{ (\xi{+}1)^2, (\xi{-}1)^2
\}$. We can now compare $\calM$ and $\calM_\eff$ by plotting them over
the $(\xi,h)$-plane, and indeed Figure \ref{fig:Meff-M} shows that the
conjecture holds for this simple case.
\begin{figure}
\includegraphics[width=0.45\textwidth]{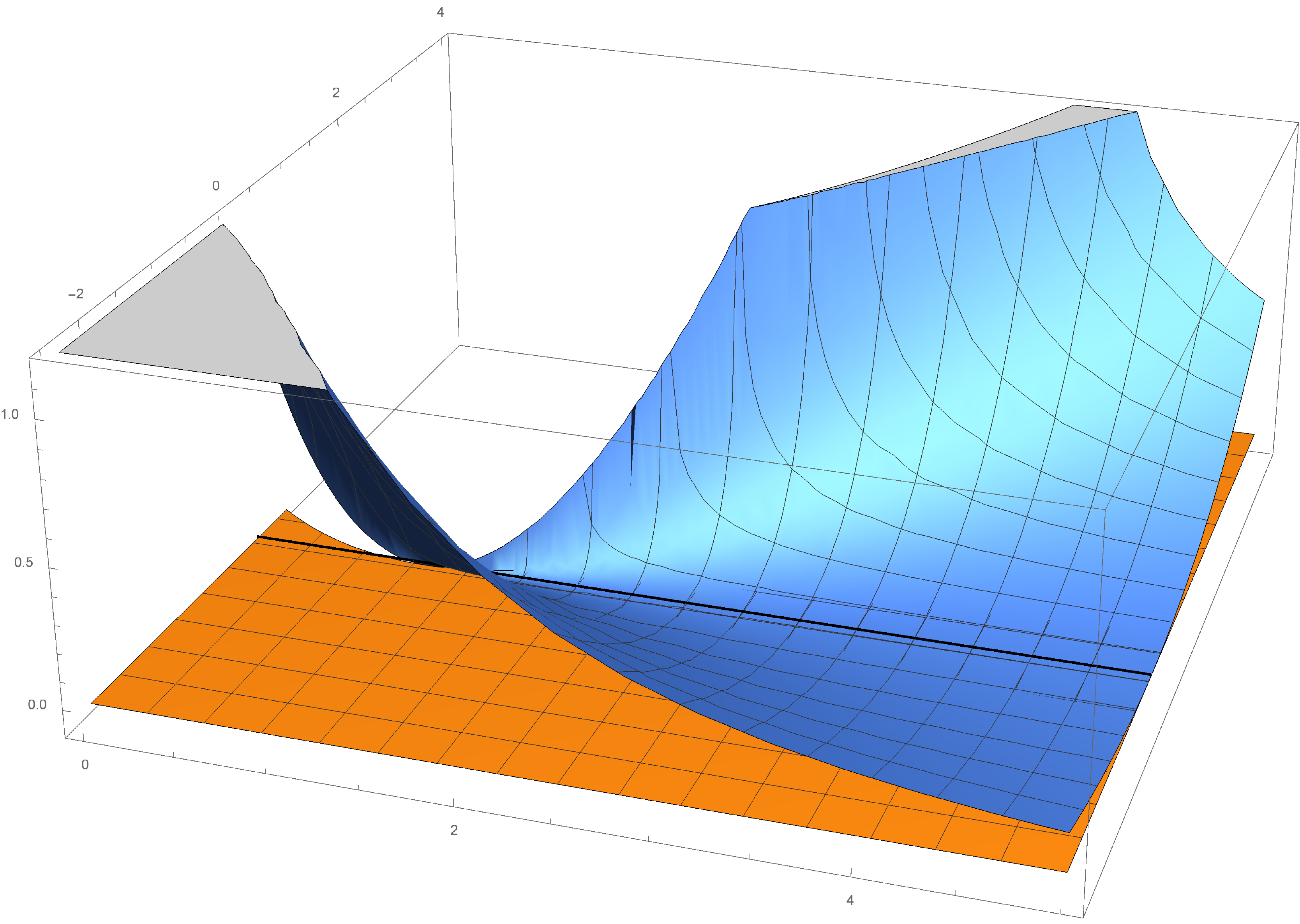}\qquad
\includegraphics[width=0.5\textwidth]{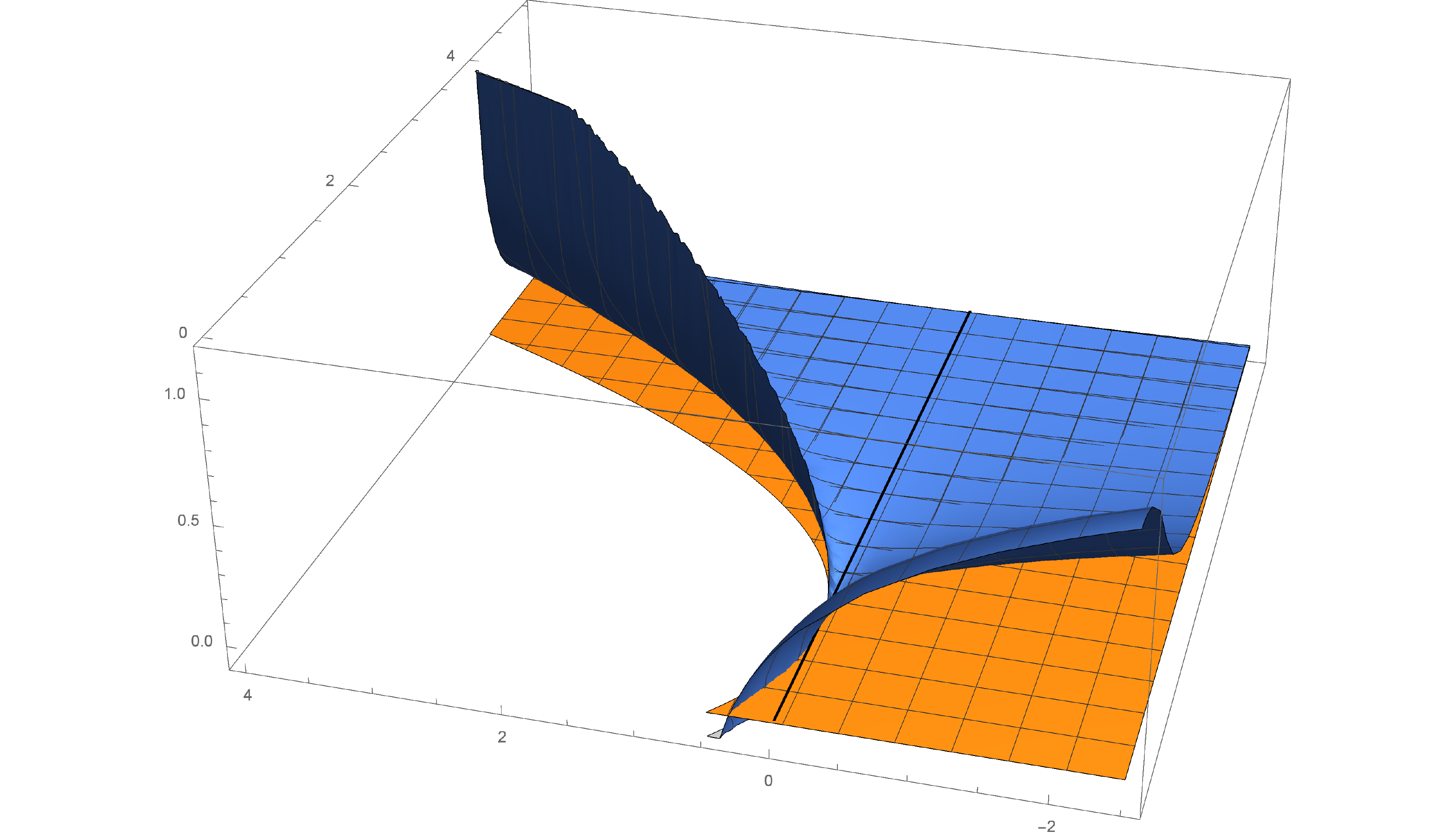}
\caption{Confirmation of the conjecture $\calM\leq \calM_\eff$ for a
  special case: The left figure shows that the graph of $(\xi,h)\mapsto
  \calM(V(\xi,h),h)- \xi V(\xi,h)$ lies above $0$. The left figure
  shows that $(\xi,h)\mapsto
  \calM_\eff(V(\xi,h),h)- \calM(V(\xi,h)$  lies above $0$. In both
  cases the functions equal $0$ the bold line
  $\set{(\xi,h)}{h=0,\;\xi\geq 1}$. The values for $|\xi|<1$ and $h>0$
are not relevant, because they do not correspond to any $v$.}
\label{fig:Meff-M}
\end{figure} 
\end{example}

\subsection{Non-convergence of primal and dual dissipation parts}
\label{su:PartitionD}

The main observation is that EDP-convergence, and
even more relaxed EDP-convergence, are able to work in cases where the
nature of the dissipation potential can change its structure. 
In our wiggly-energy model we found that even though $\calR_\eps =
\calR$ we have $\calR_\eff \neq \calR$. Moreover, for quadratic
$\calR$ we obtain an $\calR_\eff$ that behaves like $v \mapsto \ol\mfp
|v|$ for small $v$. 

Such nontrivial changes in the dissipation structure were already
observed in \cite{LMPR17MOGG}. For instance it is shown that 
the diffusion through a 
layer of thickness $\eps$ with a mobility $a\eps$ has a EDP-limit that
describes the jump conditions at a membrane with transmission
coefficient $a>0$. The natural gradient structure for diffusion is
$(\rmL^1(\Omega),\calE,\calR_\eps)$ with the relative entropy
$\calE(u)=\int_\Omega ( u\log u -u +1)\dd x$ and the quadratic
dissipation potentials of Wasserstein-Kantorovich type, namely 
\[
\calR_\eps^*(u,\xi)= \int_\Omega \tfrac{A_\eps(x)}{2} 
 \,|\nabla \xi(x)|^2 \, u(x) \dd x 
\]
The mobility $A_\eps$ equals $1$ except for the small layer. 
It is shown in \cite{Lier12VME,LMPR17MOGG} that we have
EDP-convergence to $(\rmL^1(\Omega),\calE,\calR_\eff)$, and the
surprising fact is that $\calR_\eff$ is non-quadratic in $\xi$,
because it involves exponential function of the jump of $\xi$ over the
limiting membrane. 

This change in the structure of the dissipation potentials highlights
a general point in EDP-convergence, even when we restrict to exact
solutions $q_\eps$ of the gradient systems
$(\bfQ,\calE_\eps,\calR_\eps)$. Clearly, we have 
\[
\calE_\eps(q_\eps(T)) + \mfD_\eps(q_\eps) = \calE_\eps(q_\eps(0)).
\]
Assume $q_\eps(0)\weak q_0(0)$ and $\calE_\eps(q_\eps(0)) \to
\calE_0(q_0(0))$ (i.e.\ well-prepared initial conditions), the
convergence $q_\eps \weak q_0$ in $\rmW^{1,p}(0,T;\bfQ)$ implies 
\[
\calE_\eps(q_\eps(t)) \to \calE_0(q_0(t)) \text{ for all }t\in [0,T]
\quad \text{and }\mfD_\eps(q_\eps) \to \mfD_0(q_0).
\]
This means that $q_\eps(t)$ is a recovery sequence for the energies
$\calE_\eps$ and $q_\eps(\cdot)$ is a recovery sequence for the
dissipation functionals.

However, the dissipation potential $\mfD_\eps$ can be understood as the sum
of a primal part $\mfD_\eps^\prim$ given via $\calR_\eps$ and of a dual part
$\mfD_\eps^\dual$ given via $\calR^*_\eps$:
\[
\mfD_\eps^\prim (q)=\int_0^T\calR_\eps(q(t),\dot q(t))\dd t \quad
\text{and} \quad \mfD_\eps^\dual (q)=\int_0^T\calR^*_\eps
 \big(q(t),{-}\rmD\calE_\eps(q(t))\big)\dd t .\quad
\]

To understand how the effective dissipation potential $\calR_\eff$
differs from the limits of $\calR_\eps$ we may consider the 
separate limits 
\[
\mathsf D^\prim(q_0) := \lim_{\eps \to 0} \mfD^\prim_\eps(q_\eps) \quad 
\text{and} \quad 
\mathsf D^\dual(q_0) := \lim_{\eps \to 0} \mfD^\dual_\eps(q_\eps) 
\]
along solutions $q_\eps$ of $(\bfQ,\calE_\eps,\calR_\eps)$ converging
to a solution $q_0$ of $(\bfQ,\calE_0,\calR_\eff)$. Setting
\[
 \mfD_\eff^\prim(u_0):=\int_0^T
  \calR_\eff(q_0,\dot q_0) \dd t \quad 
\text{and} \quad   \mfD_\eff^\dual(u_0):=
\int_0^T \calR^*_\eff(q_0,{-}\rmD\calE_0(q_0)) \dd t.
\]
we emphasize that, in general,
(relaxed) EDP-convergence does not imply the identities 
\begin{equation}
  \label{eq:sfD=Deff}
  \mathsf D^\prim(q_0)=\mfD_\eff^\prim(u_0) \quad 
\text{and} \quad  \mathsf D^\dual(q_0) = \mfD_\eff^\dual(u_0).  
\end{equation}
However, for the case considered in Section \ref{su:Mosco} these
identities are established in \cite[Eqn.\,(4.29c)]{MiRoSa13NADN}
based on the Mosco convergences $\calR_\eps\Mto
\calR_0=\calR_\eff$, cf.\ \eqref{eq:calR.Mto}. 

The problems in more general cases are most easily understood when
considering $p$-homogenous dissipation potentials $\calR$ with $p>1$. Then, 
Euler's formula gives $\langle \pl\calR(v),v\rangle= p
\calR(v)$ and $\langle \xi, \pl\calR^*(\xi)\rangle= p'
\calR^*(\xi)$. Moreover, we have 
\[
\xi\in \partial\calR(v) \ \Longrightarrow \  p\calR(v)=
\calR(v){+} \calR^*(\xi)= \langle \xi,v\rangle= p'\calR^*(\xi).
\] 
Thus, if all dissipation potentials $\calR_\eps$ are $p$-homogeneous,
we have $\mfD^\prim_\eps(q_\eps)=\frac1p \mfD_\eps(q_\eps)$  and 
$\mfD^\dual_\eps(q_\eps)=\frac1{p'} \mfD_\eps(q_\eps)$, and 
the convergence of $\mfD_\eps(q_\eps)\to \mfD_0(q_0)$ yields 
\[
\mathsf D^\prim(q_0)=\frac1p \mfD_0(q_0)\text{ and } 
\mathsf D^\dual(q_0)=\frac1{p'} \mfD_0(q_0). 
\]
Of course, by (relaxed) EDP-convergence we have the representation 
\[
\mfD_0(q_0) = \int_0^T \calM(q_0,\dot q_0,{-}\rmD\calE(q_0)) \dd t = \mfD_\eff^\prim(q_0)+ \mfD_\eff^\dual(q_0).
\]
Here the second identity follows since $q_0$ is a solutions such that
$(\dot q_0,{-}\rmD\calE_0(q_0))$ lies in $\mathsf C_\calM$, where
$\calM$ equals $\calM_\eff$,  as both functional equal $\langle \xi,v\rangle$
on $\mathsf C_\calM$.

The question as to whether the two identities in \eqref{eq:sfD=Deff}
hold is now reduced to the question whether $\calR_\eff(q,\cdot)$ is
still $p$-homogeneous. Thus, in the Sandier-Serfaty approach, where
$p=2$ for $\eps>0$ as well as for $\eps=0$, we have the desired identity. 

However, in our wiggly-energy model we can start with arbitrary
$p>1$ for $\eps>0$ but end up with $\calR_\eff$ satisfying $\langle
\pl\calR_\eff (u,v),v\rangle = \alpha(u,v)\calR_\eff(v)$ with
$\alpha(u,v)\in {[1,p[}$, see Proposition
\ref{pr:pHom}. Hence, we obtain a strict inequality, namely 
\begin{align*}
&\mathsf D^\prim(u_0)\   =\ \frac1p\mfD_0(u_0)\ =\ \frac1p 
 \int_0^T\big( \calR_\eff(u_0,\dot u_0) + 
     \calR^*_\eff(u_0,{-}\rmD\calE_0(u_0))\big) \dd t\\
 & = \ \frac1p\int_0^T \!\!
\pl_v\calR_\eff(u_0,\dot u_0)\dot u_0 \dd t\ = \ 
\int_0^T \!\!\frac{\alpha(u_0,\dot u_0)}p \calR_\eff(u_0,\dot u_0)\dd t 
 \ \lneqq \ \int_0^T\!  \calR_\eff(u_0,\dot u_0)\dd t.
\end{align*}
Because $\alpha(u,0)=1$ the effect is stronger if $\dot u_0$ is small,
i.e.\ when we are close to the rate-independent case. 

In the membrane limit of thin layers discussed in
\cite{Lier12VME,LMPR17MOGG} we have quadratic dissipation potentials
for $\eps>0$, i.e.\ $p=2$. However, for $\eps=0$ one obtains
$\calR_\eff$ with a growth like $|v| \log |v|$ for $|v|\gg 1$. Again
we have $\langle \pl\calR_\eff (\dot q),\dot q\rangle = b(\dot
q)\calR_\eff(\dot q)$, where $b(\dot q)\leq 2$ and $b(\dot q)<2$ for
certain $\dot q$.  However, there the effect is stronger for large
$\dot q$ and disappears for $\dot q\to 0$.

For both cases we see that in the limiting primal part
of the dissipation
functional $\int_0^T\calR_\eff(q_0,\dot q_0) \dd t  $  is larger than
the limit $\mathsf D^\prim(q_0) = \lim_{\eps \to 0}
\mfD^\prim_\eps(q_\eps)$. This is also seen in the inequality
$\calR_\eps \Gweak \calR_0 \leq \calR_\eff$. 
We interpret this as the effect of
microscopic dissipative processes that need to be modeled on the
macroscale for the limit system $(\bfQ,\calE_0,\calR_\eff)$. 

It is an interesting question to understand whether relaxed
EDP-convergence always leads to an increase for the primal part of the
dissipation functional; more precisely, do we always have 
$\mathsf D^\prim(q_0) \leq  \int_0^T\calR_\eff(q_0,\dot q_0) \dd t$?

{\small
\bibliographystyle{my_alpha}
\bibliography{alex_pub,bib_alex,bib_pwd}}

\newcommand{\etalchar}[1]{$^{#1}$}
\def\cprime{$'$}
\begin{thebibliography}{11}\itemsep0.1em

\bibitem[AbK88]{AbKn_88}
{\scshape R.~Abeyaratne {\upshape and} J.~K.~Knowles}.
\newblock On the dissipative response due to discontinuous strains in bars of
  unstable elastic material.
\newblock {\em Internat. J. Solids Structures}, 24(10), 1021--1044, 1988.

\bibitem[AbK97]{AbKn_97}
{\scshape R.~Abeyaratne {\upshape and} J.~K.~Knowles}.
\newblock On the kinetics of an austenite $\rightarrow$ martensite phase
  transformation induced by impact in a {C}u{A}l{N}i shape-memory alloy.
\newblock {\em Acta Materialia}, 45(4), 1671 -- 1683, 1997.

\bibitem[ACJ96]{AbChJa96KMWE}
{\scshape R.~Abeyaratne, C.-H.~Chu, {\upshape and} R.~James}.
\newblock Kinetics of materials with wiggly energies: theory and application to
  the evolution of twinning microstructures in a {C}u-{A}l-{N}i shape memory
  alloy.
\newblock {\em Phil. Mag. A}, 73, 457--497, 1996.

\bibitem[AGS05]{AmGiSa05GFMS}
{\scshape L.~Ambrosio, N.~Gigli, {\upshape and} G.~Savar{\'e}}.
\newblock {\em Gradient flows in metric spaces and in the space of probability
  measures}.
\newblock Lectures in Mathematics ETH Z\"urich. Birkh\"auser Verlag, Basel,
  2005.

\bibitem[Att84]{Atto84VCFO}
{\scshape H.~Attouch}.
\newblock {\em Variational Convergence of Functions and Operators}.
\newblock Pitman Advanced Publishing Program. Pitman, 1984.

\bibitem[BdV08a]{BuDeVa08ECBG}
{\scshape M.~Buliga, G.~{de Saxc\'e}, {\upshape and} C.~Valle\'e}.
\newblock Existence and construction of bipotentials for graphs of multivalued
  laws.
\newblock {\em J. Convex Anal.}, 15(1), 87--104, 2008.

\bibitem[BdV08b]{BuDeVa08?NMCM}
{\scshape M.~Buliga, G.~{de Saxc\'e}, {\upshape and} C.~Valle\'e}.
\newblock Non maximal cyclically monotone graphs and construction of a
  bipotential for the {C}oulomb's dry friction law.
\newblock {\em arXiv:0802.1140v1}, 2008.

\bibitem[Bha99]{Bhattacharya_99}
{\scshape K.~Bhattacharya}.
\newblock Phase boundary propagation in a heterogeneous body.
\newblock {\em R. Soc. Lond. Proc. Ser. A Math. Phys. Eng. Sci.}, 455(1982),
  757--766, 1999.

\bibitem[BoP16]{BonPel16QRIL}
{\scshape G.~A.~Bonaschi {\upshape and} M.~A.~Peletier}.
\newblock Quadratic and rate-independent limits for a large-deviations
  functional.
\newblock {\em Contin. Mech. Thermodyn.}, 28(4), 1191--1219, 2016.

\bibitem[Bra02]{Brai02GCB}
{\scshape A.~Braides}.
\newblock {\em $\Gamma$-Convergence for Beginners}.
\newblock Oxford University Press, 2002.

\bibitem[Bra13]{Brai13LMVE}
{\scshape A.~Braides}.
\newblock {\em Local minimization, Variational Evolution and
  Gamma-convergence}.
\newblock Lect.\ Notes Math.\ Vol.\ 2094. Springer, 2013.

\bibitem[Bud17]{BulDes17SBEN}
{\scshape M.~Buliga {\upshape and} G.~de~Saxc\'e}.
\newblock A symplectic {B}rezis-{E}keland-{N}ayroles principle.
\newblock {\em Math. Mech. Solids}, 22(6), 1288--1302, 2017.

\bibitem[{Dal}93]{Dalm93IGC}
{\scshape G.~{Dal Maso}}.
\newblock {\em An Introduction to {$\Gamma$}-Convergence}.
\newblock Birkh\"auser Boston Inc., Boston, MA, 1993.

\bibitem[DKW17]{Dondl:2017ut}
{\scshape P.~W.~Dondl, M.~W.~Kurzke, {\upshape and} S.~Wojtowytsch}.
\newblock {The Effect of Forest Dislocations on the Evolution of a Phase-Field
  Model for Plastic Slip}.
\newblock {\em arXiv.org}, July 2017.

\bibitem[DMT80]{DeMaTo80PEMS}
{\scshape E.~{De Giorgi}, A.~Marino, {\upshape and} M.~Tosques}.
\newblock Problems of evolution in metric spaces and maximal decreasing curve.
\newblock {\em Atti Accad. Naz. Lincei Rend. Cl. Sci. Fis. Mat. Natur. (8)},
  68(3), 180--187, 1980.

\bibitem[EfM06]{EfeMie06RILS}
{\scshape M.~Efendiev {\upshape and} A.~Mielke}.
\newblock On the rate--independent limit of systems with dry friction and small
  viscosity.
\newblock {\em J. Convex Anal.}, 13(1), 151--167, 2006.

\bibitem[EkT76]{EkeTem76CAVP}
{\scshape I.~Ekeland {\upshape and} R.~Temam}.
\newblock {\em Convex Analysis and Variational Problems}.
\newblock North Holland, 1976.

\bibitem[EsC93]{EsCl_93}
{\scshape J.~Escobar {\upshape and} R.~Clifton}.
\newblock On pressure-shear plate impact for studying the kinetics of
  stress-induced phase transformations.
\newblock {\em Materials Science and Engineering: A}, 170(1), 125 -- 142, 1993.

\bibitem[Fen49]{Fenc49CCF}
{\scshape W.~Fenchel}.
\newblock On conjugate convex functions.
\newblock {\em Canadian J. Math.}, 1, 73--77, 1949.

\bibitem[GaM05]{Garroni:2005ve}
{\scshape A.~Garroni {\upshape and} S.~M{\"u}ller}.
\newblock {$Gamma$-limit of a phase-field model of dislocations}.
\newblock {\em SIAM Journal on Mathematical Analysis}, 36(6), 1943--1964
  (electronic), 2005.

\bibitem[GaM06]{Garroni:2006tn}
{\scshape A.~Garroni {\upshape and} S.~M{\"u}ller}.
\newblock {A variational model for dislocations in the line tension limit}.
\newblock {\em Archive For Rational Mechanics And Analysis}, 181(3), 535--578,
  2006.

\bibitem[GiD17]{GidDes17GDFB}
{\scshape P.~Gidoni {\upshape and} A.~DeSimone}.
\newblock On the genesis of directional friction through bristle-like mediating
  elements crawler.
\newblock {\em ESAIM Control Optim. Calc. Var.}, 23(3), 1023--1046, 2017.

\bibitem[Jam96]{Jame96HPT}
{\scshape R.~D.~James}.
\newblock Hysteresis in phase transformations.
\newblock In {\em ICIAM 95 (Hamburg, 1995)}, volume~87 of {\em Math. Res.},
  pages 135--154. Akademie Verlag, Berlin, 1996.

\bibitem[Lie12]{Lier12VME}
{\scshape M.~Liero}.
\newblock {\em Variational Methods for Evolution}.
\newblock PhD thesis, Institut f\"ur Mathematik, Humboldt-Universit\"at zu
  Berlin, 2012.

\bibitem[LM{\etalchar{*}}17]{LMPR17MOGG}
{\scshape M.~Liero, A.~Mielke, M.~A.~Peletier, {\upshape and} D.~R.~M.~Renger}.
\newblock On microscopic origins of generalized gradient structures.
\newblock {\em Discr. Cont. Dynam. Systems Ser.~S}, 10(1), 1--35, 2017.

\bibitem[LMS17]{LiMiSa14?OETP}
{\scshape M.~Liero, A.~Mielke, {\upshape and} G.~Savar\'e}.
\newblock Optimal entropy-transport problems and the {H}ellinger--{K}antorovich
  distance.
\newblock {\em Invent. math.}, 2017.
\newblock Accepted. WIAS preprint 2207, arXiv:1508.07941v2.

\bibitem[Men02]{Meno02GSWE}
{\scshape G.~Menon}.
\newblock Gradient systems with wiggly energies and related averaging problems.
\newblock {\em Arch. Rat. Mech. Analysis}, 162, 193--246, 2002.

\bibitem[Mie11]{Miel11CDEB}
{\scshape A.~Mielke}.
\newblock Complete-damage evolution based on energies and stresses.
\newblock {\em Discr. Cont. Dynam. Systems Ser.~S}, 4(2), 423--439, 2011.

\bibitem[Mie12]{Miel12ERID}
{\scshape A.~Mielke}.
\newblock Emergence of rate-independent dissipation from viscous systems with
  wiggly energies.
\newblock {\em Contin. Mech. Thermodyn.}, 24(4), 591--606, 2012.

\bibitem[Mie16]{Miel16EGCG}
{\scshape A.~Mielke}.
\newblock On evolutionary {$\Gamma$}-convergence for gradient systems
  {(Ch.\,3)}.
\newblock In A.~Muntean, J.~Rademacher, {\upshape and} A.~Zagaris, editors,
  {\em Macroscopic and Large Scale Phenomena: Coarse Graining, Mean Field
  Limits and Ergodicity}, Lecture Notes in Applied Math. Mechanics Vol.\,3,
  pages 187--249. Springer, 2016.
\newblock Proc. of Summer School in Twente University, June 2012.

\bibitem[MiT12]{MieTru12DVEC}
{\scshape A.~Mielke {\upshape and} L.~Truskinovsky}.
\newblock From discrete visco-elasticity to continuum rate-independent
  plasticity: rigorous results.
\newblock {\em Arch. Rational Mech. Anal.}, 203(2), 577--619, 2012.

\bibitem[MoM77]{ModMor77EGC}
{\scshape L.~Modica {\upshape and} S.~Mortola}.
\newblock Un esempio di {$\Gamma$}-convergenza.
\newblock {\em Boll. Un. Mat. Ital. B}, 14, 285--299, 1977.

\bibitem[MoP12]{Monneau:2010te}
{\scshape R.~e.~Monneau {\upshape and} S.~Patrizi}.
\newblock {Homogenization of the Peierls-Nabarro model for dislocation
  dynamics}.
\newblock {\em Journal of Differential Equations}, 253(7), 2064--2105, 2012.

\bibitem[MPR14]{MiPeRe14RGFL}
{\scshape A.~Mielke, M.~A.~Peletier, {\upshape and} D.~R.~M.~Renger}.
\newblock On the relation between gradient flows and the large-deviation
  principle, with applications to {M}arkov chains and diffusion.
\newblock {\em Potential Analysis}, 41(4), 1293--1327, 2014.

\bibitem[MRS09]{MiRoSa09MSJR}
{\scshape A.~Mielke, R.~Rossi, {\upshape and} G.~Savar\'{e}}.
\newblock Modeling solutions with jumps for rate-independent systems on metric
  spaces.
\newblock {\em Discr. Cont. Dynam. Systems Ser.~A}, 25(2), 585--615, 2009.

\bibitem[MRS12]{MiRoSa12BVSV}
{\scshape A.~Mielke, R.~Rossi, {\upshape and} G.~Savar\'{e}}.
\newblock {B}{V} solutions and viscosity approximations of rate-independent
  systems.
\newblock {\em ESAIM Control Optim. Calc. Var.}, 18(1), 36--80, 2012.

\bibitem[MRS13]{MiRoSa13NADN}
{\scshape A.~Mielke, R.~Rossi, {\upshape and} G.~Savar\'{e}}.
\newblock Nonsmooth analysis of doubly nonlinear evolution equations.
\newblock {\em Calc. Var. Part. Diff. Eqns.}, 46(1-2), 253--310, 2013.

\bibitem[MRS16]{MiRoSa16BVSI}
{\scshape A.~Mielke, R.~Rossi, {\upshape and} G.~Savar\'{e}}.
\newblock Balanced-viscosity ({B}{V}) solutions to infinite-dimensional
  rate-independent systems.
\newblock {\em J. Europ. Math. Soc.}, 18, 2107--2165, 2016.

\bibitem[PoG12]{PopGra12PTMH}
{\scshape V.~L.~Popov {\upshape and} J.~A.~T.~Gray}.
\newblock Prandtl-{T}omlinson model: {H}istory and applications in friction,
  plasticity, and nanotechnologies.
\newblock {\em Z.\ angew.\ Math.\ Mech. (ZAMM)}, 92(9), 692--708, 2012.

\bibitem[Pra28]{Pran28GKTF}
{\scshape L.~Prandtl}.
\newblock Gedankenmodel zur kinetischen {T}heorie der festen {K}\"orper.
\newblock {\em Z.\ angew.\ Math.\ Mech. (ZAMM)}, 8, 85--106, 1928.

\bibitem[PRV14]{PeReVa14LDSH}
{\scshape M.~A.~Peletier, F.~Redig, {\upshape and} K.~Vafayi}.
\newblock Large deviations in stochastic heat-conduction processes provide a
  gradient-flow structure for heat conduction.
\newblock {\em J. Math. Physics}, 55, 093301/19, 2014.

\bibitem[PuT02a]{PugTru02MTP}
{\scshape G.~Puglisi {\upshape and} L.~Truskinovsky}.
\newblock A mechanism of transformational plasticity.
\newblock {\em Contin. Mech. Thermodyn.}, 14, 437--457, 2002.

\bibitem[PuT02b]{PugTru02RIHB}
{\scshape G.~Puglisi {\upshape and} L.~Truskinovsky}.
\newblock Rate independent hysteresis in a bi-stable chain.
\newblock {\em J. Mech. Phys. Solids}, 50(2), 165--187, 2002.

\bibitem[PuT05]{PugTru05TRIP}
{\scshape G.~Puglisi {\upshape and} L.~Truskinovsky}.
\newblock Thermodynamics of rate-independent plasticity.
\newblock {\em J. Mech. Phys. Solids}, 53, 655--679, 2005.

\bibitem[Roc70]{Rock70CA}
{\scshape R.~T.~Rockafellar}.
\newblock {\em Convex Analysis}.
\newblock Princeton University Press, 1970.

\bibitem[SaS04]{SanSer04GCGF}
{\scshape E.~Sandier {\upshape and} S.~Serfaty}.
\newblock Gamma-convergence of gradient flows with applications to
  {G}inzburg-{L}andau.
\newblock {\em Comm. Pure Appl. Math.}, LVII, 1627--1672, 2004.

\bibitem[Ser11]{Serf11GCGF}
{\scshape S.~Serfaty}.
\newblock Gamma-convergence of gradient flows on {H}ilbert spaces and metric
  spaces and applications.
\newblock {\em Discr. Cont. Dynam. Systems Ser.~A}, 31(4), 1427--1451, 2011.

\bibitem[SK{\etalchar{*}}09]{SKTO09BDSC}
{\scshape T.~J.~Sullivan, M.~Koslowski, F.~Theil, {\upshape and} M.~Ortiz}.
\newblock On the behaviour of dissipative systems in contact with a heat bath:
  Application to andrade creep.
\newblock {\em J. Mech. Phys. Solids}, 57(7), 1058--1077, 2009.

\bibitem[Ste08]{Stef08BEPD}
{\scshape U.~Stefanelli}.
\newblock The {B}rezis-{E}keland principle for doubly nonlinear equations.
\newblock {\em SIAM J. Control Optim.}, 47(3), 1615--1642, 2008.

\bibitem[Sul09]{Sull09AGDR}
{\scshape T.~J.~Sullivan}.
\newblock {\em Analysis of Gradient Descents in Random Energies and Heat
  Baths}.
\newblock PhD thesis, Dept.\ of Mathematics, University of Warwick, 2009.

\bibitem[Tom29]{Toml29MTF}
{\scshape G.~A.~Tomlinson}.
\newblock A molecular theory of friction.
\newblock {\em Philos. Mag.}, 7, 905--939, 1929.

\bibitem[Vis13]{Visi13VFSS}
{\scshape A.~Visintin}.
\newblock Variational formulation and structural stability of monotone
  equations.
\newblock {\em Calc. Var. Part. Diff. Eqns.}, 47, 273--317, 2013.

\bibitem[Vis15]{Visi15?SSFE}
{\scshape A.~Visintin}.
\newblock Structural stability of flows via evolutionary {$\Gamma$}-convergence
  of weak-type.
\newblock {\em arxiv:1509:03819}, 2015.

\bibitem[Vis17a]{Visi17?EGCW}
{\scshape A.~Visintin}.
\newblock Evolutionary {$\Gamma$}-convergence of weak type.
\newblock {\em arXiv:1706.02172}, pages 1--8, 2017.

\bibitem[Vis17b]{Visi17?SCSP}
{\scshape A.~Visintin}.
\newblock Structural compactness and stability of pseudo-monotone flows.
\newblock {\em arXiv:1706.02176}, 2017.

\end{thebibliography}
\end{document}